\input amstex 

\pageno=1

\def\1{\hbox{\rm\rlap {1}\hskip .03in{\rm I}}} 
 \def\id{{\text {id}}} \def\tr{{\text {tr}}}
\def\card{{\text {card}}}

\hsize=12.5cm \vsize=19.0cm \parindent=0.5cm \parskip=0pt \baselineskip=12pt
\topskip=12pt \def\skipaline{\vskip 12pt plus 1pt} 

\def\H#1\par{\skipaline\noindent\bf 
#1\rm\par\nobreak\skipaline\nobreak\noindent}

\def\HH#1\par{\skipaline\noindent\bf
#1\rm\par\nobreak\skipaline\nobreak\noindent} 
\def\Hom{{\text {Hom}}}

\def\Spin {{\text {Spin}}}

\def\End {{\text {End}}}
 \def\Aut {{\text {Aut}}}
 \def\Rep {{\text {Rep}}}

\def\Tr {{\text {Tr}}}

\def\Tr {{\text {Tr}}}

\def\mod {{\text {mod}}}

\def\Dim {{\text {Dim}}}
 \def\dim  {{\text {dim}}}

\def\mod {{\text {mod}}\,}

\def\id {{\text {id}}} 
 
\def\pr {{\text {pr}}}

\def\Perm {{\text {Perm}}}

 \centerline {  \bf   Homotopy  field  theory in dimension 3 and  crossed group-categories}

\skipaline
\centerline {Vladimir  Turaev}

\skipaline
 \skipaline
 \skipaline
 \centerline {  \bf   Abstract}
\skipaline
 
A $3$-dimensional homotopy quantum field
theory (HQFT)   can be   described as a TQFT   for  surfaces and $3$-cobordisms endowed with homotopy
classes of maps into a given space.   For a group $\pi$,  we
 introduce a notion of a modular crossed $\pi$-category and show that
 such a category gives rise to a 3-dimensional HQFT with target space $K(\pi,1)$.  
This includes numerical invariants of 3-dimensional $\pi$-manifolds
and a  2-dimensional homotopy modular functor. We also introduce and discuss 
a parallel notion of a quasitriangular crossed Hopf $\pi$-coalgebra.

\skipaline

 \skipaline
\skipaline
\centerline {\bf Contents}
\skipaline

\noindent  {\bf Introduction}

\noindent {\bf 1. Group-categories} 

\noindent {\bf 2.  Crossed, braided, and ribbon $\pi$-categories}

\noindent  {\bf 3.  Colored $\pi$-tangles and their invariants}

\noindent {\bf 4.  Colored $\pi$-graphs and their invariants}

\noindent {\bf 5.  Trace,  dimension, and algebra of colors  }

\noindent {\bf 6.  Modular crossed $\pi$-categories }

\noindent  {\bf 7.  Invariants of 3-dimensional $\pi$-manifolds}

\noindent {\bf 8.   A 2-dimensional homotopy modular functor}

\noindent {\bf 9.   A 2-dimensional HQFT}

\noindent {\bf 10.   A 3-dimensional HQFT}

\noindent {\bf 11.   Hopf group-coalgebras}

\noindent {\bf 12.  Canonical extensions}

\noindent  {\bf 13.  Transfer of categories}

\noindent {\bf   Appendix 1. Quasi-abelian  cohomology of groups}

\noindent {\bf   Appendix 2. State sum invariants of 3-dimensional
$\pi$-manifolds}

\noindent  {\bf   Appendix 3. Open problems}

\skipaline
\skipaline
\skipaline
\centerline {\bf    Introduction}

\skipaline

A homotopy quantum field theory (HQFT)  
 is a version of a   topological quantum field theory (TQFT) 
for   manifolds endowed with  maps into  a fixed 
topological space.  The   HQFT's were introduced in [Tu3]  where   the author gave an
  algebraic characterization of   2-dimensional HQFT's
 whose target space is  the Eilenberg-MacLane space $K(\pi,1)$ determined by
a group $\pi$.
In this  paper we  focus on 3-dimensional HQFT's
with target space    $K(\pi,1)$.  A manifold $ M$ endowed with a homotopy class  of maps   
$M\to K(\pi,1)$ is called a {\it $\pi$-manifold}. The homotopy classes of maps    
$M\to K(\pi,1)$ classify principal $\pi$-bundles over $M$ and (for
connected $M$) bijectively correspond to the homomorphisms    $\pi_1(M)\to \pi$.   
A 3-dimensional HQFT with target space   $K(\pi,1)$
  comprises  two   ingredients:   a homotopy  modular functor assigning
$K$-modules to    $\pi$-surfaces and an invariant of 3-dimensional 
$\pi$-manifolds taking value in the module associated with the boundary. In
particular,  the HQFT  provides numerical invariants of closed 3-dimensional 
$\pi$-manifolds.   Our main  aim   is to
introduce an algebraic technique allowing to construct   3-dimensional HQFT's.

Our approach is based on a deep connection between the theory of braided
categories and   invariants of knots,  links   and 3-manifolds. 
This connection has been essential in the construction of \lq\lq quantum"
invariants of knots and 3-manifolds from quantum groups,  see [RT],  [Tu2], [KRT]. 
Here we extend these ideas to  links 
$\ell\subset  S^3$  endowed with homomorphisms  $\pi_1(S^3\backslash \ell) \to
\pi$ and to  3-dimensional  $\pi$-manifolds. To this end we introduce a notion
of a crossed $\pi$-category and study braidings and twists in such categories. This
leads us to the notion of a modular crossed $\pi$-category.  We show that 
each modular crossed $\pi$-category gives rise to a 3-dimensional HQFT
with target     $K(\pi,1)$.
In the case $\pi=1$ we recover the usual construction of 
  3-dimensional TQFT's from modular categories, see [Tu2]. 

The  crossed $\pi$-categories are quite delicate algebraic objects. 
We discuss a few general methods producing such categories.
In particular we introduce {\it quasitriangular   Hopf $\pi$-coalgebras}
and show that they give rise to
crossed $\pi$-categories. 
Other methods are based on 
a study  of self-equivalences of a braided category,  a study of quasi-abelian
cohomology of $\pi$, and    a transfer-type  construction.
 This gives  several   examples of modular crossed
$\pi$-categories.
However, the problem of  systematic finding  of 
modular crossed $\pi$-categories is largely open. It would be most
interesting to extend the   quantum groups   associated with
semisimple finite dimensional complex Lie algebras   to
quasitriangular   Hopf $\pi$-coalgebras.

The content of the paper is as follows. In Sections 1 and 2 we introduce
crossed $\pi$-categories and various additional structures on them (braiding,
twist etc.). In Sections 3 and 4 we introduce  $\pi$-links, $\pi$-tangles and
$\pi$-graphs in $\bold R^3$ colored over  a ribbon crossed 
$\pi$-category $\Cal C$. We also define their canonical   functorial invariant
taking values in $\Cal C$. In
Section 5 we   study traces of morphisms in $\Cal C$.
 In Section 6 we introduce modular crossed 
$\pi$-categories. They are used in Section 7 to define invariants of 3-dimensional
$\pi$-manifolds. In Sections 8-10 we  introduce the 2-dimensional
  and   3-dimensional HQFT's
derived from a modular crossed 
$\pi$-category. In Sections 11-13 we discuss     algebraic
  constructions of crossed $\pi$-categories.  

The paper is ended with three appendices. In Appendix 1 we
briefly discuss  quasi-abelian 3-cohomology of groups and their relations to 
crossed categories. In Appendix 2 we outline a state sum approach to invariants of
3-dimensional $\pi$-manifolds. In Appendix 3 we discuss a few open problems.

Throughout the paper, the symbol $K$ denotes a  commutative ring with unit.
The symbol $\pi$ denotes a group. By a category we   mean a small  category.

  We shall work  in the smooth set up although all our definitions have 
  topological and piecewise-linear versions. Thus, by   manifolds   
 we   mean   smooth
manifolds and by homeomorphisms we   mean smooth homeomorphisms.

\skipaline 
 \centerline  {\bf 1. Group-categories}

\skipaline \noindent {\bf 1.1.  Generalities on monoidal categories.}  Let ${\Cal C}$ 
be
 a  monoidal category with   unit object $\1$.    Recall (see for instance [Ma]) that
we have
 invertible associativity morphisms $$\{a_{U,V,W}:  (U\otimes V) \otimes W \to
 U\otimes (V\otimes W)\}_{U,V,W\in {\Cal C}} \leqno (1.1.a) $$ and invertible 
morphisms
 $$\{l_U:  U\to U\otimes \1, \,\, \,r_U:U\to \1 \otimes U\}_{U\in {\Cal C}} \leqno
(1.1.b)$$
 satisfying the pentagon identity $$(\id_U \otimes a_{V,W,X})\,a_{U,V\otimes
 W,X}\,(a_{U,V,W} \otimes \id_X) =a_{U,V,W\otimes X}\,a_{U\otimes V,W, X}\leqno
 (1.1.c)$$ and the triangle identities $$a_{U,\1,V} (l_U\otimes \id_V)=\id_U
 \otimes r_V \leqno (1.1.d)$$ for any $U,V,W,X\in {\Cal C}$, where the inclusion $U\in
 {\Cal C}$ means that $U$ is an object of ${\Cal C}$.  The morphisms $l,r$ 
should satisfy $l_{\1}=r_{\1}$ and be
 natural in the sense that for any morphism $f:U\to V$ we have $l_V f=(f\otimes
 \id_{\1}) l_U, r_V f= (\id_{\1}\otimes f) r_U$.  The associativity morphisms
 (1.1.a) should be natural in a similar sense.

A left duality in ${\Cal C}$ associates to any object $U\in {\Cal C}$ an object $U^*\in {\Cal C}$ and
two morphisms $$ b_U:  \1\to U\otimes U^*, \,\,d_U:U^*\otimes U\to \1 \leqno 
(1.1.e)
$$ such that $$ (l_U)^{-1} (\text {id}_U\otimes d_U) \,a_{U,U^*,U} \,(b_U 
\otimes
\text {id}_U) r_U=\text {id}_U, \leqno (1.1.f)$$ $$ (r_{U^*})^{-1}(d_U\otimes
\text {id}_{U^*}) (a_{U^*,U,U^*})^{-1} (\text {id}_{U^*}\otimes b_U)
l_{U^*}=\text {id}_{U^*}.  \leqno (1.1.g) $$

We   call the morphisms (1.1.a), (1.1.b), (1.1.e) the {\it structural morphisms}
of ${\Cal C}$.

A monoidal category ${\Cal C}$ is {\it strict} if $(U\otimes V) \otimes W = U\otimes
(V\otimes W)$, $U= U\otimes \1=\1 \otimes U$ for any $U,V,W\in {\Cal C}$ and the
morphisms $\{a_{U,V,W}\}_{U,V,W\in {\Cal C}} $ and $\{l_U, r_U\}_{U\in {\Cal C}} $ are the
identity morphisms.  It is well-known that each monoidal category is equivalent
to a strict monoidal category in a canonical way.

 \skipaline \noindent {\bf 1.2.   $\pi$-categories.}  
A monoidal category ${\Cal C}$ 
is said to be 
 {\it $K$-additive} if the $\Hom's$ in ${\Cal C}$ are $K$-modules and
both the
 composition and  the tensor product of morphisms are bilinear over $K$.
We say that a $K$-additive category ${\Cal C}$ {\it splits as a disjoint union of
subcategories} $\{{\Cal C}_\alpha\}$ numerated by certain $\alpha $ if:  

- each ${\Cal C}_\alpha$
is a full subcategory of ${\Cal C}$;

- each object of ${\Cal C}$ belongs to ${\Cal C}_\alpha$ for a unique
$\alpha $;

- if $U\in {\Cal C}_\alpha$ and $V\in {\Cal C}_\beta$ with $\alpha\neq \beta$ then
$\Hom_{\Cal C} (U,V)=0$.

For a group $\pi$,  a {\it  
$\pi$-category over $K$} is a $K$-additive monoidal category with left duality
${\Cal C}$ which splits as a disjoint union of subcategories $\{{\Cal C}_\alpha\} $ numerated
by $\alpha \in \pi$   such that

(i) if $U\in {\Cal C}_\alpha$ and $V\in {\Cal C}_\beta$ then $U\otimes V\in {\Cal C}_{\alpha\beta}$;

(ii) if $U\in {\Cal C}_\alpha$ then $U^*\in {\Cal C}_{\alpha^{-1}}$.

We shall write   ${\Cal C}=\amalg_\alpha {\Cal C}_\alpha$ and call the subcategories
$\{{\Cal C}_\alpha\} $ of ${\Cal C}$ the {\it components} of ${\Cal C}$.  The category ${\Cal C}_1$
corresponding to the neutral element $1\in \pi$ is called the 
{\it neutral component} of ${\Cal C}$. Conditions (i) and (ii) show that ${\Cal C}_1$ is
closed under tensor multiplication and taking the dual object. 
Condition (i) implies the inclusion $\1\in {\Cal C}_1$. Thus,
${\Cal C}_1$ is a $K$-additive  monoidal category with left duality.

 \skipaline \noindent {\bf 1.3.  Example:    $\pi$-categories from
3-cocycles.}  It is well-known that 3-cocycles give rise to associativity
morphisms in categories.  Here we elaborate this construction in our context.
Let $a=\{a_{\alpha,\beta,\gamma} \in K^*\}_{\alpha,\beta,\gamma\in \pi}$ be a
3-cocycle of the group $\pi$ with values in the multiplicative group $K^*$
consisting of invertible elements of $K$.   Thus
$$ a_{\alpha \beta,\gamma,\delta}\,a_{\alpha, \beta,\gamma \delta} = a_{\alpha,
\beta,\gamma} \, a_{\alpha, \beta \gamma, \delta} \, a_{ \beta,\gamma, \delta}
\leqno (1.3.a)$$ for any $\alpha,\beta, \gamma,\delta\in \pi$.  Let $b=
\{b_{\alpha} \in K^*\}_{\alpha \in \pi}$ be a family of elements of $K^*$
numerated by $\pi$.   With the pair $(a,b)$ we associate a  
$\pi$-category ${\Cal C} $ as follows.

For $\alpha \in \pi$, we define ${\Cal C}_{\alpha}$ to be a category with one object
$V_\alpha$.  For $\alpha, \beta \in \pi$, set $$ \Hom (V_\alpha,V_\beta)=\cases
K,~ { {if}}\,\,\, \alpha=\beta, \\ 0,~ { {if}}\,\,\, \alpha\neq \beta.\endcases $$ The
tensor  product
is given by $V_\alpha\otimes V_\beta=V_{\alpha\beta}$.  The composition and
tensor product of morphisms is given by multiplication in $K$.  Clearly, 
$\1=V_1$
is the unit object of ${\Cal C}$.  Now we define the structural morphisms  
in ${\Cal C}$.   The
associativity morphism (1.1.a) for $U=V_\alpha, V=V_\beta, W=V_\gamma$ is by
definition $ a_{\alpha, \beta,\gamma}\in K^* \subset K= \End
(V_{\alpha\beta\gamma})$.  The pentagon equation follows from (1.3.a).  The
morphisms $l_\alpha:  V_\alpha\to V_\alpha\otimes \1=V_\alpha$ and $
r_\alpha:V_\alpha\to \1 \otimes V_\alpha=V_\alpha$ are defined by
$$l_\alpha=(a_{\alpha,1,1})^{-1} \in K^*\subset K= \End (V_{\alpha }),\,\, \,\
r_\alpha= a_{ 1,1,\alpha} \in K^*\subset K= \End (V_{\alpha }).$$ The triangle
identity follows from the equality $a_{\alpha,1, \beta}=a_{\alpha,1, 1} 
\,a_{1,1,
\beta}$ obtained from (1.3.a) by the substitution $\beta=1, \gamma=1,
\delta=\beta$.  The dual of $V_\alpha$ is by definition $V_{\alpha^{-1}}$.  The
duality morphism $ \1\to V_\alpha\otimes V_{\alpha^{-1}}=\1$ is defined to be
$b_\alpha\in K^*\subset K= \End (\1)$.  The morphism $d_\alpha:V_{\alpha^{-1}}
\otimes V_\alpha\to \1$ is determined uniquely from   (1.1.f) and
(1.1.g).  In fact these equations give two different expressions for 
$d_\alpha$:
$$d_\alpha=(b_\alpha)^{-1}\, (a_{\alpha, \alpha^{-1},\alpha}\, a_{\alpha , 
1,1}\,
a_{1,1,\alpha })^{-1}, \leqno (1.3.b)$$  
 $$ d_\alpha=(b_\alpha)^{-1}\,
a_{\alpha^{-1}, \alpha,\alpha^{-1}}\, a_{\alpha^{-1}, 1,1}\,
a_{1,1,\alpha^{-1}}.\leqno (1.3.c)$$ To show that the right-hand sides are 
equal
we substitute $\beta=\delta=\alpha^{-1}, \gamma=\alpha$ in (1.3.a).  This gives
$$ a_{1, \alpha,\alpha^{-1}}\,a_{ \alpha,\alpha^{-1},1}= a_{\alpha,
\alpha^{-1},\alpha }\, a_{\alpha , 1,\alpha^{-1}} \,a_{\alpha^{-1},
\alpha,\alpha^{-1}}$$ $$=a_{\alpha, \alpha^{-1},\alpha }\, a_{\alpha , 1,1} \,
a_{1 , 1,\alpha^{-1}} \ a_{\alpha^{-1}, \alpha,\alpha^{-1}}.$$ Substituting
$\alpha=\gamma=1 $ in (1.3.a) we obtain that $a_{1,\alpha , 1 }=1$ for all
$\alpha$.  Substituting $\gamma=\delta=1, \beta= \alpha^{-1}$ (resp.  
$\alpha=1,
\delta=\beta=\gamma^{-1}$) in (1.3.a) we obtain that $a_{ 
\alpha,\alpha^{-1},1}=
(a_{\alpha^{-1}, 1, 1 })^{-1}$ (resp.  $a_{1, \alpha,\alpha^{-1}}= (a_{1,
1,\alpha })^{-1}$).  Combining these equalities we obtain that the right-hand
sides of (1.3.b) and (1.3.c) coincide.
It is clear that
${\Cal C}$ is a   $\pi$-category.

 \skipaline \noindent {\bf 1.4.  Operations on group-categories.}  We define
here a few elementary operations on group-categories.  Note first that the
group-categories can be pulled back along group homomorphisms.  Having a 
group
homomorphism $q:\pi'\to \pi$ we can transform any   $\pi$-category ${\Cal C}$ 
into
a  $\pi'$-category ${\Cal C}'=q^*({\Cal C})$ defined by ${\Cal C}'_{\alpha}= {\Cal C}_{q(\alpha)}$ 
for
any $\alpha\in \pi'$.  Composition, tensor multiplication, and the structural
morphisms
  in ${\Cal C}$ induce the corresponding operations  in ${\Cal C}'$ in the obvious way.  

The
group-categories can be pushed forward along group homomorphisms.  Having a 
group
homomorphism $q:\pi'\to \pi$ we can transform any   $\pi'$-category ${\Cal C}'$ 
into
a  $\pi$-category ${\Cal C}=q_*({\Cal C}')$ defined by ${\Cal C}_{\alpha}=
\amalg_{\beta\in q^{-1} (\alpha)} {\Cal C}'_{\beta}$  for
any $\alpha\in \pi$.  Composition, tensor multiplication, and the structural
morphisms
  in ${\Cal C}'$ induce the corresponding operations  in ${\Cal C}$ in the obvious
way.

For any family of  $\pi$-categories $\{{\Cal C}^i\}_{i\in I}$, we define a 
direct product ${\Cal C}=\prod_i {\Cal C}^i$.  The category ${\Cal C}$ is a disjoint union of
categories $\{{\Cal C}_\alpha\}_{\alpha\in \pi}$.  The objects of ${\Cal C}_\alpha$ are
families $\{U_i\in {\Cal C}^i_\alpha\}_{ i\in I}$.  The operations on the objects and
the unit object  are defined by $$\{U_i \}_{ i}\otimes \{U'_i \}_{
i}=\{U_i\otimes U'_i \}_{ i},  \,\,\, (\{U_i \}_{ i})^*=\{U^*_i 
\}_{
i},\,\,\, \1 =\{\1_{{\Cal C}^i} \}_{ i},  \leqno (1.4.a)$$ where   $i$ runs over $I$.  A
morphism $\{U_i\in {\Cal C}^i \}_{ i}\to \{U'_i\in  {\Cal C}^i
\}_{ i}$ in ${\Cal C} $ is a  family $\{f_i:U_i \to U'_i \}_{ i}$ where each $f_i$ is 
a
morphism in ${\Cal C}^i $.  The $K$-additive structure in ${\Cal C}$ is defined
coordinate-wisely so that $$\Hom_{{\Cal C}} (\{U_i \}_{ i},\{U'_i \}_{ i})= \prod_{
i\in I} \Hom_{{\Cal C}^i} (U_i,U'_i).$$ The composition of morphisms is 
coordinate-wise,
i.e., $\{f'_i:U'_i \to U''_i \} \circ \{f_i:U_i \to U'_i \}= \{f'_i f_i:U_i \to
U''_i \}$.  The tensor product for morphisms and  the structural morphisms $a,l,r,d,b$
are defined coordinate-wisely.  All the axioms of a   $\pi$-category follow
from the fact that they are satisfied coordinate-wisely.

For a finite family of   $\pi$-categories $\{{\Cal C}^i\}_{i\in I}$, we  
define
a tensor product ${\Cal C}'=\otimes_i {\Cal C}^i$.  The category ${\Cal C}'$ is a disjoint union 
of
categories $ \{{\Cal C}'_\alpha\}_{\alpha\in \pi}$.  The objects of ${\Cal C}'_\alpha$ are 
the
same as the objects of the category ${\Cal C}_\alpha\subset \prod_i {\Cal C}^i$ above. 
 The operations on
the objects and the unit object  are defined by   (1.4.a).  By definition, $$\Hom_{{\Cal C}'}
(\{U_i \}_{ i},\{U'_i \}_{ i})=  \bigotimes_{
i\in I} \Hom_{{\Cal C}^i} (U_i,U'_i).$$ This   $K$-module is additively generated by
elements of type $\otimes_{ i} (f_i:U_i \to U'_i)$.  The composition of 
morphisms
is defined on the generators by $$ \bigotimes_{ i} (f'_i:U'_i \to U''_i)\, 
\circ
\,\bigotimes_{ i} (f_i:U_i \to U'_i)= \bigotimes_{ i} (f'_if_i:U_i \to 
U''_i).$$
This extends to arbitrary morphisms by $K$-linearity and makes ${\Cal C}'$ a
$K$-additive category.  The tensor product of morphisms is defined on the
generators by $ (\otimes_{ i} f_i ) \otimes (\otimes_{ i} g_i ) =\otimes_{ i} (
f_i\otimes g_i)$ and extends to arbitrary morphisms by $K$-linearity.  Observe
that there is a canonical functor $  \prod_i {\Cal C}^i \to {\Cal C}'$ which is the
identity on the objects and sends a morphism $\{f_i \}_{ i}$ in $\prod_i {\Cal C}^i$
into the morphism $\bigotimes_{ i} f_i $ in ${\Cal C}'$.  This functor is by no means
$K$-linear but  does
preserve the tensor product.  Applying  this functor to the structural
morphisms in $\prod_i {\Cal C}^i$
 defined above we obtain structural morphisms in ${\Cal C}'$ 
satisfying
all the conditions of Section 1.1.     In this way ${\Cal C}'$ becomes a  
$\pi$-category.

\skipaline \centerline {\bf 2.  Crossed, braided, and ribbon $\pi$-categories}

 \skipaline \noindent {\bf 2.1.  Crossed $\pi$-categories.} 
Let ${\Cal C}$ 
be
 a $K$-additive monoidal category with left duality.  By
an {\it automorphism} of ${\Cal C}$ we   mean an invertible $K$-linear (on the
morphisms) functor $\varphi:{\Cal C}\to {\Cal C}$ which preserves the tensor product, the unit
object,  the
duality and the structural morphisms $a,l,r,b,d$.  Thus, $$\varphi(\1)=\1,\,\,
\varphi(U\otimes V)=\varphi(U)\otimes \varphi( V),\,\, \varphi(U^*)=
(\varphi(U))^*,\,\, \varphi(f\otimes g)=\varphi(f)\otimes \varphi( g)$$ for any
objects $U,V$ and any morphisms $f,g$ in ${\Cal C}$ and $$
\varphi 
(a_{U,V,W})=a_{\varphi
(U),\varphi (V),\varphi (W)},\,\, \varphi (l_{U})=l_{ \varphi
(U) },\,\, \varphi (r_{U})=r_{ \varphi (U) },$$ $$
 \varphi (b_U)=
b_{\varphi(U)}, \,\,\, \varphi (d_U)= d_{\varphi(U)} $$ for any objects $U,V,W\in {\Cal C}$. 
The group of automorphisms of ${\Cal C}$ is denoted by $\Aut ({\Cal C})$.

A {\it  crossed 
$\pi$-category over $K$} is a $\pi$-category 
${\Cal C}$   endowed with a group homomorphism $\varphi:\pi\to 
\Aut ({\Cal C})$ such that
  for all $ \alpha, \beta \in \pi$ the functor $
\varphi_{\alpha}=\varphi(\alpha):{\Cal C}\to {\Cal C}$  maps    ${\Cal C}_\beta $
into ${\Cal C}_{\alpha \beta  \alpha^{-1}}$.
  For objects $U\in {\Cal C}_\alpha, V\in {\Cal C}_\beta$,  set
$${}^U\!V=\varphi_{\alpha}  (V)\in
{\Cal C}_{\alpha\beta\alpha^{-1}}.$$ 
In particular,
${}^U\!U=\varphi_{\alpha} 
(U)\in
{\Cal C}_{\alpha}$ for any $U\in {\Cal C}_\alpha$. Note the identities $${}^{U }\!  (V\otimes 
W)={}^U\!
V\otimes {}^{U }\!  W \leqno (2.1.a)$$ $${}^{(U\otimes V)}\!  W= {}^{U }\!
({}^V\!  W), \leqno (2.1.b)$$ $${}^{U }\!  (V^*)=({}^{U }\!  V)^*, \leqno
(2.1.c)$$ $${}^{\1 }\!  V={}^{U }\!  ({}^{U^* }\!  V)={}^{U^* }\!  ({}^{U}\!
V)=V, \,\,\, {}^{U}\!  \1=\1, \leqno (2.1.d)$$ for any $U,V,W\in {\Cal C}$.  
Similarly,
for an object $U\in {\Cal C}_\alpha$ and a morphism $f:V\to V'$ in   ${\Cal C} $ 
set
$${}^U\!f=\varphi_{\alpha} (f):  {}^U\!V\to {}^U\!(V').$$ Note the identities
$${}^{U}\!  (f' \circ f)= {}^{U}\!  (f') \circ {}^{U}\!  f, \leqno (2.1.e)$$
$${}^{U }\!  (f\otimes g)={}^U\!  f\otimes {}^{U }\!  g, \leqno (2.1.f)$$ 
$${}^{U
}\!  (id_V)=id_{({}^{U }\!  V)}, {}^{U }\!  (b_V)=b_{({}^{U }\!  V)}, {}^{U }\!
(d_V)=d_{({}^{U }\!  V)}, \leqno (2.1.g) $$ $$ {}^{U }\!  (l_V)=l_{({}^{U }\!
V)}, {}^{U }\!  (r_V)=r_{({}^{U }\!  V)}, {}^{U }\!  (a_{V,W,X})=a_{({}^{U }\!
V), ({}^{U }\!  W), ({}^{U }\!  X)},$$ $${}^{(U\otimes V)}\!  f= {}^{U }\!
({}^V\!  f), \,\,\, {}^{\1 }\!  f = {}^{U }\!  ({}^{U^* }\!  f)={}^{U^* }\!
({}^{U}\!  f)=f.\leqno (2.1.h)$$

Examples of crossed $\pi$-categories will be given in Section 2.6 and in further
sections.

 \skipaline \noindent {\bf 2.2.  Braiding in   $\pi$-categories.}  Let ${\Cal C}$ be a
crossed $\pi$-category.  A {\it braiding} in ${\Cal C}$ is a system of invertible
morphisms $$\{c_{U,V}:U\otimes V \to {}^U\!  V \otimes U\}_{U,V\in {\Cal C}} \leqno
(2.2.a)$$ satisfying the following three conditions:

(2.2.1) for any morphisms $f:U\to U', g:V\to V'$ such that $U,U'$ lie in the 
same
component of ${\Cal C}$, we have $$c_{U',V'} (f\otimes g)= ({}^U\!  g \otimes f)
\,c_{U,V}; \leqno 
(2.2.b)$$

(2.2.2) for any objects $U,V,W\in {\Cal C}$ we have $$c_{U\otimes V,W} \leqno 
(2.2.c)$$
$$=a_{{}^{(UV)}\!  W,U,V} (c_{U,{}^V\!  W} \otimes \id_V) (a_{U,{}^V\!
W,V})^{-1} (\id_U \otimes c_{V,W})\, a_{U,V,W}, $$ $$ c_{U,V\otimes W}\leqno
(2.2.d)$$ $$= (a_{({}^U\!  V), ({}^U\!  W), U})^{-1} (\id_{({}^U\!  V)} \otimes
c_{U,W})\, a_{({}^U\!  V),U, W} (c_{U,V} \otimes \id_W) (a_{U,V,W})^{-1};$$

(2.2.3) the action of $\pi$ on ${\Cal C}$ preserves the braiding, i.e., for any
$\alpha\in \pi$ and any $V,W\in {\Cal C}$ we have $$\varphi_{\alpha}
(c_{V,W})=c_{\varphi_{\alpha} (V), \varphi_{\alpha}(W)}.$$

  Note that if in   (2.2.1) the objects $U,U'$ do not
lie in the same component of ${\Cal C}$ then both sides of (2.2.b)  are equal to 0 and have
the same source $U\otimes V$ but may have different targets.  Formulas (2.1.a)
and (2.1.b)  imply
that the targets of the morphisms on the left and right hand sides of (2.2.c) 
and
(2.2.d) are the same so that these equalities make sense.

A crossed $\pi$-category endowed with a braiding is said to be {\it braided}.
For $\pi=1$,  we obtain the standard 
definition of a braided monoidal category.

A braiding in a crossed $\pi$-category ${\Cal C}$ satisfies a version of the
Yang-Baxter identity. Assume for simplicity that
${\Cal C}$ is strict. Then for any braiding (2.2.a)
in ${\Cal C}$ and any objects $U,V,W \in {\Cal C}$,
$$(c_{({}^U\!  V),({}^U\!  W)} \otimes \id_U)\,
(\id_{({}^{U}\!  V)} \otimes c_{U,W})\,
(c_{U,V} \otimes \id_W)
 \leqno (2.2.e)$$
$$
=(\id_{({}^{UV}\!  W)} \otimes c_{U,V})\,   (c_{U,{}^V\!  W} \otimes \id_V)
\, (\id_{U} \otimes c_{V,W}).$$
Indeed,
by (2.2.d) and (2.2.3),
$$
(c_{({}^U\!  V),({}^U\!  W)} \otimes \id_U)\,
(\id_{({}^{U}\!  V)} \otimes c_{U,W})\,
(c_{U,V} \otimes \id_W)
= ({}^U\! (c_{V,W}) \otimes \id_U)\, c_{U,V \otimes W}.$$
Applying (2.2.b) to $f=\id_U, g= c_{V,W}$ and using (2.2.d),
we obtain
$$({}^U\! (c_{V,W}) \otimes \id_U)\, c_{U,V \otimes W}=c_{U,{}^V\!  W \otimes V}
(\id_{U} \otimes c_{V,W})  $$
$$
=(\id_{({}^{UV}\!  W)} \otimes c_{U,V})\,   (c_{U,{}^V\!  W} \otimes \id_V)
\, (\id_{U} \otimes c_{V,W}).$$

If ${\Cal C}$ is strict, then applying (2.2.c), (2.2.d) to  $U=V=\1$ and $V=W=\1$  and using
the invertibility of $ c_{U,\1}, c_{\1,U}$,  we obtain 
$$
c_{U,\1}=c_{\1,U}=\text {id}_U\leqno (2.2.f)
$$
for any object $U\in {\Cal C}$.  

 \skipaline \noindent {\bf 2.3.  Twist in   $\pi$-categories.}    A {\it
twist} in a braided (crossed) $\pi$-category ${\Cal C}$ is a family of invertible
morphisms $\{\theta_U:U\to {{}^U\!  U}\}_{U\in {\Cal C}}$ satisfying the following
conditions:

(2.3.1) for any morphism $f:U\to V$ with $U,V$ lying in the same component of 
${\Cal C}$
we have $\theta_V f=({}^U\!  f)\, \theta_U$;

(2.3.2) for any $U\in {\Cal C}$ we have $(\theta_U\otimes \id_{U^*})b_U 
=(\id_{({{}^U\!
U})} \otimes \theta_{({{}^U\!  U})^*}) b_{({{}^U\!  U})}$;

(2.3.3) for any objects $U,V\in {\Cal C}$, we have $$\theta_{U\otimes V} =c_{{}^{(UV)}\!
  V,{}^U\!  U}\, c_{({{}^U\!  U}),({{}^V\!  V})}\, (\theta_{U}\otimes
\theta_V); \leqno (2.3.a)$$

(2.3.4) the action of $\pi$ on ${\Cal C}$ preserves the twist, i.e., for any $\alpha\in
\pi$ and any $V\in {\Cal C}$, we have $\varphi_{\alpha}
(\theta_V)=\theta_{\varphi_{\alpha}(V)}$.

 As an exercise, the reader may check that the morphisms on both sides of the
equations (2.3.1) - (2.3.4) have the same source and target.
If ${\Cal C}$ is strict, then it follows from (2.2.f) and (2.3.3) that
$\theta_{\1}=\id_{\1}$.

A braided crossed $\pi$-category endowed with a twist is called a {\it
ribbon crossed $\pi$-category}. For $\pi=1$,  we obtain the standard 
definition of a ribbon monoidal category.

The neutral component ${\Cal C}_1$ of a ribbon crossed $\pi$-category ${\Cal C}$ is a
ribbon category  in the usual sense of the word.
Note also that every ribbon crossed $\pi$-category is equivalent to a strict 
ribbon
crossed $\pi$-category in a canonical way.

 \skipaline \noindent {\bf 2.4.  Dual morphisms.}   Condition (2.3.2) is better
understood when it is rewritten in terms of dual morphisms.
For a  morphism $f:U\to V$ in
 a   monoidal category with left duality,    the dual (or transpose)
morphism $f^*:V^*\to  U^*$  is defined  by  
$$
f^*=(r_{U^*})^{-1} (d_{V}\otimes 
\text {id}_{U^*}) \,(a_{V^*, V, U^*})^{-1} (\text {id}_{V^*}\otimes  (f\otimes  \text
{id}_{U^*})) \,(\text {id}_{V^*}\otimes  
b_U)\, l_{V^*}. $$
It follows   from (1.1.g) that
$(\text {id}_U)^*=\text {id}_{U^*}$. It is well-known  that  $(f g)^*=g^* f^*$ for   
composable morphisms $f,g$.  
Condition (2.3.2) can be shown to be  equivalent to  
$$(\theta_{U})^*= \theta_{{}^U\!  (U^*)}. \tag 2.4.a  $$ 

 \skipaline \noindent {\bf 2.5.  Operations on   ribbon group-categories.}
The operations on group-categories defined in Sections 1.4   can be adopted
to the setting of crossed (resp. braided, ribbon) group-categories. 
   Having a 
group
homomorphism $q:\pi'\to \pi$ we can pull back any  crossed $\pi$-category ${\Cal C}$ 
into
a  crossed $\pi'$-category ${\Cal C}'=q^*({\Cal C})$   as in Section 1.4 
with action of $\pi'$  
defined by $\varphi_{\alpha}=\varphi_{q(\alpha)}:  {\Cal C}'_{\beta}\to
{\Cal C}'_{\alpha \beta  \alpha^{-1}} $ where  $\alpha, \beta\in \pi'$.  
A braiding (resp. twist) in  ${\Cal C}$ induces  a braiding (resp. twist) in  ${\Cal C}'$    in  the 
obvious way.  In particular, if $\pi'\subset \pi$ is a subgroup of $\pi$, then any
crossed (resp. braided, ribbon) $\pi$-category 
${\Cal C}=\amalg_{\alpha\in \pi} {\Cal C}_\alpha$ induces a 
crossed (resp. braided, ribbon) $\pi'$-category 
$ \amalg_{\alpha\in \pi'} {\Cal C}_\alpha$.

A crossed (resp. braided, ribbon)
group-category can be pushed forward along   group epimorphisms whose kernel
acts trivially on the category. 
Consider a  group
epimorphism $q:\pi'\to \pi$ whose kernel   acts as the identity
on a  crossed (resp. braided, ribbon) 
$\pi'$-category ${\Cal C}'$.  Then the action of $\pi'$ on  ${\Cal C}'$ induces an
action of $\pi$ on the push-forward  $\pi$-category $ q_*({\Cal C}')$
 defined
in Section 1.4.  A braiding (resp. twist) in  ${\Cal C}'$ induces  a braiding (resp.
twist) in  $q_*({\Cal C}')$    in  the  obvious way.

Given a  family of  crossed $\pi$-categories $\{{\Cal C}^i\}_{i\in I}$,   the
direct product $\prod_i {\Cal C}^i$ is a crossed $\pi$-category. The action of
$\alpha\in \pi$
 on objects  and morphisms is defined 
by
 $$   \varphi_\alpha(\{U_i \}_{
i\in I})=\{\varphi_\alpha ( U_i) \}_{ i\in I},\,\,\, 
 \varphi_\alpha(\{f_i \}_{
i\in I})=\{\varphi_\alpha ( f_i) \}_{ i\in I}.
 \leqno (2.5.a)$$  
If $\{{\Cal C}^i\}_{i\in I}$ are 
braided (resp.   ribbon)
$\pi$-categories then $\prod_i {\Cal C}^i$ is a braided (resp.  ribbon) $\pi$-category: the
braiding and twist are defined coordinate-wisely and their coordinates are
the braiding and twist  in $\{{\Cal C}^i\}_{i\in I}$, respectively.

For a finite family of  crossed  $\pi$-categories $\{{\Cal C}^i\}_{i\in I}$, 
the tensor product $ \otimes_i {\Cal C}^i$ is a crossed  $\pi$-category. The action
of $\alpha\in \pi$  on objects is defined 
as in   (2.5.a).
  The action of $\alpha\in \pi$ on morphisms 
is
defined on the generators by $\varphi_\alpha(\bigotimes_{ i} f_i )= 
\bigotimes_{i}
\varphi_\alpha (f_i )$ and extends to arbitrary morphisms   by
$K$-linearity.   
If $\{{\Cal C}^i\}_{i\in I}$ are 
braided (resp.   ribbon)
$\pi$-categories then $\otimes_i  {\Cal C}^i$ is a braided (resp.  ribbon) $\pi$-category:
  the braiding and twist   are obtained from the corresponding morphisms in 
$\prod_i {\Cal C}^i$
via the canonical functor $\prod_i {\Cal C}^i\to \otimes_i {\Cal C}^i $.

We define a  transformation  of crossed  $\pi$-categories  called
{\it reflection}.   Let  ${\Cal C}=\amalg_{\alpha\in \pi} {\Cal C}_\alpha$
be  a crossed  $\pi$-category with tensor product $\otimes$, duality ${\ast}$,
structural morphisms $a,l,r,b,d$ and $\pi$-action $\varphi:\pi \to \Aut ({\Cal C})$.
We define a   crossed $\pi$-category ${\overline {\Cal C}}=\amalg_{\alpha\in \pi} {\overline {\Cal C}}_\alpha$
with tensor product ${\overline \otimes}$, duality ${\star}$, structural
morphisms ${\overline a},{\overline l},{\overline r},{\overline b},{\overline d}$
 and $\pi$-action ${\overline \varphi}:\pi \to \Aut ({\overline {\Cal C}})$
as follows:

- ${\overline {\Cal C}}={\Cal C}$ as categories (but not monoidal categories);

-   ${\overline {\Cal C}}_\alpha={\Cal C}_{\alpha^{-1}}$ as categories for all $\alpha\in \pi$;

- $\1_{\overline {\Cal C}}=\1_{\Cal C}$;

- for objects $U\in {\overline {\Cal C}}_\alpha, V\in {\overline {\Cal C}}_\beta$, set 
 $U{\overline \otimes}
V=\varphi_{\beta^{-1}} (U)\otimes V\in {\overline {\Cal C}}_{\alpha\beta}$;

- for morphisms $f:U\to U', g:V\to V'$ in ${\overline {\Cal C}}$ 
with $U\in {\overline {\Cal C}}_\alpha, U'\in
{\overline {\Cal C}}_{\alpha'}, V\in {\overline {\Cal C}}_\beta, 
V' \in {\overline {\Cal C}}_{\beta'}$, set 
 $$f{\overline \otimes} g=\cases
\varphi_{\beta^{-1}} (f)\otimes g \in \Hom_{\overline {\Cal C}}(U{\overline \otimes} V, U'{\overline \otimes} V'),~ {
{if}}\,\,\,  \beta=\beta', 
\\ 0  \in \Hom_{\overline {\Cal C}}(U{\overline \otimes} V, U'{\overline \otimes} V'),~ { {if}}\,\,\,
\beta\neq \beta';\endcases $$ 

- for $U\in {\overline {\Cal C}}_\alpha$, set 
$U^{\star}=\varphi_\alpha (U^*)\in {\overline {\Cal C}}_{\alpha^{-1}}$
and $${\overline l}_U=l_U,  \,\,\,{\overline r}_U=r_U, \,\,\, {\overline d}_U=d_U,
\,\,\,
{\overline b}_U=\varphi_\alpha(b_U);$$

-  for objects $U\in {\overline {\Cal C}}_\alpha, V\in {\overline {\Cal C}}_\beta, W\in {\overline {\Cal C}}_\gamma$, set
$${\overline a}_{U,V,W}=a_{\varphi_{\gamma^{-1} \beta^{-1}} (U), 
\varphi_{\gamma^{-1}  } (V), W};$$

- for $\alpha\in \pi$ set ${\overline \varphi}_\alpha= \varphi_\alpha$.

A routine check   shows that ${\overline {\Cal C}}$ is a crossed $\pi$-category.
Moreover, if $c, \theta$ are a braiding and a twist in ${\Cal C}$ then the formulas
$${\overline c}_{U,V}=(c_{V,U})^{-1},\,\,\,\, {\overline \theta}_{U}=
(\theta_{\varphi_\alpha(U)})^{-1}  $$
(where $U\in {\overline {\Cal C}}_\alpha$) define a braiding and 
a twist in ${\overline {\Cal C}}$. We call ${\overline {\Cal C}}$ the  
{\it mirror} of ${\Cal C}$.  Its neutral component 
${\overline {\Cal C}}_1$ is the mirror of
${{\Cal C}}_1$ in the sense of [Tu2, Section I.1.4].  It is easy to see that  
$\overline {\overline {\Cal C}}={\Cal C}$.

 \skipaline \noindent {\bf 2.6.  Example.}  
Consider the $\pi$-category ${\Cal C}$ defined in Section  1.3 and 
assume that   both $a$ and $b$ are invariant under
conjugation, i.e., $$
a_{\delta\alpha\delta^{-1}, \delta \beta \delta^{-1}, \delta\gamma\delta^{-1}}=
a_{\alpha, \beta,\gamma }, \leqno (2.6.a)$$ 
and  $b_{\delta\alpha\delta^{-1} }= b_{\alpha }$ for any $\alpha,
\beta,\gamma,\delta\in \pi$. Then  ${\Cal C}$ is  a crossed $\pi$-category as
follows. For $\alpha,\beta\in \pi$,  set $\varphi_\alpha(V_\beta)= V_{\alpha
\beta  \alpha^{-1}}$.  This extends to  morphisms in ${\Cal C}$ in the obvious way since
all non-zero morphisms in ${\Cal C}$ are proportional to the identity endomorphisms of
objects.  The resulting functor $\varphi_{\alpha}:{\Cal C}\to {\Cal C}$ preserves all
the structural morphisms  in ${\Cal C}$ since $a$ and $b$ are  conjugation
invariant.   To construct specific examples we can take $b=1$. Finding conjugation
invariant 3-cocylces is a   delicate task.  Obvious examples:  the trivial
cocycle $a=1 $; any 3-cocycle in the case of abelian $\pi$.

It is clear that a braiding   in ${\Cal C}$ is given by a family $\{c_{\alpha,\beta}\in
K^*\}_{\alpha,\beta  \in
\pi}$  where 
$c_{\alpha,\beta}$
determines the braiding morphism $$ V_{\alpha \beta}= V_{\alpha }\otimes V_{
\beta} \to \varphi_\alpha(V_\beta) \otimes V_\alpha=V_{\alpha \beta \alpha^{-1}}
\otimes V_{\alpha }=V_{\alpha \beta}. \leqno (2.6.b)$$
  The conditions on the braiding  
can  be reformulated as    the following idenitites:  $$c_{\delta\alpha\delta^{-1},
\delta \beta \delta^{-1}}= c_{\alpha, \beta }, \leqno (2.6.c)$$ $$c_{\alpha \beta,
\gamma}= c_{ \beta, \gamma} \, c_{\alpha, \beta \gamma\beta^{-1}} \, a_{\alpha,
\beta, \gamma}\, (a_{\alpha, \beta \gamma\beta^{-1}, \beta})^{-1}
\,a_{\alpha\beta \gamma\beta^{-1}\alpha^{-1}, \alpha,\beta}, \leqno (2.6.d)$$
$$c_{\alpha, \beta \gamma}= c_{\alpha, \beta } \,c_{\alpha, \gamma}\, (a_{\alpha,
\beta, \gamma})^{-1} \,a_{\alpha\beta \alpha^{-1}, \alpha,\gamma}\, (a_{ \beta , 
\gamma
, \alpha})^{-1},\leqno (2.6.e)$$ 
for all $\alpha, \beta,
\gamma,\delta\in \pi$.
 The equality (2.6.d) can be rewritten in a more convenient form using
(2.6.c).  Namely, observe that $c_{\alpha, \beta \gamma\beta^{-1}}= c_{\delta ,
\gamma}$ with $\delta=\beta^{-1}\alpha\beta$. Now, $\alpha\beta =\beta \delta$
which gives the following equivalent form of (2.6.d):  $$c_{\beta \delta ,
\gamma}=c_{\beta, \gamma}\, c_{\delta , \gamma}\, a_{ \beta \delta\beta^{-1},
\beta,\gamma } \,( a_{ \delta, \gamma,\beta})^{-1}  \,
a_{ \delta\gamma\delta^{-1}, \delta, \beta }.  \leqno (2.6.f)$$
 A direct
computation   shows that (2.6.c)  follows
from (1.3.a), (2.6.a,e,f).

The definition of a twist in ${\Cal C}$ considerably simplifies since ${{}^U\!  U}=U$ for any  
$U\in {\Cal C}$. Given a braiding $\{c_{\alpha,\beta}\in
K^*\}_{\alpha,\beta  \in
\pi}$ in ${\Cal C}$, a twist in ${\Cal C}$ is determined by a family 
 $\{\theta_{\alpha }\in K^*\}_{\alpha \in \pi}$ (where $\theta_{\alpha }$
is the twist $V_{\alpha }\to V_{\alpha }=\varphi_\alpha ( V_{\alpha })$)
such that $$\theta_{\alpha\beta}=
c_{\alpha, \beta } c_{ \beta , \alpha } \theta_\alpha \theta_\beta,\leqno (2.6.g)$$
$$ \theta_{{\alpha}^{-1}}=\theta_{\alpha}\leqno (2.6.h)$$ for all $\alpha, \beta \in
\pi$.   The equality (2.6.g) implies that  $\theta_{\alpha\beta}=\theta_{ \beta
\alpha }$ so that $\theta$ is conjugation invariant.
 
To sum up, a conjugation invariant tuple $(a,b,c,\theta)$ satisfying 
(1.3.a), (2.6.e-h)  gives rise to a ribbon crossed $\pi$-category ${\Cal C}={\Cal C}
(a,b,c,\theta)$. 
Such tuples $(a,b,c,\theta)$ 
  form a  group
under  pointwise multiplication.  This group operation corresponds to 
tensor multiplication of the $\pi$-categories ${\Cal C}
(a,b,c,\theta)$.
The   reflection of ribbon $\pi$-categories defined in Section 2.5 
corresponds to the following involution   in the set of tuples $(a,b,c,\theta)$:
$${\overline a}_{\alpha, \beta,\gamma }=a_{\beta^{-1}\alpha^{-1} \beta,
\beta^{-1},\gamma^{-1}},\,\,\, {\overline c}_{\alpha, \beta
}=(c_{\beta^{-1},\alpha^{-1} })^{-1}, \leqno (2.6.i)$$ $${\overline b}_{\alpha  }=b_{
\alpha^{-1}},\,\,\,  {\overline \theta}_{\alpha  }=(\theta_{ \alpha  })^{-1} $$ where
$\alpha, \beta, \gamma \in \pi$.  We shall further discuss   equations 
(2.6.e-h) in Appendix 1. 

If $G\subset \pi$ is a subgroup of the center of $\pi$ then the action of $G$ on
${\Cal C} (a,b,c,\theta)$ is trivial so that we can push ${\Cal C}
(a,b,c,\theta)$ forward along the projection
$\pi\to \pi/G$. This gives   a ribbon crossed
$(\pi/G)$-category.

\skipaline \noindent {\bf 2.7.  Remarks.}  1. The   objects of a crossed
$\pi$-category
$(\Cal C, \varphi:\pi\to 
\Aut ({\Cal C}))$  form a $\pi$-automorphic set in terminology of Brieskorn [Br] or a
$\pi$-rack in terminology of Fenn and Rourke [FR2]. Recall that a $\pi$-rack
 is a set $X$
equipped with a left  action of $\pi$ and a map $\partial:X\to \pi$ such that
$\partial(\alpha a)= \alpha \partial (a) \alpha^{-1}$ for all $\alpha\in \pi, a\in X$.
The  underlying $\pi$-rack of $\Cal C$  comprises the set of objects of $\Cal C$,
the action of $\pi$ on this set induced by $\varphi$ and the map assigning to any
object of ${\Cal C}_\alpha$ the element $\alpha\in \pi$.

2. It is clear from definitions that a braiding  in a crossed $\pi$-category in general is not a braiding   
in the underlying monoidal category in the usual sense of the word. There is  one exceptional case. Namely, assume that $\pi$ is
abelian and $\Cal C$ is a $\pi$-category  as in Section 1.2. The trivial homomorphism $\varphi=1:\pi \to \Aut (\Cal
C)$ makes
$\Cal C$  a crossed $\pi$-category. It is clear that a braiding (resp.\ twist) in this crossed $\pi$-category is  a braiding
(resp.\ twist) in $\Cal C$ in the usual sense.
 
\skipaline \centerline {\bf 3.  Colored $\pi$-tangles and their invariants}

\skipaline \noindent {\bf 3.1.  Colored $\pi$-links.}  It is well-known that a
framed oriented link in $S^3$ whose components are colored with objects of a
ribbon category gives rise to an invariant lying in the ground ring of this
category, see [Tu2].  To adapt this theory to our present setting we 
introduce
$\pi$-links and their colorings.  For the sake of future references, we  
consider links in an arbitrary connected oriented 3-manifold $M$.

  Let $\ell=\ell_1\cup ...  \cup\ell_n \subset M$ be an oriented $n$-component
link in $M$ with $n\geq 0$.  Denote the open 3-manifold $M\backslash \ell$ by
$C_{\ell}$ where $C$ stands for complement.  We say that $\ell$ is {\it framed} if
each its component $\ell_i$ is provided with a  {\it longitude} $\tilde {\ell}_i
\subset C_{\ell}$ which goes very closely along $\ell_i$
 (it may wind around
$\ell_i$ several times).  Set $\tilde \ell= \cup_{i=1}^n
\tilde {\ell}_i$.  For a path $\gamma:[0,1] \to C_{\ell}$ connecting a point
$z=\gamma(0) $ to a point   $\gamma(1)\in \tilde {\ell}_i$, denote by
$\mu_{\gamma}\in \pi_1(C_{\ell},z) $ the (homotopy) meridian of ${\ell}_i$
represented by the loop $\gamma m_i \gamma^{-1}$, where $m_i $ is a small
loop encircling $\ell_i$ with linking number $+1$.  
We
similarly define a (homotopy) longitude  $\lambda_{\gamma}=[\gamma \tilde
{\ell}_i \gamma^{-1}]\in \pi_1(C_{\ell},z)$ where the square brackets denote 
the
homotopy class of a loop and the circle $\tilde {\ell}_i$ is viewed as a loop
beginning and ending in $\gamma(1)$.  Both $\mu_{\gamma}$ and $
\lambda_{\gamma}$ are invariant under homotopies of $\gamma$ fixing
$\gamma(0)=z$ and keeping $\gamma(1)$ on $\tilde {\ell}$.  Clearly,
$\mu_{\gamma}$ and $ \lambda_{\gamma}  $ commute in $\pi_1(C_{\ell},z) $.
If $\beta$ is a  loop in $(C_{\ell},z)$ (i.e., a path
$[0,1] \to C_{\ell}$ beginning and ending in $z$) then
$$\mu_{\beta \gamma}=[\beta] \,\mu_{  \gamma} \, [\beta]^{-1}\,\,\,\, {\text {and}}
\,\,\,\, \lambda_{\beta \gamma}=[\beta]\, \lambda_{  \gamma} \, [\beta]^{-1}.$$

By a {\it $\pi$-link} in $M$ we shall mean a triple (a framed oriented link
$\ell\subset M$, a base point $z \in C_{\ell} $, a group homomorphism $
g:\pi_1(C_{\ell},z) \to \pi$).

Fix a crossed $\pi$-category $({\Cal C},\varphi:\pi \to \Aut ({\Cal C}))$ which we   call
the {\it  category of colors}.  A {\it coloring}
 of a $\pi$-link $(\ell ,z , g )$ is a function  $u$  which
assigns to every path $\gamma:[0,1] \to C_{\ell}$ with  $ \gamma(0)=z,
\gamma(1)\in \tilde {\ell}$ an object $u_\gamma\in {\Cal C}_{g(\mu_{\gamma})}$ such
that

(i)  $u_\gamma$ is preserved under homotopies of $\gamma$ fixing
$\gamma(0)=z$ and keeping $\gamma(1)$ on $\tilde {\ell}$;

(ii) if $\beta $ is a  loop in $(C_{\ell},z)$, then $u_{\beta
\gamma}=\varphi_{g([\beta])}(u_\gamma)$.

Pushing the endpoint $\gamma(1)\in \tilde \ell$ of a path $\gamma$ as above along the
corresponding component of $\tilde \ell$ we can deform $\gamma$ into a path
homotopic to $ \lambda_{\gamma} \gamma$.  Conditions (i), (ii) imply that $$
\varphi_{g(\lambda_{\gamma})} (u_\gamma)=u_{\gamma
}.  \leqno (3.1.a)$$
We shall see below (Lemma 3.2.1) that a coloring of an
$n$-component $\pi$-link $(\ell=\ell_1\cup  ...
\cup\ell_n,z , g )$ is uniquely determined by the objects associated to
any
given $n$ paths  $\gamma_1,...,\gamma_n$
connecting $z $ to $\tilde {\ell}_1,...,\tilde {\ell}_n$, respectively. In the role of
these objects we can take any objects of   ${\Cal C}_{g(\mu_{\gamma_1})},...,{\Cal C}_{
g(\mu_{\gamma_n})}$
satisfying (3.1.a).

A $\pi$-link endowed with a coloring is said to be {\it
${\Cal C}$-colored} or briefly {\it
colored}. The notion of an ambient isotopy in $M$ applies to $\pi$-links and colored
$\pi$-links in $M$ in the obvious way. This allows us to  consider the (ambient)
isotopy classes of such links.

The structure
of a colored
$\pi$-link  $(\ell ,z , g,u )$ in $M$  can be transferred  along 
paths   in $M\backslash \ell$ relating various base points. 
  More precisely,  let
  $\rho:[0,1] \to M\backslash \ell$ be a path  with $\rho (0)=z$.
We define a new 
colored $\pi$-link $(\ell' =\ell, z' , g':\pi_1(C_{\ell},z')\to \pi,u' )$
by  $z'=\rho(1)$, 
$g' ([\alpha])= g( [\rho \alpha \rho^{-1}])$ for any loop $\alpha$ in $(M\backslash
\ell,z')$,
$u'_\gamma=u_{\rho \gamma} $ for any path 
$\gamma$ in $C_{\ell}$ leading from $z'$ to  $\tilde \ell $. It is clear that the
transfers along homotopic paths (with the same endpoints) give the same results.
The  transfer preserves the ambient isotopy class of the colored $\pi$-link
$(\ell ,z , g,u )$: it is ambiently isotopic to 
  $(\ell,z' , g' ,u' )$
via an isotopy of the identity map $\id_M:M\to M$ which pushes $z$ along
$\rho$ and is constant in a neighborhood of $\ell$.

Although we shall not need it, note that a ${\Cal C}$-colored $\pi$-link can be
defined  in terms of $\pi$-racks (cf. Section 2.7) as a
framed oriented link endowed with a homomorphism of its fundamental rack
(as defined in [FR2]) into the 
underlying $\pi$-rack of $\Cal C$.  
The notion of colored $\pi$-links can be  reformulated also in terms of principal
$\pi$-bundles over link complements.

\skipaline \noindent {\bf 3.2.  Colored $\pi$-tangles.}  By a tangle with 
$k\geq
0$ inputs and $ l\geq 0$ outputs we   mean a  tangle $T \subset
{\bold R}^2\times [0,1]$ with bottom endpoints (inputs) $(r,0,0), r=1,...,k$ 
and
top endpoints (outputs) $(s,0,1), s=1,...,l$.  The tangle $T$ consists of a
finite number of mutually disjoint oriented embedded circles and arcs lying in
the open strip $ {\bold R}^2\times ]0,1[$ except the endpoints of the arcs.  At the
endpoints, $T$ should be orthogonal to the planes ${\bold R}^2\times 0, {\bold
R}^2\times1$.  We   denote the open 3-manifold $({\bold R}^2\times [0,1] )
\backslash T$ by $C_{T}$.

We say that $T$ is {\it framed} if each its component $t$ is provided with a 
longitude
$\tilde {t}\subset C_{T}$ which goes very closely along $t$.  
Clearly, $\tilde {t}$ is an arc (resp.  a circle) if  $t$ 
is
an arc (resp.  a circle).  We   always  assume that the   longitudes
of the arc components of $T$ have the endpoints $(r,-\delta,0),
r=1,...,k$ and $(s,-\delta,1), s=1,...,l$ with small $\delta>0$.  Set $\tilde 
T=
\cup_{t} \tilde {t} $ where $t$ runs over all the components of $T$.

 As the base point of $C_{T}$ we   always choose a point $z$ with a big
negative second coordinate $z_2 \ll  0$ so that $T \subset {\bold R}\times
[z_2+1,\infty] \times [0,1]$.  The set of such $z$ is contractible.  This 
allows
us to supress the base point from the notation for the fundamental group of $
C_{T}$.

If $T$ is oriented, then for each path $\gamma:[0,1] \to C_{T}$ connecting the base
point to $\tilde T$,  we
introduce a meridian $\mu_{\gamma}\in \pi_1(C_{T}) $ as in Section 3.1.  If
$\gamma(1)$ lies on a circle component of $\tilde T$ then we also have a
longitude $\lambda_{\gamma}\in \pi_1(C_{T})$.  Both $\mu_{\gamma}$ and $
\lambda_{\gamma}$ are invariant under homotopies of $\gamma$ fixing $\gamma(0)$
and keeping $\gamma(1)$ on $\tilde {T}$.  Clearly, $\mu_{\gamma}$ and $
\lambda_{\gamma} $ commute in $\pi_1(C_{T}) $.

A {\it $\pi$-tangle} is a pair (a framed oriented tangle $T$, a group 
homomorphism
$ g:\pi_1(C_{T}) \to \pi$).  The definition of a coloring of 
a $\pi$-link  extends to $\pi$-tangles word for word.  
A $\pi$-tangle endowed with a coloring is said to be {\it
colored}.

As usual, we shall consider
$\pi$-tangles and colored $\pi$-tangles up to ambient isotopy in
${\bold R}^2\times [0,1]$ constant on the
endpoints.  The  $\pi$-tangles (resp.  colored $\pi$-tangles) with 
0
inputs and 0 outputs are nothing but $\pi$-links (resp.  colored $\pi$-links) 
in
${\bold R}^2\times ]0,1[$.

Let $ T=(T,g,u)$ be a colored $\pi$-tangle.  With the $r$-th input $(r,0,0)$ of $T$ 
we associate a triple $(\varepsilon_r=\pm, \alpha_r\in \pi, U_r\in {\Cal
C}_{\alpha_r})$ as follows.  Set $\varepsilon_r=+$ if the arc component of $T$
incident to the $r$-th input is directed  out of  $ {\bold R}^2\times [0,1]$  and set
$\varepsilon_r=-$ otherwise.  Let $\gamma_r$ be the  path in $C_T$ leading from
the base point $z=(z_1,z_2,z_3)$ of $C_T$ to the $r$-th input of $\tilde T$ and
defined as   composition of the linear
path from   $z $ to $(r,z_2,0)$ with the linear path
from $(r,z_2,0)$ to $(r,-\delta,0)$.   We call  $\gamma_r$ the {\it canonical path}
associated to  the $r$-th input of $ T$. The meridian
$\mu_r=\mu_{\gamma_r}\in \pi_1(C_T)$ is called the {\it canonical meridian} 
associated to the $r$-th input of $T$. 
Set $\alpha_r=g(\mu_r) \in \pi$  and  $U_r=u_{\gamma_r} \in {\Cal C}_{\alpha_r}$. 
The sequence $(\varepsilon_1,\alpha_1,U_1),...,(\varepsilon_k,\alpha_k,U_k)$
 (where $k$ is the number of inputs of $T$) is called the {\it source} of $T$. 
Similarly, with the $s$-th
output $(s,0,1)$ of $T$  we associate a triple
$(\varepsilon^s=\pm, \alpha^s\in \pi, U^s\in {\Cal C}_{\alpha^s})$ as follows.  Set $\varepsilon^s=+$
if the arc component of $T$ incident to the $s$-th output is directed  
inside   $ {\bold R}^2\times [0,1]$  and set $\varepsilon^s=-$ otherwise. 
Let $\gamma^s$ be the canonical path  leading from the base
point $z=(z_1,z_2,z_3)$ of $C_T$ to the $s$-th output $(s,-\delta,1)$ of $\tilde T$:
It is defined as the composition of the linear path  leading from $z$  to $(s,z_2,1)$
with the linear path  leading from
$(s,z_2,1)$ to $(s,-\delta,1)$.  Set  $\alpha^s=g(\mu_{\gamma^s})\in
\pi $ and $U^s=u_{\gamma^s} \in {\Cal C}_{\alpha^s}$.
 The sequence $(\varepsilon^1,\alpha^1,U^1),...,(\varepsilon^l,\alpha^l,U^l)$
 (where $l$ is the number of outputs of $T$) is called the {\it target} of $T$. 
 Note the obvious equality
$$\prod_{r=1}^k (\alpha_r)^{\varepsilon_r}= \prod_{s=1}^l
(\alpha^s)^{\varepsilon^s}.$$

A  colored $\pi$-tangle $(T,g:\pi_1(C_T)\to \pi, u)$ is said to be a {\it colored $\pi$-braid} if   $T$
is a  framed oriented  braid. Let $k$  be  the number of  strings of $T$.  Observe
that the orientation of $T$, the homomorphism $g$ and the coloring
$u$ are uniquely determined by the source 
$(\varepsilon_1,\alpha_1,U_1),...,(\varepsilon_k,\alpha_k,U_k)$ of $T$.  Indeed, the
signs $\varepsilon_1,...,\varepsilon_k$ determine the orientation of the strings.
The group $\pi_1(C_T)$ is free on   $k$ generators represented by the canonical
meridians  $\mu_1,...,\mu_k$ corresponding to the inputs  of
$T$.  Hence,  the homomorphism $g$ is 
  determined by   $\alpha_1=g(\mu_1), ...,
\alpha_k=g(\mu_k)$.  
It follows from definitions that the coloring $u$
is determined by $U_1,...,U_k$. The next lemma implies that, conversely,  given a
framed   braid $T$ on $k\geq 0$ strings and a finite sequence 
$(\varepsilon_1,\alpha_1,U_1),...,(\varepsilon_k,\alpha_k,U_k)$  with
$ \varepsilon_r=\pm, \alpha_r\in \pi, U_r\in {\Cal C}_{\alpha_r}$ for  
$r=1,...,k$, we can extend $T$ uniquely to a colored $\pi$-braid 
with source $(\varepsilon_1,\alpha_1,U_1),...,(\varepsilon_k,\alpha_k,U_k)$.

 \skipaline \noindent {\bf 3.2.1.  Lemma.  } {\sl Let 
  $(T=t_1\cup ... \cup t_n,  g )$ be an $n$-component $\pi$-tangle. Let
$\rho_i:[0,1] \to C_{T} $ be a path connecting the base point of $C_T$  to
a point of $\tilde {t}_i$ where $i=1,...,n$.  Let $U_{ i }\in {\Cal C}_{g(\mu_{\rho_i})}$ for
$i=1,...,n$. Assume that for each circle component $t_i$ of $T$ we have 
$\varphi_{g(\lambda_{\rho_i} )} (U_i)=U_{i}$. Then there is a unique   coloring $u$
of $T$ such that   $u_{\rho_i}=U_i$ for  $i=1,...,n$. }

\skipaline {\sl Proof.} 
Any path $\gamma:[0,1] \to  C_{T}$
connecting the base point $z $ to $\tilde {t}_i$ can be deformed fixing
$\gamma(0)=z$ and keeping $ \gamma(1) \in \tilde {t}_i$ so that $\gamma(1)=
\rho_i(1)$.  Then 
 $$u_{\gamma}=u_{\gamma \rho_i^{-1} \rho_i}
=\varphi_{g([ \gamma \rho_i^{-1}])}(u_{\rho_i}) 
=\varphi_{g([ \gamma \rho_i^{-1}])}(U_{i}). $$
  This proves the
uniqueness. Conversely, our assumptions imply that the object  $$u_{\gamma}=
\varphi_{g([ \gamma \rho_i^{-1}])}(U_{i}) \in  {\Cal C}_{g(\mu_\gamma)} \leqno
(3.2.a)$$ does not depend on the choice of the deformation of $\gamma$ used 
to ensure $\gamma(1)= \rho_i(1)$. It is  easy to check that 
formula   (3.2.a) defines a 
coloring $u$ of $\ell$ such that   $u_{\rho_i}=U_i$ for all $i=1,...,n$.

 \skipaline \noindent {\bf 3.3.  Diagrams for colored 
$\pi$-tangles.}  It is standard in knot theory to present oriented framed tangles by
plane pictures called tangle diagrams.  A tangle diagram lies in a horizontal strip
 on the page of the picture  identified with ${\bold R}\times [0,1] =
{\bold R}\times 0\times [0,1] \subset {\bold R}^2\times [0,1]$. 
We agree that the first axis is a horizontal line  on the page of the picture
directed to the right, the second axis is orthogonal to the plane of the picture and
is  directed from the reader towards this plane, the third axis is a vertical line on
the plane of the picture   directed from the bottom to the top. 
A tangle diagram  consists  of
oriented immersed arcs and circles   lying in general
position with indication of over/undercrossings in all double points. 
The diagram  has
the same inputs and outputs as the corresponding tangle.   The framing
is given by shifting the tangle along the vector $(0,-\delta,0)$ with small
$\delta>0$.  Note that 
the points with negative second coordinate lie above the picture.

 It is easy to extend the technique of tangle diagrams to present  a $\pi$-tangle
$(T, g:\pi_1(C_{T}) \to \pi)$.  We first present $T$ by a 
 tangle diagram, $D$.  The
undercrossings of $D$  split
$D$ into disjoint oriented embedded arcs  in ${\bold R}\times [0,1]$.  (We do not
break  $D$
at the overcrossings).  For each of these arcs, say $e$, consider the linear 
path
in $ {\bold R}^2\times [0,1]$ connecting the base point $z\in C_{T}$ to a point
of $e+(0,-\delta,0)$.  In the pictorial language, the point $z$ lies high above
  $D$ and the path in question is obtained by rushing from $z$
straight to a point lying slightly above $e$.  Denote this path by
$\gamma(e)$.  We label $e$ with   $g_e=g(\mu_{\gamma(e)})\in \pi$.  
In
pictures, one usually puts   $g_e$ on a small arrow drawn
beneath $e$ and crossing $e$ from right to left.  Knowing $g_e$
for all arcs $e$ of $D$, we can   recover the homomorphism $g$ because the
meridians $\{\mu_{\gamma(e)}\}_e$ generate $\pi_1(C_{T})$.  We say that the
$\pi$-tangle $(T,g)$ is {\it presented by the diagram} $D$   
 whose arcs are labeled by
elements of $\pi$ as above.  

A  tangle diagram $D$  whose arcs are  labeled by elements of $\pi$ presents a
$\pi$-tangle if and only if the following   local condition is
satisfied:

$(\ast)$.  Encircling a double point of $D$ and multiplying the corresponding
four elements of $\pi$ we always obtain $1\in \pi$.

 It is understood that crossing an arc $e$ from right to left we read $g_e$ 
while
crossing $e$ from left to right we read $g_e^{-1}$.

To present a colored $\pi$-tangle $(T,g,u)$ we additionally endow each arc $e$
of $D$ with the object $u_{\gamma(e)}\in {\Cal C}_{g_e}$. 
 This data uniquely determines
$(T,g,u)$.

Conversely, consider a tangle diagram $D$ whose arcs $e$ are labeled with 
pairs
$(g_e\in \pi, U_e\in {\Cal C}_{g_e})$.  Assume that  condition $(\ast)$ is met.  
Then
$D$ presents a colored $\pi$-tangle if and only if the following local 
condition
is satisfied:

$(\ast \ast)$.  For any double point $d$ of $ D$, consider the three arcs
$e,f,h$ of $D$ incident to $d$ such that in a neighborhood of $d$ they appear 
as
the overcrossing, the incoming undercrossing and the outgoing undercrossing,
respectively.  Then $U_f=(\varphi_{g_e})^{\varepsilon} (U_h)$ where
$\varepsilon =\pm 1$ is the sign of $d$.

The necessity of $(\ast \ast)$ follows from  the
definition of a coloring and the obvious equality $\gamma
(f)=(\mu_{\gamma (e)})^{\varepsilon} \gamma (h)$. 
 The sufficiency of $(\ast \ast)$
 follows from Lemma 3.2.1 and the fact that   $(\ast
\ast)$ implies   (3.1.a) for all circle components of $T$.

In the sequel by a diagam of a  colored
$\pi$-tangle we mean its  diagram labeled as above so that conditions 
$(\ast), (\ast \ast)$
are satisfied. The   technique of Reidemeister moves on diagrams of framed
oriented tangles (see, for instance, [Tu2, Section I.4]) extends to the   diagrams of
colored $\pi$-tangles. The key point is that  any Reidemeister move on the
underlying unlabeled diagram  extends   uniquely to a move on the labels keeping  all
the labels outside of the  2-disc where the diagram is modified.   As in the standard
theory,  two diagrams  of colored $\pi$-tangles represent isotopic   colored
$\pi$-tangles if and only if they can be related     by a finite sequence of (labeled)
Reidemeister moves. 

 \skipaline \noindent {\bf 3.4.  Category of colored $\pi$-tangles.}  We 
define
a category of colored $\pi$-tangles ${\Cal T}={\Cal T}(\pi, {\Cal C},K) $ as follows.  
The objects of ${\Cal T}$
are finite sequences $\{(\varepsilon_r,\alpha_r,U_r)\}_{r=1}^k$ 
where
$k\geq 0, \varepsilon_r=\pm, \alpha_r\in \pi, U_r\in {\Cal C}_{\alpha_r}$ for  
$r=1,...,k$.  A morphism of such sequences
$\{(\varepsilon_r,\alpha_r,U_r)\}_{r=1}^k\to \{(\varepsilon^s,\alpha^s,U^s)\}_{s=1}^l$ in
$\Cal T$ is a
finite formal sum $\sum_i z_i T_i$  where
  $z_i\in K$ and  $T_i$ is a colored $\pi$-tangle with $k$
inputs and  $l$ outputs such that  
$\{(\varepsilon_r,\alpha_r,U_r)\}_{r=1}^k$ is the source of $T_i$ and
$\{(\varepsilon^s,\alpha^s,U^s)\}_{s=1}^l$ is the target of $T_i$.  If the sum $\sum z_i 
T_i$ consists of  one term $1\,T$ then we say that the corresponding morphism 
in $\Cal T$ is
represented  by $T$.  The composition of morphisms $T\circ T'$ represented by 
colored $\pi$-tangles $T,T'$ is obtained by the gluing of $T$ on the top of $T'$; this
extends by $K$-linearity to all morphisms.  The identity morphism of an object
$\{(\varepsilon_r,\alpha_r,U_r)\}_{r=1}^k$ is   represented by the trivial colored
$\pi$-braid with zero framing and source (and target) 
$\{(\varepsilon_r,\alpha_r,U_r)\}_{r=1}^k$.
 This completes the definition of
  ${\Cal T}$.

  \skipaline \noindent {\bf 3.4.1.  Lemma.  } {\sl The category ${\Cal T}$ is a
 strict  ribbon   crossed $\pi$-category.}

\skipaline {\sl Proof.}  The tensor product for the objects of ${\Cal T}$ is the
juxtaposition of sequences.  The unit object is the empty sequence.  The tensor
product of the morphisms in ${\Cal T}$ represented by colored $\pi$-tangles $T,T'$
is obtained by  placing any diagram of $T'$ on the right of any diagram 
of
$T$.  Then the union of these two diagrams represents the colored $\pi$-tangle 
$T
\otimes T'$.  This extends to arbitrary morphisms in ${\Cal T}$ by linearity.  In
this way, ${\Cal T}$ becomes a strict monoidal category.

The dual of an object $U=((\varepsilon_1,\alpha_1,U_1),...,(\varepsilon_k,\alpha_k,U_k)
)\in {\Cal T}$ is by definition the object $U^*=((-\varepsilon_k,\alpha_k,U_k),
 (-\varepsilon_{k-1},\alpha_{k-1},U_{k-1}),...,(-\varepsilon_{1},\alpha_{1},U_{1}))$.  The
 duality morphisms $ b_U:  \1\to U\otimes U^*$ and $d_U:U^*\otimes U$  are
defined in the same way as in the usual theory
 where $\pi=1$.  They are represented by  tangle diagrams
 consisting  of $k$ disjoint concentric 
cups (resp.  caps).  The orientations and labels on  these 
cups (resp.  caps) are uniquely determined by the 
data
 at the outputs (resp. inputs). Formulas (1.1.f) and (1.1.g) are straightforward.
 Note that the associativity morphisms and the structural morphisms $l,r$ in
${\Cal T}$ are the identities so that formulas (1.1.f) and (1.1.g) simplify to
$$   (\text {id}_U\otimes d_U)  (b_U 
\otimes
\text {id}_U)  =\text {id}_U, \,\, (d_U\otimes
\text {id}_{U^*})   (\text {id}_{U^*}\otimes b_U)
 =\text {id}_{U^*}.  \leqno (3.4.a) $$

Given $\alpha \in \pi$,   consider a full subcategory ${\Cal T}_\alpha$ of 
${\Cal T}$
whose objects are sequences 
$((\varepsilon_1,\alpha_1,U_1),...,(\varepsilon_k,\alpha_k,U_k))$
with $\prod_{r=1}^k (\alpha_r)^{\varepsilon_r}=\alpha$.  It is obvious that
${\Cal T}=\amalg_{\alpha \in \pi} {\Cal T}_\alpha$. 
 It follows from definitions that ${\Cal T}$ is a $\pi$-category.

For each $\alpha\in \pi$ we define an automorphism $\varphi_\alpha$
 of ${\Cal T}$ as
follows.  The action on the objects is given by
$$\varphi_\alpha(\{(\varepsilon_r,\alpha_r,U_r)\}_{r=1}^k) =\{(\varepsilon_r,\alpha
\alpha_r\alpha^{-1}, \varphi_\alpha(U_r))\}_{r=1}^k$$ where on the right-hand side we
use the action of $\alpha$ on the category of colors ${\Cal C}$.  The automorphism
$\varphi_\alpha:\Cal T \to \Cal T$ transforms a colored $\pi$-tangle $(T,
g:\pi_1(C_T)\to \pi, u)$ into the colored $\pi$-tangle $(T, \alpha g
\alpha^{-1}:\pi_1(C_T)\to \pi, {}^{\alpha}\! u)$ where ${}^{\alpha}\! u$
 is the coloring
of the $\pi$-tangle $(T, \alpha g \alpha^{-1})$ assigning to any path
$\gamma:[0,1]\to C_T$ leading from the base point to $\tilde T$ the object
$\varphi_\alpha(u_\gamma)\in {\Cal C}_{\alpha  g(\mu_\gamma)
\alpha^{-1}}$.  This   extends to arbitrary morphisms in ${\Cal T}$ by 
linearity
and yields an automorphism $\varphi_\alpha$ of ${\Cal T}$.
In this way ${\Cal T}$ becomes a
crossed $\pi$-category.

For  objects $U=\{(\varepsilon_r,\alpha_r,U_r)\}_{r=1}^k \in \Cal T_{\alpha}$ and
$V=\{(\varepsilon^s,\alpha^s,U^s)\}_{s=1}^l \in \Cal T_{\beta}$ with $\alpha, \beta \in
\pi$, the braiding $ U\otimes V \to {}^U\!  V \otimes U$ is represented by the same
(framed oriented) braid $T=T_{k,l}$ on $k+l$ strings as in the standard theory.  The
braid $T$ is given by a diagram  consisting
of two families of parallel linear intervals:  one family consists of $k$
intervals going from $k$   leftmost inputs to $k$ rightmost outputs,
the second family consists of $l$ intervals going from $l$ rightmost inputs
to
$l$ leftmost outputs  and meeting from below all  the
intervals of the first family.  We     extend $T$ uniquely to a colored $\pi$-braid 
with source 
$$U\otimes V= ((\varepsilon_1,\alpha_1,U_1),...,(\varepsilon_k,\alpha_k,U_k),
(\varepsilon^1,\alpha^1,U^1),..., (\varepsilon^l,\alpha^l,U^l)).$$
  It follows from definitions that the target of $T$ is then
$$ ((\varepsilon^1, \alpha \alpha^1 \alpha^{-1}, \varphi_\alpha(U^1)),...,
(\varepsilon^l,\alpha \alpha^l \alpha^{-1}, \varphi_\alpha(U^l)),
(\varepsilon_1,\alpha_1,U_1),...,(\varepsilon_k,\alpha_k,U_k))$$
$$=\varphi_\alpha( V ) \otimes U={}^U\!  V \otimes U.$$
Hence
$T$
 represents a
morphism $U\otimes V \to {}^U\!  V \otimes U$ in ${\Cal T}$   
  denoted $c_{U,V}$.  The axioms   of a braiding are
straightforward. 

For   $U=\{(\varepsilon_r,\alpha_r,U_r)\}_{r=1}^k \in \Cal T_{\alpha}$, the twist 
$$  \theta_U:U \to {}^U\!U
= \varphi_\alpha(U) =\{(\varepsilon_r,\alpha \alpha_r\alpha^{-1},
\varphi_\alpha(U_r))\}_{r=1}^k$$    is
represented by  the
same  framed   braid $ t_k$ as in the standard theory.  It  can
be
obtained from the trivial braid on $k$ strings with constant framing  by one full
right-hand twist.  
 We    extend $t_k$ uniquely to a colored $\pi$-braid 
with source $U$.   It is easy to deduce from definitions that the target of 
this colored $\pi$-braid  is $\varphi_\alpha(U)$ so that it 
represents a morphism $U \to {}^U\!  U$ in ${\Cal T}$. 
 The axioms  
of
a twist are straightforward.  
Thus, ${\Cal T}$ is a ribbon crossed $\pi$-category.

\skipaline \noindent {\bf 3.5.  Ribbon functors.}  To 
formulate our main theorem concerning the category of colored $\pi$-tangles we
need the notion of a ribbon $\pi$-functor. 
Consider   two   crossed $\pi$-categories  ${\Cal C},{\Cal C}'$.  A {\it monoidal 
$\pi$-functor} from ${\Cal C}'$ to ${\Cal C}$  is a covariant functor $F:{\Cal C}'\to {\Cal C}$ such that

(i) $F$ is $K$-linear on morphisms;

(ii) $F(\1_{{\Cal C}'})=\1_{{\Cal C}}$ and $F(f\otimes g)= F(f) \otimes F(g)$ for any morphisms
or objects $f,g$ in ${\Cal C}'$;

(iii) $F$ transforms the structural morphisms $a,l,r$ 
in
${\Cal C}'$ into the corresponding morphisms in ${\Cal C}$:  for any $U,V,W\in {\Cal C}'$,
$$F(a_{U,V,W})=a_{F(U), F(V),F(W)}, \,\,\,\,F(l_{U})=l_{F(U)} , \,\,\,\,
F(r_{U})=r_{F(U)};$$

(iv) $F$ maps ${\Cal C}'_\alpha$ into ${\Cal C}_\alpha$ for all $\alpha\in \pi$;

(v) $F$ is equivariant with respect to the action of $\pi$, i.e.,
$F\circ \varphi_\alpha=\varphi_\alpha \circ F$ for any $\alpha\in \pi$.

Note that we do not require monoidal functors to preserve duality.
If ${\Cal C}$ and ${\Cal C}'$ are strict, then condition (iii)
is superfluous since the structural  morphisms in question are the
identities.

Assume now that ${\Cal C},{\Cal C}'$ are ribbon $\pi$-categories. A monoidal $\pi$-functor
$F:{\Cal C}'\to {\Cal C}$ is said to be {\it ribbon} if it transforms the  
 braiding and twist  in
${\Cal C}'$ into the braiding and twist in ${\Cal C}$, i.e.,   for any $U,V \in {\Cal C}'$,
$$F(c_{U,V})=c_{F(U),F(V)},\,\,\,\,\,
F(\theta_{U})=\theta_{F(U)} .$$

\skipaline \noindent {\bf 3.6.  Theorem.  } {\sl If   ${\Cal C}$ is a 
strict ribbon  crossed 
  $\pi$-category,  then there is a unique ribbon monoidal $\pi$-functor
$F:  {\Cal T}={\Cal T}(\pi, {\Cal C},K) \to {\Cal C}$ such that

(i) for any length 1 object $(\varepsilon,\alpha, U)$ of ${\Cal T}$,
we have $ F((\varepsilon,\alpha, U))=U^{\varepsilon}$ where 

$$ U^{\varepsilon}=\cases
U\in {\Cal C}_\alpha, ~ { {if}}\,\,\, \varepsilon=+, \\ U^*\in {\Cal C}_{\alpha^{-1}} ,~ { {if}}\,\,\,
\varepsilon=-;\endcases $$
 
(ii) for any ${\alpha}\in  \pi, U\in {\Cal C}_{\alpha}$,  
$$F(b_{(+,\alpha, U)})=b_U:  \1_{\Cal C} \to U\otimes U^*\,\,\, {\text and} \,\,\,  F(d_{(+,\alpha, U)})=
d_U:U^*\otimes U \to \1_{\Cal C}.$$}

Theorem 3.6 generalizes the results of [Tu2, Chapter I] where 
$\pi=1$.

The unicity in Theorem 3.6 is quite straightforward.  
The assumption that $F$ is monoidal and condition 
(i)
  determine $F$ on all objects:
$$F((\varepsilon_1,\alpha_1,U_1),...,(\varepsilon_k,\alpha_k,U_k))
 =
\bigotimes_{r=1}^k (U_r)^{\varepsilon_r}.$$
    Next we observe that the morphisms
of type $$ (c_{(\varepsilon,\alpha, U),   (\varepsilon',\alpha',U')})^{\pm 1}, \,\,
(\theta_{(\varepsilon,\alpha, U)})^{\pm 1}, \,\,  b_{(+,\alpha, U)}, \,\,  d_{(+,\alpha,
U)} \leqno (3.6.a) $$ (where $\varepsilon, \varepsilon'=\pm,\, \alpha, \alpha'\in
\pi, \, U\in {\Cal C}_\alpha,  U'\in {\Cal C}_{\alpha'}$) generate the category ${\Cal T}$ in the
sense that any morphism in ${\Cal T}$ can be obtained from these generators by
taking tensor product and composition.  Given the value of $F$ on these generators,
we can compute $F$ on any morphism in ${\Cal T}$.  This proves the uniqueness of
$F$. 

The values of $F$ on the colored $\pi$-tangles $b_{(-,\alpha, U)}, d_{(-,\alpha, U)}$
with ${\alpha}\in  \pi, U\in {\Cal C}_{\alpha}$ can be computed from the   following
equalities in $\Cal T$: $$d_{(-,\alpha, U)}=d_{(+,\alpha,  {}^U\! U)}\, c_{(+,\alpha, 
{}^U\!  U), (-,\alpha,  U)}\,  (\theta_{(+,\alpha, U)}\otimes \id_{(-,\alpha, U)}),\leqno
(3.6.b)$$ $$b_{(-,\alpha, U)}= (\id_{(-,\alpha, U)} \otimes (\theta_{(+,\alpha,
U)})^{-1})
 \,(c_{(-,\alpha, U), (+,\alpha,  {}^U\! U)})^{-1} \, b_{(+,\alpha, U)}.\leqno (3.6.c)$$

A colored  $\pi$-link $ (\ell,z,g,u)$ in
${\bold R}^2\times ]0,1[$  represents an
endomorphism  of the unit object $\1_{\Cal T}$ in ${\Cal T}$ (the empty sequence)
and   is mapped by $F$  into    $F(\ell,z,g,u )\in \End_{\Cal C}
(\1_{\Cal C})$.  It is obvious that any 
colored  $\pi$-link   in $S^3=\bold R^3 \cup \{\infty\}$ is ambiently isotopic
to a unique (up to ambient  isotopy) 
colored  $\pi$-link   in
${\bold R}^2\times ]0,1[$. This allows us to apply  $F$ to colored  $\pi$-links in
$S^3$. For instance,  an empty $\pi$-link $\emptyset$ represents the identity
endomorphism  of   $\1_{\Cal T}$   and therefore
  $F(\emptyset)=1\in \End_{\Cal C} (\1_{\Cal C})$. 

In Section 4 we  generalize   Theorem 3.6 to  colored $\pi$-graphs
and briefly discuss a proof of these theorems.

\skipaline \centerline {\bf 4.  Colored $\pi$-graphs and their invariants}

 \skipaline   The standard 
theory of   colored tangles  generalizes to
so-called ribbon graphs. They are instrumental in the
construction of   3-dimensional TQFT's from modular categories, see
[Tu2].  Here we discuss a similar generalization of colored $\pi$-tangles.

 \skipaline \noindent {\bf 4.1.  Colored $\pi$-graphs.}   We   recall  the
definition of a ribbon graph referring to [Tu2] for details. A  {\it coupon} is a   
rectangle   with a  distinguished  side  called the bottom base; the opposite side is
called   the top base.  A {\it  ribbon graph} $\Omega$ in ${\bold R}^2\times [0,1]$
with $k\geq 0$ inputs $\{(r,0,0)\}_{ r=1}^k$ and $ l\geq 0$ outputs $\{(s,0,1)\}_{
s=1}^l$ consists of a finite family of    arcs, circles and
coupons  embedded in ${\bold R}^2\times [0,1]$.  We call these
arcs, circles and
coupons the {\it strata} of $\Omega$.  The inputs
and   outputs of $\Omega$ should be among the endpoints of the arcs, all the
other endpoints of the arcs should lie on the   bases  of the coupons. 
Otherwise the strata  of $\Omega$ should be disjoint.  At the inputs and outputs,
the arcs should be orthogonal to the planes ${\bold R}^2\times 0, {\bold
R}^2\times1$.  All the  strata  of $\Omega$ should be
oriented and framed so that their framings form  a continuous non-singular vector
field on $\Omega$ transversal to $\Omega$. This vector field is called the {\it
framing} of $\Omega$.
  At the inputs and outputs, the framing should be given by the vector $(0,-1,0)$. 
 On each coupon the framing
should  be transversal to the coupon and  yield together with the orientation of the
coupon the right-handed orientation of ${\bold R}^2\times [0,1]$. 
Slightly pushing  $\Omega$ along its framing we obtain a
disjoint   copy $\tilde \Omega$ of $\Omega$.
Pushing an arc (resp. a
circle, a coupon)  $t$  of $\Omega$ along the framing we obtain an arc (resp. a
circle, a coupon)  $\tilde t\subset \tilde \Omega$.    We
use here a language slightly different but equivalent to the one in [Tu2], where
instead of the framings on arcs and circles we considered orthogonal
2-dimensional bands and annuli.

Let $\Omega\subset {\bold R}^2\times [0,1]$ be a ribbon graph.  We provide its
complement $ C_\Omega=({\bold R}^2\times [0,1])\backslash \Omega$
with the \lq\lq canonical" base point  
$z$ with a big
negative second coordinate $z_2 \ll  0$. As in Section 3, we can  supress the base
point from the notation for   the fundamental group of $  C_{\Omega}$. 

 We   define  meridians of an arc or a
circle   of $\Omega$ exactly as   in Sections 3.1 and 3.2.   With the inputs and
outputs of $\Omega$ we associate canonical meridians as in Section 3.2.   We also
define   meridians of a coupon $Q \subset \Omega$ as follows.  Choose an 
oriented interval $q \subset Q$ leading from the top base of $Q$ to its bottom
base. 
For a path
$\gamma:[0,1] \to C_{\Omega}$ connecting the base point
 to a point   $\gamma(1) \in \tilde Q$,  the   meridian 
$\mu_{\gamma}\in \pi_1(C_{\Omega}) $  is 
represented by the loop $\gamma m \gamma^{-1}$ where $m \subset C_\Omega $
is a small loop    encircling $Q$ (in a plane transversal to  $q$) and having 
linking number $+1$ with $q$.

  A {\it $\pi$-graph}
 is a ribbon graph $\Omega\subset {\bold R}^2\times [0,1]$ 
 endowed with a homomorphism $g:  \pi_1(C_\Omega) \to \pi$. Clearly, a
$\pi$-graph without coupons is nothing but a $\pi$-tangle.

Fix a crossed $\pi$-category $({\Cal C},\varphi:\pi \to \Aut ({\Cal C}))$.
  A {\it coloring}  of  a $\pi$-graph
$(\Omega,g:  \pi_1(C_\Omega) \to \pi)$ consists of two functions $u$ and $v$.
The    
function  $u$ assigns to every arc or circle $t$ of $\Omega$ and to every path
$\gamma:[0,1] \to C_{\Omega}$ connecting the base point  
$\gamma(0) =z\in C_\Omega$ to a
point     $\gamma(1) \in \tilde t\subset \tilde \Omega$ a certain  object
$u_\gamma\in {\Cal C}_{g(\mu_{\gamma})}$  where
$\mu_{\gamma}\in \pi_1(C_{\Omega})$ is the   meridian of $t$ determined by
$\gamma$. This function should satisfy the same conditions as in Section 3, i.e., 

(i)  $u_\gamma$ is preserved under homotopies of $\gamma$ fixing
$\gamma(0)  $ and keeping $\gamma(1)$ on $\tilde {t}$;

(ii) if $\beta$ is a  loop in $(C_{\Omega},z)$ then
$u_{\beta\gamma}=\varphi_{g([\beta])}(u_\gamma)$.

\noindent The   function  $v$
assigns to every coupon  $Q\subset \Omega$ and every  path $\gamma$ in
$C_\Omega$ connecting the base point 
$\gamma(0) =z\in C_\Omega$ to a
point     $\gamma(1) \in \tilde Q \subset \tilde
\Omega$  a certain  morphism $v_\gamma$ in $
{\Cal C}_{g(\mu_{\gamma})}$ satisfying three conditions:

(i)  $v_\gamma$ is preserved under homotopies of $\gamma$ fixing
$\gamma(0) $ and keeping $\gamma(1)$ on $\tilde {Q}$;

(ii) if $\beta$ is a  loop in $(C_{\Omega},z)$ then
$v_{\beta\gamma}=\varphi_{g([\beta])}(v_\gamma)$.

To formulate the third condition on $v_\gamma$ we need some preparations.
Let  $t_1,..., t_m$   be the   arcs of
$\Omega$ incident to the bottom  side of $Q$ 
in the order  determined by the direction of this side induced by the
orientation of    $Q$. Set $\varepsilon_i =+1$   if $t_i$ is directed
  out of $Q$ and $\varepsilon_i =-1$ otherwise. 
Let $t^1,...,t^n$ be the   arcs of $\Omega$
incident to the top  side of $Q$ 
in the order  determined by the direction of this side opposite to the one
induced by the orientation of    $Q$. Set
  $\varepsilon^j =-1$   if $t^j$ is directed
  out of $Q$ and $\varepsilon^j =+1$ otherwise. 
 For $i=1,...,m$,
we can compose the path $\gamma$ with a   path in $\tilde Q$
leading from $\gamma (1)\in \tilde Q$ to the endpoint of $\tilde t_i$ lying on the
bottom side of $\tilde Q$.  The resulting path, $\gamma_i$, leads from the base
point of $C_\Omega$ to a point of  $\tilde t_i$.  
Similarly, for $j=1,...,n$,
we can compose   $\gamma$ with a   path in $\tilde Q$
leading from $\gamma (1)$ to the endpoint of $\tilde t^j$ lying on the
top side of $\tilde Q$.  The resulting path, $\gamma^j$, leads from the base
point of $C_\Omega$ to a point of  $\tilde t^j$.   Note that 
$$ \mu_\gamma=\prod_{i=1}^m  (\mu_{\gamma_i})^{\varepsilon_i} =
 \prod_{j=1}^n  (\mu_{\gamma^j})^{\varepsilon^j} \in \pi_1(C_\Omega).$$
Now we can state the third condition on   $v$:

(iii) for any coupon $Q$ of $\Omega$ and  any path $\gamma$ in $C_\Omega$
connecting the base point   to a point   of   $\tilde Q$,  $$ v_\gamma \in
\Hom_{\Cal C}(\bigotimes_{i=1}^m (u_{\gamma_i})^{\varepsilon_i},\,
 \bigotimes_{j=1}^n (u_{\gamma^j})^{\varepsilon^j}) $$
where   $u_{\gamma_i} \in {\Cal C}_{g(\mu_{\gamma_i})}$
and  $u_{\gamma^j} \in {\Cal C}_{g(\mu_{\gamma^j})}$ are the objects provided by
the function $u$.
 
A $\pi$-graph endowed with a coloring is said to be {\it
colored}.
The definition of the source and the target given in Section 3.2 applies to colored
$\pi$-graphs word for word. 
The technique of  diagrams discussed in Section 3.3 generalizes to colored
$\pi$-graphs in the obvious way. The coupons
should be presented by small  rectangles in ${\bold R}\times [0,1]$ with horizontal
bottom and top bases such that the bottom base has smaller second coordinate
than the top base and the orientation of the coupon is counter-clockwise.

 We define a category of colored $\pi$-graphs $\Cal G= \Cal G(\pi,  {\Cal C},K)$
 following the
lines  of Section 3.4. This category has the same objects as $\Cal T(\pi,  {\Cal C},K)$.
The morphisms in $\Cal G $ are formal linear combinations over $K$ of (the
isotopy classes of) colored $\pi$-graphs 
with the given  source and target. The composition and tensor product are defined
exactly as in $\Cal T$. It is clear that $\Cal T$ is a monoidal subcategory of $\Cal
G$. 

\skipaline \noindent {\bf 4.2.  Remark.}  Sometimes it is  convenient 
to formulate the notions of a $\pi$-graph and a colored $\pi$-graph using the graph
exteriors. The {\it exterior} $E_\Omega$ of a ribbon graph $\Omega\subset
\, {\bold R}^2\times [0,1]  $ is the complement in 
 $  {\bold R}^2\times [0,1]  $ of an open regular neighborhood of $\Omega$.
We    choose  this neighborhood 
so that $\tilde \Omega \subset \partial E_\Omega$.
It is clear that $E_\Omega\subset C_\Omega$ is a deformation retract of
$C_\Omega$
so that $\pi_1(E_\Omega)=\pi_1(C_\Omega)$.  Repeating the definitions of 
Section  4.1
 but considering only    loops and  paths in $E_\Omega$ we
obtain an equivalent  definition  of a colored $\pi$-graph  in terms of the exterior.
  
  \skipaline \noindent {\bf 4.3.  Lemma.  } {\sl The category ${\Cal G}$ is a
 strict  ribbon   crossed   $\pi$-category.}

\skipaline The proof follows the lines of the proof of Lemma 3.4.1. In particular, 
given $\alpha \in \pi$,   the category ${\Cal G}_\alpha$ is defined as the full
subcategory of  ${\Cal G}$
whose objects are sequences 
$((\varepsilon_1,\alpha_1,U_1),...,(\varepsilon_k,\alpha_k,U_k))$
with $\prod_{r=1}^k (\alpha_r)^{\varepsilon_r}=\alpha$.
The
action of   $\alpha\in \pi$ on $\Cal G$ is  obtained by applying $\varphi_\alpha$ to
the objects and morphisms forming the colorings of $\pi$-graphs.  
The braiding, twist and duality in $\Cal G$ come from the corresponding
morphisms in $\Cal T$ via the inclusion $\Cal T \hookrightarrow \Cal G$. 

 \skipaline \noindent {\bf 4.4.  Elementary  colored $\pi$-graphs.  }
Let $V=((\varepsilon_1,\alpha_1,V_1),...,(\varepsilon_k,\alpha_k,V_k))$
and $V'=((\varepsilon^1,\alpha^1,V^1),...,(\varepsilon^l,\alpha^l,V^l))$ be objects of $\Cal
G_\alpha$ with $\alpha\in \pi$.  Let
$$ f  \in
\Hom_{\Cal C}(\bigotimes_{i=1}^m (V_{i})^{\varepsilon_i},\,
 \bigotimes_{j=1}^n (V^{j})^{\varepsilon^j}). $$
With this data we associate a  colored $\pi$-graph $\Omega=\Omega (V,V',f)$ as
follows. Geometrically, $\Omega$ consists of one coupon, $k$ short vertical intervals
attached to its bottom base and $l$ short vertical intervals attached to its top base.
The coupon and the intervals lie in   $\bold R \times [0,1]$, the framing is
constant and orthogonal to $\bold R \times [0,1]$. 
The orientation of the bottom   (resp. top) intervals is 
determined by the signs $\varepsilon_1,...,\varepsilon_k$ 
(resp. $\varepsilon^1,...,\varepsilon^l$).
Note that the group
$\pi_1(C_\Omega)$ is generated by the canonical   meridians $\mu_{1},
...,\mu_k, \mu^{1}, ...,\mu^l \in \pi_1(C_\Omega)$ 
corresponding to the inputs and  
outputs of $\Omega$. They are  subject to   only one relation
$$(\mu_1)^{\varepsilon_1}...(\mu_k)^{\varepsilon_k}
=(\mu^1)^{\varepsilon^1}...(\mu^l)^{\varepsilon^l}. \leqno (4.4.a)$$
Therefore the formulas $g(\mu_r)=\alpha_r , g(\mu^s)=\alpha^s $ where $r=1,...,k$
and $s=1,...,l$ define  a group homomorphism $g:\pi_1(C_\Omega) \to \pi$  sending
both sides of (4.4.a) into $\alpha$.  We define  $ \Omega (V,V',f)$ to be the 
$\pi$-graph $(\Omega,g)$ endowed with   the coloring $(u,v)$ such
that 

(i) the source of $ \Omega (V,V',f)$ is $V$ and the target of $ \Omega (V,V',f)$ is
 $V'$;

(ii)  for the linear path $\gamma$ in $C_\Omega$ leading from the base point 
to the point lying slightly above $\Omega$, we have $v_\gamma=f$.

Such a coloring of $\Omega$ exists and is unique.
We call $ \Omega (V,V',f)$  the {\it elementary colored $\pi$-graph}
 associated with
$V,V',f$.

 \skipaline \noindent {\bf 4.5.  Theorem.  } {\sl If   ${\Cal C}$ is a 
strict ribbon  crossed  
  $\pi$-category, then there is a unique ribbon monoidal $\pi$-functor
$F:  {\Cal G}={\Cal G}(\pi, {\Cal C},K) \to {\Cal C}$ extending the functor
${\Cal T}(\pi, {\Cal C},K) \to {\Cal C}$ of Theorem 3.6 and such that for any
elementary colored $\pi$-graph $ \Omega (V,V',f)$ as above we have $F(\Omega
(V,V',f))=f$. }

 \skipaline

Theorem 4.5 generalizes Theorem I.2.5 in [Tu2] where $\pi=1$. 
 
 The unicity in Theorem 4.5 is   straightforward: it suffices to observe that 
  the morphisms (3.6.a) and the  elementary colored $\pi$-graphs generate 
 ${\Cal G}$.  The proof of
the existence 
 follows the proof of [Tu2, 
 Theorem I.2.5] with natural changes. The key role in the proof  given in [Tu2] is
played by certain commutative diagrams in the target category ${\Cal C}$. We give
here versions of these diagrams in our setting:  for any objects $V,W\in {\Cal C}$
the following two diagrams are commutative  $$\CD V=V\otimes \1@>\text
{id}_V\otimes b_W>>   V\otimes W\otimes W^*@>c_{V,W}\otimes \text {id}_{W^*}>>
{}^V\! W\otimes V\otimes W^*\\  @V\text c_{V,\1}={id}_{V}VV    @Vc_{V,W\otimes
W^*}VV @V\text {id}_{({}^V\! W)} \otimes c_{V,W^*}VV\\
  V=\1\otimes V @> b_{({}^V\! W)}\otimes \text {id}_V >>
{}^V\! W\otimes {}^V\! (W^*)\otimes V @=
{}^V\!  W\otimes {}^V\!  (W^*)\otimes V,
\endCD  $$
$$\CD
 {}^V\!  (W^*)\otimes V\otimes W @<c_{V,W^*} \otimes \text {id}_W<< V\otimes
W^*\otimes W @>\text {id}_V\otimes  d_{W}>> V\otimes
\1=V\\
  @V\text {id}_{{}^V\!  (W^*)} \otimes c_{V,W}VV    @Vc_{V,W^*\otimes W}VV @V 
 \text c_{V,\1}={id}_{V}VV\\
 {}^V\!  (W^*)\otimes {}^V\! W \otimes V @=
{}^V\!  (W^*)\otimes {}^V\! W \otimes V @>   d_{({}^V\! W)}\otimes \text {id}_V>>
\1\otimes V =V.
\endCD $$
One should also use two similar diagrams where the middle vertical arrow is
replaced with $c_{ W^*\otimes W, V}:W^*\otimes W\otimes V \to V\otimes
W^*\otimes W$ and also the diagrams obtained by replacing $V,W$ with $W^*, V$,
respectively.

 \skipaline  \noindent {\bf  4.6. Corollary. }
 {\sl  For any object $U$ of   
 a 
strict ribbon  crossed  
  $\pi$-category  ${\Cal C}$,  there is a canonical isomorphism 
$\eta_U:U\to U^{**}=(U^*)^*$ in ${\Cal C}$. For any morphism $f:U\to V$ in ${\Cal
C}$,
the following diagram is commutative:}
$$\CD U @>\eta_U>>   U^{**}\\  @VfVV    @VVf^{**}V\\
  V  @> \eta_V >> V^{**} . \endCD  $$

\skipaline {\sl Proof.}  If  $U\in {\Cal C}_{\alpha}$ with ${\alpha}\in  \pi$ 
then the morphisms
$$\eta_U= (F(d_{(-,\alpha, U)})\otimes \text {id}_{U^{**}}) (\text {id}_U\otimes
b_{U^*}):U\to U^{**},$$
$$ (d_{U^*}\otimes \text {id}_{U}) (\text {id}_{U^{**}}\otimes
F(b_{(-,\alpha, U)})): U^{**}\to  U $$
    are mutually inverse isomorphisms.
 This can be deduced directly from definitions
using (3.4.a). A pictorial proof  uses Theorem 4.5   and
followes the   lines of the proof of Corollary I.2.6.1 in [Tu2]. 
The second claim is an exercise.

 \skipaline \noindent {\bf 4.7.  Properties of $F$.} We discuss here
three  operations on  a  colored $\pi$-graph $(\Omega,g,u,v)$
preserving   $F(\Omega,g,u,v)$.
Fix 
a circle stratum $\ell$ of
$\Omega$.

1.  The invariant  $F (\Omega,g,u,v)$ does not change when we replace the
 colors of the paths leading to $\ell$ by    isomorphic 
colors. More precisely, fix  a path  $\gamma$
leading from the canonical base point to the longitude $\tilde \ell$ of $\ell$. Let
$V$ be an object of the category ${\Cal C}_{g(\mu_\gamma)}$
 isomorphic to $u_\gamma
\in {\Cal C}_{g(\mu_\gamma)}$. We define a new coloring $(u',v)$ of
$\Omega$ as follows. For the paths leading from the base point to the  strata    
of $\Omega$ distinct from $\ell$, the associated colors are the  same as in   $(u,v)$.
Set $u'_\gamma=V$. This data extends uniquely to a  coloring $(u',v)$
of $(\Omega,g)$. Then $F(\Omega,g,u,v)= F(\Omega,g,u',v)$. To prove this equality  
observe  that   $F(\Omega,g,u,v)$ does not change when we  insert into $\ell$ two
coupons whose colors corresponding to $\gamma$ are    mutually inverse
isomorphisms $u_\gamma \to V$ and $V\to u_\gamma$. Pushing one of the coupons
along $\ell$ and  eventually cancelling it with the second coupon we obtain
$F(\Omega,g,u',v)$.

2. The invariant  $F  (\Omega,g,u,v)$ does not change when we invert the orientation
of 
  $\ell$   and simultaneously replace the
colors associated with   $\ell$  by dual objects. The latter means that for
any path    $\gamma$ leading from the
  base point to   $\tilde \ell$,   we set
$u'_\gamma =(u_\gamma)^*$. For the paths leading  to other
strata of $\Omega$, the color is preserved. This yields a
coloring $(u',v)$ of the $\pi$-graph $(\Omega',g)$ obtained from $(\Omega,g)$ by
inverting the orientation of $\ell$.
Then   $F(\Omega,g,u,v)= F(\Omega',g,u',v)$. The proof follows the lines of the proof
of [Tu2, Corollary I.2.8.1].

3.  Let  ${\ell_1}, {\ell_2}$ be
  two parallel copies of $\ell$ 
determined by the   framing and going very closely to $\ell$.
Consider the ribbon graph ${\Omega'}=(\Omega \backslash \ell) \cup {\ell_1}\cup
{\ell_2}$. To describe  the relevant  colorings of $\Omega'$   we use the language of
 exteriors. The exterior $E'=E_{\Omega'}$ of $\Omega'$  can be obtained  from 
the exterior $E=E_{\Omega}$ of $\Omega$ by
gluing  
  (2-disc with two holes) $\times S^1$ along the 2-torus in $\partial E$
corresponding to $\ell$. Denote the inclusion homomorphism 
$\pi_1(E) \to  \pi_1(E')$ by $h$.
 Assume that there are 
a homomorphism $g':\pi_1(E') \to \pi$ and a coloring
$(u',v')$ of the $\pi$-graph $(\Omega', g')$ satisfying the following conditions:

(i)  $g=g'h:\pi_1(E) \to \pi$;

(ii) on the paths leading from the base point to all strata of $\Omega
\backslash \ell$ the colorings  $(u,v)$ and $(u',v')$ coincide;

 (iii) if $\gamma$ is a  path   in $
E$ leading from the  base point to   $\tilde \ell \subset \partial
E$, then composing $\gamma$ with short paths in
$E'\backslash E$ we obtain    paths $\gamma_1,\gamma_2$ in $
E'$   leading   to   ${\tilde \ell}_1, {\tilde \ell}_2$,
respectively, such that $h(\mu_\gamma) = \mu_{\gamma_1} \mu_{\gamma_2} $
and
$u_\gamma=u'_{\gamma_1}\otimes u'_{\gamma_2}$ where 
 $u'_{\gamma_i}\in  {\Cal C}_{g' (\mu_{\gamma_i})}$ for
$i=1,2$. 

We say   then that the colored $\pi$-graph $({\Omega'}, g',u',v')$ is obtained from
$(\Omega,g,u,v)$ by {\it doubling} of $\ell$. It is easy to check that  
$F(\Omega,g,u,v)= F({\Omega'}, g',u',v')$, cf. 
[Tu2, Corollary I.2.8.3].

\skipaline \centerline {\bf 5.  Trace,  dimension, and algebra of colors  }

 \skipaline We   use Theorem 4.5 to define   the trace,  dimension, and
algebra of colors associated with a ribbon crossed
$\pi$-category  ${\Cal C}$.

 \skipaline \noindent {\bf 5.1.  Trace.}   We define the trace of an  endomorphism
$f:U\to U$ of an object    $U\in {\Cal C}$ by 
$$\tr(f)= d_{ ({}^U\! U)} \,  c_{({}^U\!  U), U^*} \, 
(\theta_{ U} f\otimes \id_{U^*})  \, 
b_{U} \in  \End_{\Cal C}(\1) =\Hom_{\Cal C}(\1,\1).$$
It is clear that for any $k\in K$, we have $\tr(kf)=k\,\tr(f)$.
For any $\beta\in \pi$, we have $\tr(\varphi_\beta(f))= \varphi_\beta (\tr(f))$
where on the right hand side  $\varphi_\beta$ acts on $\End_{\Cal C}(\1)$.

 If 
$U\in {\Cal C}_\alpha$ with $\alpha\in \pi$ then we   can rewrite the definition of
$\tr(f:U\to U)$   using (3.6.b):   
 $$\tr(f)= F(d_{(-,\alpha,U)})\,  (f\otimes
\id_{U^*}) \,  b_{U}  \in  \End_{\Cal C}(\1).$$
Theorem 4.5  gives a geometric interpretation of 
$\tr(f)$. Consider the elementary colored $\pi$-graph $\Omega_f=\Omega
((+,\alpha,U), (+,\alpha,U),f)$ with one input and one output. 
As in the standard theory we can close $\Omega_f$ to obtain  a colored $\pi$-graph
$\hat\Omega_f$ with 0 inputs and 0 outputs. A diagram of $\hat\Omega_f$ is obtained
from any (labeled) diagram of $\Omega_f$
 by connecting the top and bottom endpoints by an arc  
disjoint from the rest of the diagram.  More formally,
$$\hat \Omega_f= d_{(-,\alpha,U)} \, (\Omega_f \otimes
\id_{(-,\alpha,U)}) \, b_{(+,\alpha,U)}.$$
Hence, $$F(\hat \Omega_f) =F(d_{(-,\alpha,U)}) \,  (f\otimes
\id_{U^*}) \,  b_{U}  =\tr(f).$$
Using this geometric interpretation and following the lines of [Tu2, Section I.2.7]
we obtain the following lemma.
 
\skipaline  \noindent {\bf  5.1.1. Lemma.}
Ê (i)  {\sl  For any   morphisms $f:U\to  V$, $g:V\to  U$ in ${\Cal C}$, we
have  } $\text {tr}(fg)=\text {tr}(gf);$ 

(ii)  {\sl for any endomorphisms $f,g$ of   objects
of ${\Cal C}$, we have}     $\text {tr}(f^*)=\text {tr}(f)  $  {\sl and} 
$\text {tr}(f\otimes  g)=\text {tr}(f)\, \text {tr}(g). $ 
  \skipaline

 \skipaline \noindent {\bf 5.2. Dimension.} We   define the dimension of
an object $U\in   {\Cal C}$ by  $$\dim  (U)= \tr (\id_U)  =  d_{ ({}^U\! U)} \, 
c_{({}^U\!  U), U^*} \,  (\theta_{ U}  \otimes \id_{U^*})  \, 
b_{U} \in  \End_{\Cal C}(\1).$$
If $U\in {\Cal C}_\alpha $, then  $\dim  (U)= F(d_{(-,\alpha,U)}) \, 
b_{(+,\alpha,U)}$ is the value of
$F$ on a colored $\pi$-knot represented by a diagram consisting of an embedded
oriented circle labeled with $(\alpha,U)$. This knot is an (oriented) unknot
$K\subset S^3$ endowed with homomorphism $\pi_1(C_K)=\bold Z\to \pi$
sending the meridian to $\alpha$. 

 The properties of $F$ established in Section 4.7 (or Lemma 5.1.1) imply 
that
 isomorphic objects have equal dimensions and 
for any objects $U,V$ and $\beta\in \pi$, we have $\dim  (U^*)=\dim  (U)$,
$\dim  (\varphi_\beta(U))= \varphi_\beta (\dim  (U))$ and 
$\text {dim}(U\otimes  V)=\text {dim}(U)\,\text {dim}(V)$.

Note that for morphisms and objects of the neutral component ${\Cal C}_1\subset {\Cal C}$, the
definitions above coincide with the standard definition of the 
trace and dimension in a
ribbon category. This implies that 
  for    any  $f\in \End_{\Cal C}(\1)$, we have   $ \text {tr}(f)=f$. In particular, 
 $\text {dim}({\1_{\Cal C}})=\tr (\id_{\1})=\id_{\1}$.

 \skipaline \noindent {\bf 5.3.  Algebra of colors.}  We define an {\it algebra of
colors}  or {\it Verlinde algebra} $L=L({\Cal C})$ of  the
ribbon crossed
$\pi$-category  ${\Cal C}$ as follows. 
Consider   the   $K$-module  $\oplus_{U\in {\Cal C}} \End_{\Cal C} (U)$
where $U$ runs over all objects of ${\Cal C}$. The additive generator of
this module represented by $f\in  \End_{\Cal C} (U)$ will be denoted by
$\langle U, f \rangle$ or briefly by $\langle f \rangle$.   We factorize this module by
the following relations:  for any  morphisms $f:U\to V$,  $g:V\to U$ in ${\Cal C}$,
  $$\langle V, fg \rangle=\langle U, gf \rangle. \leqno (5.3.a)$$
Denote the   quotient $K$-module
by $L$. We provide   $L $
with multiplication    by the formula $\langle f \rangle \,\langle f' \rangle= \langle
f\otimes f'\rangle$. Clearly, 
  $L$
is an associative   $K$-algebra with unit $\langle \id_{\1} \rangle  \in L $. 
Every object $U\in {\Cal C} $ determines an element   $
\langle U \rangle =\langle \id_U \rangle\in L $.

The algebra $L$ is $\pi$-graded:    $L=\oplus_{\alpha\in \pi} L_\alpha$ where
$L_\alpha$ is the submodule of $L$ additively generated by elements $\langle U, f
\rangle$ 
 with $U\in {\Cal C}_\alpha$.
We have  $L_\alpha L_\beta\subset L_{\alpha \beta}$ for all $\alpha, \beta \in \pi$.
The formula $ \varphi_\alpha (\langle f \rangle)=
\langle \varphi_\alpha( f)\rangle  $ defines an action of $\pi$ on $L$ by
algebra endomorphisms. Clearly,  $\varphi_\alpha(L_\beta)
=L_{\alpha\beta\alpha^{-1}}$ for all $\alpha, \beta \in \pi$.
The existence of the braiding implies that  $ab=\varphi_\alpha(b)  a$ for
any $a\in L_\alpha, b\in L_\beta$. The existence of the twist implies that
$\varphi_\alpha\vert_{L_\alpha} =\id$ for any $\alpha\in \pi$.

The trace of morphisms in ${\Cal C}$ defines a $K$-linear homomorphism 
$  L\to \End_{\Cal C}(\1)$
sending any generator $\langle f \rangle$
 as above to $\tr(f)\in \End_{\Cal C}(\1)$. We denote
this functional by $\dim$. In particular, $\dim (\langle U \rangle) =
\tr(\id_U)=\dim  (U)$ for any $U\in {\Cal C}$. It follows from the properties of the
trace that
 $\dim $ is an
algebra homomorphism.  

Sending a generator $\langle U, f \rangle$ to the generator $\langle U^*, f^*:U^*\to
U^*\rangle$ we define a $K$-linear homomorphism $L\to L$ denoted by $\ast$.  It 
follows from definitions that $\ast$ is an involutive anti-automorphism of
the algebra $L$ sending each $L_\alpha$ onto $L_{\alpha^{-1}}$ and commuting
with $ \dim: L\to \End_{\Cal C}(\1)$.

 \skipaline \noindent {\bf 5.4.  Special $\pi$-links.}  The importance
of the algebra  of colors $L=L({\Cal C})$ is due to the fact that   its elements can be
used to color so-called special $\pi$-links.  
 Let us call a   $\pi$-link $(\ell ,z , g )$ 
{\it special} if the longitudes of all the components of $\ell$ are sent by $g$ to  
$1\in \pi$. 
 A {\it generalized coloring}
 of a $\pi$-link $(\ell ,z , g )$ is a function   which
assigns to every path $\gamma:[0,1] \to C_{\ell}$ with  $ \gamma(0)=z,
\gamma(1)\in \tilde {\ell}$ an element $a_\gamma\in  L_{g(\mu_{\gamma})}$
such that

(i)  $a_\gamma$ is preserved under homotopies of $\gamma$ fixing
$\gamma(0) $ and keeping $\gamma(1)$ on $\tilde {\ell}$;

(ii) if $\beta $ is a  loop in $(C_{\ell},z)$, then $a_{\beta
\gamma}=\varphi_{g([\beta])}(a_\gamma)$.

For instance, a coloring $u$ of  $(\ell ,z , g )$ as defined in Section
3.1  gives rise to a generalized coloring by  
$a_\gamma =\langle u_\gamma\rangle$. 

 The invariant $F$ extends to 
generalized colorings  $a$  of
  a special $\pi$-link $(\ell=\ell_1\cup ... \cup \ell_n ,z , g )$. 
To define this extension,  
fix a path $\gamma_i$   leading from $z$ to a point of $\tilde \ell_i$ for $i=1,...,n$.
We  transform $\ell$ into a $\pi$-graph $\Omega$  by inserting 
one coupon into each
component   $\ell_i$ near  $\gamma_i (1)$.    If
$a_{\gamma_i}=\langle U_i,f_i\rangle \in L_{g(\mu_{\gamma_i})}$ for   $i=1,...,n$
then we set $u_{\gamma_i}=U_i, v_{\gamma_i}=f_i$ for   $i=1,...,n$. 
This extends uniquely to a coloring $(u,v)$ of  $\Omega$  (it is here that we need
$\ell$ to be special, cf. Lemma 3.2.1). Set $F(\ell,z,g,a)= F(\Omega,g, u,v)$. This
does not depend on the choice of $\gamma_i$ and extends to arbitrary $a$ by 
additivity. Note that if $a_{\gamma_i}= \langle U_i \rangle$ then there is no need
to introduce the coupon on $\ell_i$; it suffices to keep $\ell_i$ as a stratum and to
set $u_{\gamma_i}=U_i$. This gives the same invariant $F$.

Similarly,  we call a   $\pi$-graph  {\it special} if the longitudes of all
its circle strata  are sent  to   $1\in \pi$.  The elements
of $L$   can be used to color  
  circle strata of special $\pi$-graphs. The invariant $F$ extends to 
such generalized colorings of $\pi$-graphs exactly as
above.  The properties of $F$ established in Section 4.7 extend to the 
 generalized colorings in the obvious way.

\skipaline \centerline {\bf 6.  Modular crossed $\pi$-categories }

 \skipaline \noindent {\bf 6.1.  Modular crossed $\pi$-categories.}  Let ${\Cal C}$ be a 
crossed $\pi$-category.   An object $V$ of ${\Cal C}$ is said to be {\it simple}
 if $\End_{\Cal C}
(V)=K\,\id_V$. It is clear that an object isomorphic or dual to a simple object is
itself simple.  The action of $\pi$ on ${\Cal C}$ transforms   simple objects into simple
objects.  

We say that an object $U$ of  
${\Cal C}$ is {\it dominated by simple objects}  if there exist a
finite  set of simple objects $\{V_{r}\}_r$    of   ${\Cal C}$ (possibly
with repetitions) and   morphisms
 $\{f_r:V_{r}\to  U,g_r:U\to  V_{r}\}_r$ such that  
 $\id_U=\sum_r f_rg_r $.  Clearly, if $U\in {\Cal C}_\alpha$ then without loss of
generality we can  assume that  $V_{r}\in {\Cal C}_\alpha$ for all $r$. 

We say that a ribbon
crossed $\pi$-category ${\Cal C}$ is {\it modular} if
it satisfies the following five axioms:

(6.1.1) the unit object $\1_{\Cal C}$ is simple;

(6.1.2) for each $\alpha\in \pi$, the set $I_\alpha$  of the isomorphism  classes of
simple objects of ${\Cal C}_\alpha$ is finite;

(6.1.3)  for each $\alpha\in \pi$, any object  of ${\Cal C}_\alpha$ is dominated by
simple objects of ${\Cal C}_\alpha$;

(6.1.4) if $V,W$ are non-isomorphic simple objects of ${\Cal C}$ then $\Hom_{\Cal C} (V,W)=0$;

To formulate the    last axiom    we need some  notation.  
 For     $i,j \in I_1$, choose 
simple objects  $V_i,V_j\in {\Cal C}_1$ representing $i,j$, respectively, and set $$
S_{i,j}= \tr(c_{V_j,V_i}\circ c_{V_i,V_j}:V_i\otimes V_j\to V_j\otimes
V_i)\in \End_{\Cal C}(\1)=K.$$  
  It follows from the properties of the trace that $S_{i,j}$ does not depend on the
choice of $V_i, V_j$.

(6.1.5).  The square matrix  $S=[S_{i,j}]_{i,j\in I_1} $
is invertible over $K$. 

These axioms generalize the axioms of a modular 
ribbon category given in  [Tu2] where $\pi=1$.  

It follows from axioms (6.1.1)-(6.1.5) that  the neutral component ${\Cal C}_1$ of
${\Cal C}$   is a modular   category in
the sense of [Tu2]. (Note that    (6.1.5) involves only   ${\Cal C}_1$.) 
Recall   the
following property of modular categories: if $U,V$ are objects of a modular
category over $K$ then $\Hom (U,V)$ is a projective $K$-module of finite type. For
any objects $U,V$ of ${\Cal C}$ belonging to the same component of ${\Cal C}$, we have 
$\Hom_{\Cal C} (U,V)=\Hom_{\Cal C} (\1, V\otimes U^*)=\Hom_{{\Cal C}_1} (\1,
V\otimes U^*)$ so that $\Hom_{\Cal C} (U,V)$ is a projective $K$-module of finite
type.  Set $$\mu_{U,V}=\Dim  (\Hom_{\Cal C} (U,V))\in K$$ where $\Dim $ denotes
the usual dimension of  projective $K$-modules of finite type, see for instance [Tu2,
Appendix 1].   We have  $$\mu_{V,U}=\Dim  (\Hom_{\Cal C} (V,U))=\Dim  (\Hom_{\Cal
C} (U,V))=\mu_{U,V}$$ where the second equality follows from Lemma 6.4 below and
the fact that the dimensions of dual projective $K$-modules  are equal.
 Note  also  the identity
$$\mu_{U,V\otimes W} =\mu_{U\otimes W^*,V}$$
for any  $U,V,W\in {\Cal C}$.

Axiom (6.1.1) and the equality  $\varphi_\alpha (\id_{\1})=\id_{\1}$
imply that   any $\alpha\in \pi$ acts in $\End_{\Cal
C} (\1)=K$ as the identity.   Therefore the dimension of objects of ${\Cal C}$ is
invariant under the action of $\pi$: for any   $V\in {\Cal C}, \alpha\in \pi$, we
have   $\dim  (\varphi_\alpha(V))= \varphi_\alpha ( \dim  (V))=  \dim  (V)$. 

If the ground ring $K$ is a field then  axiom (6.1.4) is redundant. It is easy to show
(using for instance Lemma 6.4 below) that in this case any non-zero morphism
between simple objects is an isomorphism. 

 We   need three lemmas concerning a modular crossed $\pi$-category ${\Cal C}$.
The first lemma   computes  the algebra of colors $L({\Cal C})=\oplus_{\alpha\in
\pi} L_\alpha$ as a $K$-module. 

 \skipaline \noindent {\bf 6.2.  Lemma.  } {\sl  
Let $\alpha\in \pi$ and $\{V^\alpha_i\}_{i\in I_\alpha}$ be representatives of the
isomorphism classes of simple  objects in the category ${\Cal C}_\alpha$. Then 
 $L_\alpha$ is a free $K$-module with basis
$\{\langle V^\alpha_i \rangle   \}_{i\in I_\alpha}$.
 For any}
$U\in {\Cal C}_\alpha$,  
 $$\langle U \rangle= \sum_{i\in I_\alpha} \mu_{ V^\alpha_i, U }
\,\langle V^\alpha_i \rangle ,\leqno (6.2.a)$$
$$\dim  (U)= \sum_{i\in I_\alpha} \mu_{ V^\alpha_i, U }\,
\dim  (V^\alpha_i) ,\leqno (6.2.b)$$

\skipaline {\sl Proof.} Let $\langle  U \in {\Cal C}_\alpha,f:U\to U \rangle$ be a
generator of $L_\alpha$. By (6.1.3),
there is a finite system  of objects
 $\{V_{i(r)}\}_r$  belonging to the family $\{V^\alpha_i\}_{i\in I_\alpha}$
(possibly with repetitions) and morphisms 
 $\{f_r:V_{i(r)}\to  U,g_r:U\to  V_{i(r)}\}_r$   such that  
$\id_U=\sum_r f_rg_r$. Then $f=\sum_r f f_rg_r$
and
$$\langle U,f \rangle=\sum_r \langle U, ff_rg_r \rangle= \sum_r 
\langle V_{i(r)}, g_rff_r \rangle.$$
Since $V_{i(r)}$ is a simple object, its endomorphism $g_rff_r  $ equals
$k_r \id_{V_{i(r)}}$ for a certain $k_r \in K$. Therefore
$$\langle U,f \rangle= \sum_r \langle  V_{i(r)}, k_r \id_{V_{i(r)}}\rangle =\sum_r  k_r
\,\langle  V_{i(r)},   \id_{V_{i(r)}}\rangle  =\sum_r k_r \,\langle  V_{i(r)}\rangle .$$
Therefore  the elements 
$\{\langle V^\alpha_i \rangle  \}_{i\in I_\alpha}$ generate  $L_\alpha$ over $K$.
To prove their linear independence we define for each $i\in I_\alpha$ a linear
functional $t_i:L \to K$ as follows. Let $\langle  U\in {\Cal C},f:U\to U \rangle$ be a
generator of $L $. Denote by $f_i$ the automorphism of $\Hom_{\Cal C} (V^\alpha_i, U)$
sending each $h\in \Hom_{\Cal C} (V^\alpha_i, U)$
into $fh$. Set $t_i(f)=\Tr(f_i)\in K$ where $\Tr$ denotes the
  trace of $K$-endomorphisms of projective $K$-modules, see for instance
[Tu2, Appendix 1].  
(If $K$ is a field then $\Tr$ is the usual trace
of matrices.) It follows from the standard  properties of the trace, that $t_i$
annihilates the relation (5.3.a) and
defines thus a linear
functional $t_i:L  \to K$. Note that 
$t_i(\langle U \rangle)= \Dim  (\Hom_{\Cal C} (V^\alpha_i,
U))=\mu_{V^\alpha_i,U}$. By
(6.1.4), $t_i(\langle V^\alpha_j\rangle)= \delta^i_j\in K$  where $\delta^i_j$ is the
Kronecker delta.  This implies that  
$L_\alpha$ is a free $K$-module with basis
$\{\langle V^\alpha_i \rangle   \}_{i\in I_\alpha}$. 
 
For any object $U\in {\Cal C}$ and any $f\in \End_{\Cal C} (U)$, there is a unique
decomposition $\langle U,f \rangle=\sum_{i\in I_\alpha} r_i \,\langle   V^\alpha_i\rangle$
with $r_i\in K$. Applying   $t_i$ we obtain $r_i= \Tr(f_i)$ for all $i$. 
Substituting $f=\id_U$, we obtain
  (6.2.a). Applying the functional
$\dim:L\to \End_{\Cal C}(\1)=K$ to both sides of (6.2.a) we obtain (6.2.b). 

\skipaline  \noindent {\bf  6.3. Corollary.  }  {\sl  $L_\alpha=0$ if and only if the
category ${\Cal C}_\alpha$ is void. } \skipaline  

By a void category we mean a category with empty sets of objects and 
morphisms. Clearly, if ${\Cal C}_\alpha$ is void then $L_\alpha=0$.
Conversely, if $L_\alpha=0$ then by Lemma 6.2, the category ${\Cal C}_\alpha$ has
no simple objects. By (6.1.3),  it is void.

\skipaline  \noindent {\bf  6.4. Lemma.  }  {\sl  For any objects   $U,V$    of ${\Cal C}$,
the pairing} $(x,y)\mapsto \tr(y\circ x):\Hom_{\Cal C}(U,V)\otimes_K \Hom_{\Cal C}(V,U)\to K$ 
{\sl    is non-degenerate. } \skipaline  

We say that a $K$-bilinear pairing $H_1\otimes_K H_2\to K$ is  
non-degenerate if the adjoint homomorphisms
$H_1\to \Hom_K(H_2,K)$ and  $H_2\to \Hom_K (H_1,K)$ are isomorphisms.  

\skipaline  \noindent {\sl Proof of Lemma}  (cf. the proof of Lemma II.4.2.3 in
[Tu2]). Consider first the case $U=\1\in {\Cal C}_1$. If $V\in {\Cal C}_1$ then the claim of the
lemma follows from the standard properties of the modular category ${\Cal C}_1$. 
If $V\in {\Cal C}_\alpha$ with $\alpha\neq 1$ then the claim of the lemma
is obvious since $\Hom_{\Cal C}(U,V)=\Hom_{\Cal C} (V,U)=0$. 
In the case of  an arbitrary   $U$  we have canonical isomorphisms
$$ \alpha:\Hom_{\Cal C}(U,V) \to \Hom_{\Cal C}(\1,V\otimes U^*),\,\,\,\,\,\, \beta:\Hom_{\Cal C}(V,U)\to
\Hom_{\Cal C}(V\otimes U^*,\1)$$ 
such  that  $\tr(yx)=\tr(\beta (y)\alpha(x))$
for any $ x\in \Hom_{\Cal C}(U,V), y\in
\Hom_{\Cal C}(V,U)$. Therefore the general
case of  the lemma follows from the case $U=\1$. 

\skipaline  \noindent {\bf  6.5. Lemma.  }   {\sl  For any simple object $V$ of ${\Cal C}$,
its dimension} $\dim (V)$ {\sl    is invertible in $K$. }\skipaline  

\noindent {\sl Proof}.     With
respect to the generator $\id_V\in \Hom_{\Cal C}(V,V)=K$ 
the  bilinear form 
$(x,y)\mapsto \tr(y x)\in K$ on $\Hom_{\Cal C}(V,V) $
 is presented by the
$(1\times 1)$-matrix $[\tr(\id_V)]=[\dim (V)]$. The non-degeneracy of this form 
implies that  $\dim (V) \in K^*$.

 \skipaline \noindent {\bf 6.6. The canonical color.  }
Let ${\Cal C}$ be a modular crossed $\pi$-category over $K$. For each $\alpha\in \pi$, 
choose  representatives $\{V^\alpha_i\}_{i\in I_\alpha}$ of the
isomorphism classes of simple  objects in the category ${\Cal C}_\alpha$.  We define a
{\it canonical color} $\omega_\alpha\in L_\alpha$ by
$$\omega_\alpha=\sum_{i\in I_\alpha} \dim  (V^\alpha_i) \,\langle V^\alpha_i \rangle.\leqno
(6.6.a)$$ It is clear that $\omega_\alpha$ does not depend on the choice of
representatives  $\{V^\alpha_i\}_{i\in I_\alpha}$. The fact that the duality
$V\mapsto V^*$ transforms simple objects in ${\Cal C}_\alpha$ into simple objects
in ${\Cal C}_{\alpha^{-1}}$ and preserves the dimension
implies that
$$(\omega_\alpha)^*=\omega_{{\alpha}^{-1}}\leqno (6.6.b)$$
for all $\alpha \in \pi$.
Since the action of  any $\alpha\in \pi$ on ${\Cal C}$ transforms   simple objects in
${\Cal C}_\beta$
into simple objects in ${\Cal
C}_{\alpha \beta \alpha^{-1}}$  we have
$$\varphi_\alpha(\omega_\beta)=\omega_{\alpha \beta \alpha^{-1}}.\leqno
(6.6.c)$$ 

The next lemma generalizes a property of canonical colors well-known in the
case $\pi=1$. We follow the proof given in [Bru].

 \skipaline \noindent {\bf 6.6.1.  Lemma.  } {\sl  For any 
$\alpha, \beta\in \pi$ and}  $U\in {\Cal C}_\beta$, $ \omega_\alpha \,\langle U \rangle 
=\dim  (U)\, \omega_{\alpha\beta}$ {\sl and}  $\langle U \rangle \,\omega_\alpha  
=\dim  (U)\, \omega_{\beta\alpha }$.

\skipaline {\sl Proof.}  We  prove only the first equality, the second one is 
 similar. We have  $$ \omega_\alpha \,\langle U \rangle 
=\sum_{i\in I_\alpha} \dim  (V^\alpha_i)  \,\langle V^\alpha_i \rangle \,\langle U \rangle
=\sum_{i\in I_\alpha} \dim  (V^\alpha_i) \,\langle  V^\alpha_i \otimes U\rangle
$$
$$
=\sum_{i\in I_\alpha} 
\sum_{j\in I_{\alpha\beta}} \dim  (V^\alpha_i) \,
 \mu_{ V^{\alpha\beta}_j, V^\alpha_i \otimes U }
\,\langle   V^{\alpha\beta}_j\rangle
=\sum_{i\in I_\alpha} 
\sum_{j\in I_{\alpha\beta}} \dim  (V^\alpha_i) 
\, \mu_{  V^\alpha_i ,V^{\alpha\beta}_j \otimes U^* }
\,\langle   V^{\alpha\beta}_j\rangle$$
$$
=\sum_{j\in I_{\alpha\beta}} 
(\sum_{i\in I_\alpha}  \dim  (V^\alpha_i) \, \mu_{  V^\alpha_i ,V^{\alpha\beta}_j
\otimes U^* }) \,\langle   V^{\alpha\beta}_j\rangle
= 
\sum_{j\in I_{\alpha\beta}} \dim  ( V^{\alpha\beta}_j \otimes U^*)
\,\langle   V^{\alpha\beta}_j\rangle
$$
$$
= \sum_{j\in I_{\alpha\beta}} \dim  ( V^{\alpha\beta}_j)\, \dim  (U^*)
\,\langle   V^{\alpha\beta}_j\rangle=\dim  (U^*)\, \omega_{\alpha\beta}
=\dim  (U)\, \omega_{\alpha\beta}.$$

 \skipaline  \noindent {\bf  6.7. Elements $\Cal D , \Delta_{\pm}$ of $K$. }   In the
sequel we shall need several elements   of $K$ associated
with the neutral component ${\Cal C}_1$ of a modular crossed $\pi$-category ${\Cal C}$ over $K$.
Let $\{V_i\}_{i\in I_1}$ be  representatives of the
isomorphism classes of simple  objects in   ${\Cal C}_1$.  A {\it rank}  of ${\Cal C}_1$ is
an element $\Cal D\in K$   such that $$ {\Cal
D}^2=\sum \limits_{i\in I_1}{(\dim  (V_i))^2}\in K. $$
The existence of a 
rank is a minor technical condition which   does not reduce the range   of
our constructions.  

  Since each $V_i\in {\Cal C}_1$ is a simple object, the  twist  $\theta_{V_i}:
V_i \to \varphi_1(V_i)=V_i$  equals $v_i \,\id_{V_i}$ with  $v_i\in K$.  Since
$\theta_{V_i}$ is invertible, $v_i\in K^*$.
  Set
$$ \Delta_{\pm} =\sum \limits_{i\in I_1}{v_i^{\pm 1}(\dim  (V_i))^2}\in  K. $$
We can interpret  $\Delta_{\pm}\in K$  as the invariant $F$ of an
unknot $\ell  \subset S^3$ endowed with framing $\pm 1$,   trivial
  homomorphism  $\pi_1(C_{\ell})\to 1\subset \pi$ and
 (generalized)  color  $\omega_{1}\in L_1$. 
It is known  that   
 $\Cal D,  \Delta_{\pm} $   are invertible
in   $K$ and that $\Delta_+ \Delta_-=\Cal D^2$  (see   [Tu2, Formula II.2.4.a]).

 \skipaline  \noindent {\bf  6.8. Examples and constructions of modular
$\pi$-categories. } The ribbon $\pi$-category ${\Cal C}={\Cal C} (a,b,c,\theta)$
constructed in Section 2.6 is modular by the obvious reasons:  
each category ${\Cal C}_\alpha$ has only one object which is simple and the
matrix $S=[1]$ is the unit $(1\times 1)$-matrix. If $G\subset \pi$ is a
finite subgroup of the center of $\pi$ then   pushing  ${\Cal C} (a,b,c,\theta)$
forward along the projection $\pi\to \pi/G$ we obtain    a ribbon crossed
$(\pi/G)$-category satisfying   axioms (6.1.1)-(6.1.4) but possibly not 
(6.1.5).  This is a special case of the following 
easy fact: if the kernel of a  group
epimorphism $q:\pi'\to \pi$   is finite and acts as
the identity on a  modular
$\pi'$-category ${\Cal C}'$ then the   push-forward  $\pi$-category $ q_*({\Cal C}')$
(cf. Sections 1.4,  2.5) satisfies all axioms of a modular $\pi$-category except possibly 
(6.1.5). 
On the other hand, the pull-back of any modular
$\pi$-category ${\Cal C}$  along a 
group
homomorphism $ \pi'\to \pi$  is always a modular
   $\pi'$-category.

A tensor product of a finite family of modular crossed $\pi$-categories
is  always modular. The direct product
of $n\geq 2$ modular crossed $\pi$-categories
is not   modular because the unit object in such a direct product is not simple,
$\End(\1)=K^n$. 

It is easy to check that the mirror of a modular $\pi$-category is a modular
$\pi$-category, cf. [Tu2, Exercise II.1.9.2].

Examples of   modular $\pi$-categories with $\pi=\bold Z/2\bold Z$ are
provided by   the categories of representations of $U_q(sl_2(\bold C))$ at roots of
unity. We   use a topological description of these   categories given in  [Tu2] (cf. 
[Th]).   Let $r\geq 3$ be an odd integer and $a$ be a  primitive complex root of unity
of order $4r$. In [Tu2] the author used the tangles and the Jones-Wenzl
idempotents to define 
 a ribbon category $\Cal V(a)$ whose objects  are finite
sequences $(j_1,...,j_l)$ of integer numbers belonging to the set  $\{1,2,...,r-2\}$. We
label   such an object with  $j_1 +...+j_l (\mod 2)  \in \bold Z/2\bold Z=\pi$.   The 
category $\Cal V(a)$ splits as a disjoint union of two full subcategories comprising
objects labelled by 0 and 1, respectively. The crossing isomorphisms
$\{\varphi_\alpha\}_{\alpha \in \pi}$ are the identity maps. The standard ribbon
structure in $\Cal V(a)$  makes $\Cal V(a)$ a ribbon crossed $\pi$-category. Its 
simple objects   are the same as in the standard theory. The results of [Tu2, Section
XII.7.5] show that if $r$ is odd then  $\Cal V(a)$ a modular $\pi$-category.  More general
examples of  modular $\pi$-categories with abelian $\pi$ associated with    quantum
groups will be discussed in [LT].

\skipaline \centerline {\bf 7.  Invariants of 3-dimensional $\pi$-manifolds}

\skipaline \noindent {\bf 7.1.  Principal $\pi$-bundles and   $\pi$-manifolds.}
Let $\pi$ be a (discrete) group. 
A  {\it principal $\pi$-bundle} over a  space 
$M$ is   a   regular covering   $\tilde M \to M$ with group of automorphisms $\pi$.
The spaces $\tilde M$ and $M$ are called the {\it total space} and the {\it base} of
the bundle.  Two principal $\pi$-bundles over  
$M$ are   {\it isomorphic} if there is a homeomorphism  of their total
spaces 
 commuting with the action of $\pi$ and inducing the identity map $M\to M$.

The isomorphism classes of principal $\pi$-bundles over a  manifold   $M$
  are classified by the homotopy classes of
maps from $M$ to the Eilenberg-MacLane space $K(\pi,1)$.   The {\it monodromy} of a
principal $\pi$-bundle  $\xi$ over  $M$ at a point $z\in M$ is the homomorphism
$\pi_1(M,z)\to \pi$   induced by the classifying map  $M\to K(\pi,1)$ of $\xi$.
 The monodromy  of $\xi$ at $z$ is defined up to conjugation by an
element of $\pi$. 
 The monodromy   can be computed   geometrically as follows. Choose
a point    $\tilde z $ lying over $z$ in the total space of $\xi$.  Any loop $\alpha$ in
$(M,z)$ lifts to a path in the total space of $\xi$ beginning in  $\tilde z$ and ending
in 
 $g\tilde z$
with $g=g(\alpha)\in \pi$. The monodromy  of $\xi$ is then the group homomorphism
$\pi_1(M,z)\to \pi$  defined by $[\alpha]\mapsto
g(\alpha)$. A different choice of $\tilde z$ leads to a conjugated
homomorphism. If $M$ is    connected   then we obtain thus a bijective
correspondence between  the isomorphism classes of principal $\pi$-bundles   over 
$M$  and
the 
   group homomorphisms $\pi_1(M,z)\to \pi$ considered up to conjugation. 
The monodromies of $\xi$ at two points  $y,z\in M$
 and the isomorphism
$\pi_1(M,y)\to   \pi_1(M,z)$ induced by any  path  connecting $y$ to $z$   
 form a diagram 
$$\CD \pi_1(M,y) @>>>  \pi_1(M,z) \\ 
 @VVV    @VVV\\
  \pi @=  \pi 
\endCD  $$
  commutative up to conjugation in $\pi$.

  By a {\it  $\pi$-manifold} we mean a pair $(M,\xi)$ where $M$ is a
manifold and $\xi$ is a principal $\pi$-bundle   over  $M$.  We say that  a
$\pi$-manifold  $(M,\xi)$ is  $m$-dimensional (resp.  closed,  connected,  oriented,
etc.) if   $M$ is  $m$-dimensional (resp.  closed,  connected,  oriented,
etc.).   A
{\it homeomorphism of $\pi$-manifolds} $(M,\xi ) \to (M',\xi' )$ is  a 
homeomorphism  of the total spaces of $\xi,\xi'$
 commuting with the action of $\pi$.
Such a homeomorphism   induces a homeomorphism   $M\to M'$.
We will work in the category of oriented manifolds and consider only 
homeomorphisms of $\pi$-manifolds  $(M,\xi ) \to (M',\xi' )$   inducing
orientation preserving homeomorphisms $M\to M'$.

\skipaline \noindent {\bf 7.2.  Invariant $\tau_{\Cal C}$ of
  3-dimensional $\pi$-manifolds.} Let ${\Cal C}$
 be a modular crossed $\pi$-category over
a commutative unital ring $K$. 
We shall
derive from ${\Cal C}$    a homeomorphism invariant 
 $\tau_{\Cal C}$ of closed oriented
3-dimensional $\pi$-manifolds. Let   $L= \oplus_{\alpha\in \pi} L_\alpha$
be the algebra of colors of ${\Cal C}$ and let $\Cal D , \Delta_{\pm}$ be
the elements of $K^*$ associated with ${\Cal C}$ in Section 6.7. 

Let $(M,\xi)$ be a closed connected oriented 3-dimensional
$\pi$-manifold. Present $M$ as the result of surgery on $S^3$ along a framed
 link $\ell$ with  $\#\ell$ components.   
Recall that $M$ is obtained by gluing $\#\ell$
 solid tori to the exterior $E=E_\ell$ of $\ell$.
Take any point $z\in E\subset M$.  
 The inclusion $E\subset M$ induces an epimorphism
$\pi_1(C_\ell,z)=\pi_1(E,z) \to \pi_1(M,z)$ where $C_\ell=S^3\backslash \ell$.
Composing  it  with a homomorphism $\pi_1(M,z)\to \pi$  
representing the monodromy of
$\xi$ at $z$ we obtain a group homomorphism, $g:\pi_1(C_\ell,z) \to \pi$. We fix an
arbitrary orientation of  $\ell$.
 The triple $(\ell,z, g)$ is thus  a $\pi$-link. Clearly,  it is a special
$\pi$-link in the sense of Section 5.4.  We provide  $(\ell,z, g)$  with the following
(generalized) coloring. To every path $\gamma:[0,1] \to C_{\ell}$
with  $ \gamma(0)=z, \gamma(1)\in \tilde {\ell}$, we assign   
$ \omega_{g(\mu_{\gamma})}\in L_{g(\mu_{\gamma})}$. 
By (6.6.c), this satisfies conditions (i), (ii) of Section 5.4 and defines a 
  {\it
canonical   coloring} of $(\ell,z, g)$. Denote the resulting colored
$\pi$-link by  $ \ell_{can}$.
Set
$$\tau_{\Cal C}(M,\xi)= \Delta_{-}^{\sigma(\ell)}
\, {\Cal D}^{ -\sigma(\ell)-\#\ell-1} F(\ell_{can})\in K \leqno (7.2.a)$$
where  $\sigma(\ell)$ is the
signature of the compact oriented   4-manifold $W_\ell$ bounded by $M$
and obtained  from the 4-ball $B^4$ by attaching $n$ 2-handles    along  tubular
neighborhoods of the components of $\ell$ in $S^3=\partial B^4$.  
Here the orientation of $W_\ell$ is induced by the one of $M$.
 We use   the 
\lq\lq outward vector first" convention for the induced orientation: at any point of
$M=\partial W_\ell$ the orientation of $W_\ell$ is determined by the tuple (a
tangent vector directed outwards, a positive basis in the tangent space of $M$).

\skipaline  \noindent {\bf  7.3. Theorem.} {\sl $\tau_{\Cal C}(M,\xi)$ is a homeomorphism
invariant of   $(M,\xi)$.} \skipaline  

\skipaline {\sl Proof.} We should prove that   $\tau_{\Cal C}(M,\xi)$ does not
depend on the choices made in its definition.  
Under the conjugation of 
the   monodromy $\pi_1(M,z)\to \pi$  by an element  $\alpha\in \pi$, the
homomorphism $g$ is replaced with   $\alpha g \alpha^{-1}$. 
It follows from (6.6.c) that the  colored $\pi$-links
$\ell_{can}$ corresponding to $g$ and  $\alpha g \alpha^{-1}$ are  related by the
transformation $\varphi_\alpha$ defined in the proof of Lemma 3.4.1. By
Theorem 3.6 and (3.5.(v)), the corresponding values $F(\ell_{can})\in \End_{\Cal C}
(\1)=K$ are also related via
  $\varphi_\alpha$. 
Since  $\varphi_\alpha $ acts in $\End_{\Cal C} (\1) $ as the identity,
$F(\ell_{can}) $ does not depend on the choice of monodromy in its conjugacy
class.

The independence of the choice of   $z $ follows from the invariance of 
$F(\ell_{can})$   under transfers   described in Section 3.1.

Let $\ell'$ be obtained from $\ell$ by reversing the orientation on one of the
components. By (6.6.b),   $\ell_{can}$ and $\ell'_{can}$ are related by the
transformation described in Section 4.7.2 and  therefore
$F(\ell_{can})=F(\ell'_{can})$.  Therefore $\tau_{\Cal C}(M,\xi)$ does not depend on
the choice of orientation on $\ell$. 

To prove the independence of the choice of $\ell$ we   use   the  Kirby theory of
moves on framed links. It is shown in [Ki] that  any two framed links   in $S^3$
yielding after surgery homeomorphic 3-manifolds can be related by
certain transformations called
  Kirby moves. There are moves of two kinds.
The first    move  adds to a framed link $\ell\subset S^3$ a distant unknot
$\ell^{\pm}$  with framing $ \pm 1$; under this move the
4-manifold  $W_\ell$ is transformed into its connected sum with $ CP^2$. The
second   move preserves   $W_\ell$ and is induced  by a sliding  of a
2-handle  of $W_\ell$ across another 2-handle. We   need a more precise version of
this theory.   For a framed link $\ell\subset S^3$,   denote the result of surgery on
$\ell$ by $M_\ell$. A {\it surgery presentation} of a closed  connected oriented 
3-manifold $M$  is a pair (a framed link $\ell\subset S^3$, an isotopy
class of degree $+1$ homeomorphisms $ f:M\to M_\ell$). Note that any framing
preserving isotopy of $\ell$ onto itself induces a homeomorphism 
$j_0:M_\ell\to M_\ell$. 
Clearly, $(\ell, j_0 f)$ is  a surgery presentation of $M$; we say that
it is obtained from $(\ell, f)$
by isotopy.  The first Kirby
move $\ell\mapsto \ell'= \ell \amalg \ell^{\pm}$ induces a  homeomorphism
$j_1:M_\ell\to M_{\ell'}  $ which is the identity outside a
small 3-ball containing $\ell^{\pm}$.  The second  Kirby move ${\ell}\mapsto
{\ell}'$ induces  a 
diffeomorphisms  $W_{\ell}\to W_{{\ell}'}$ which restricts to a homeomorphism
   $j_2:M_{\ell}\to M_{{\ell}'}$.
In both cases we say that the surgery presentation $(\ell', j_k f:M\to  M_{\ell'})$
(where $k=1,2$) 
  is obtained from $(\ell, f:M\to M_\ell)$
by the $k$-th Kirby move. 
The arguments in  [Ki, Section 2] show  that  for any surgery presentations 
 $({\ell}_1, f_1:M_1\to M_{\ell_1})$  and $({\ell}_2, f_2:M_2\to M_{\ell_2})$ 
of  closed  connected oriented 
3-manifolds $M_1,M_2$ and for any 
isotopy class of degree $+1$ homeomorphisms   $f:M_{1}\to M_{2}$ there is a  
sequence of Kirby moves and isotopies transforming $({\ell}_1, f_1)$ into $
({\ell}_2,f_2 f)$.

The result $M_\ell$ of surgery  on a  special   $\pi$-link   $\ell\subset
S^3$ is a $\pi$-manifold  in the obvious way.
Any Kirby move on a special   $\pi$-link   $\ell\subset S^3$ yields a 
special   $\pi$-link   $\ell'\subset S^3$ where the homomorphism 
$\pi_1(C_{\ell'})\to \pi$ is   the composition of the inclusion
homomorphism $\pi_1(C_{\ell'})\to \pi_1(M_{\ell'})$, the isomorphism
$\pi_1(M_{\ell'}) =\pi_1(M_{\ell})$ induced by the homeomorphism 
$M_{\ell}\to M_{\ell'}$ mentioned above and the homomorphism 
$\pi_1(M_{\ell})\to \pi$ induced by the given homomorphism 
$\pi_1(C_{\ell})\to \pi$.  The results of the previous paragraph
imply that  if two special   $\pi$-links   in $S^3$ yield  after
surgery homeomorphic   $\pi$-manifolds then these 
$\pi$-links  can be
related by  a finite sequence of 
  Kirby moves and isotopies.

It is clear that  $\tau_{\Cal C}(M,\xi)$ is invariant under isotopies on $\ell$.
To prove the theorem it is thus  enough to
show that $\tau_{\Cal C}(M,\xi)$ is invariant under the Kirby moves on $\ell$.  We
begin with the first Kirby  move   $\ell\mapsto \ell'=\ell \amalg \ell^{\pm}$. The
meridian of the unknot $\ell^{\pm}$ is contractible in    $M$
and 
therefore the    $\pi$-link   $ \ell'  $
is a disjoint union of $\ell $ and the framed unknot  $\ell^{\pm} $
endowed with the trivial homomorphism  to $\pi$. The   colored
$\pi$-link   $\ell'_{can}$ is a disjoint union of $\ell_{can}$
and the framed unknot  $\ell^{\pm}_{can}$ endowed with the trivial
homomorphism  to $\pi$ and the color $\omega_1\in  L_1$.
 We have $$F(\ell'_{can})= F(\ell^{\pm}_{can})\,
F(\ell_{can})=\Delta_{\pm } F(\ell_{can}).$$   Therefore the invariance of
$\tau_{\Cal C}(M,\xi)$ under the first Kirby move follows from  the   formulas 
$\#(\ell')=\#\ell +1, \sigma (\ell')=\sigma(\ell) \pm 1$,  $\Delta_+
\Delta_-=\Cal D^2$. 

 We consider the second Kirby moves in the restricted
form studied by Fenn and Rourke, [FR1]. The  Kirby-Fenn-Rourke
moves split into   positive   and   negative ones. It is explained in [RT]
that  (modulo the first Kirby moves) it is
enough to consider only one of these families. Consider for concreteness a
negative Kirby-Fenn-Rourke move   $\ell \mapsto \ell'$.  It replaces a 
  piece $T$ of $\ell$ lying in a ball  by another piece $T'$ lying in the
same ball and having the same endpoints.  Here $T$ is 
  a system of parallel   
strings with parallel framing  and  $T'$  is obtained from $T$ by applying a
full left-hand twist and adding  an unknotted component  $t$ which encircles
$T$ and has framing $-1$.   
Note that 
$\#\ell'=\#\ell +1$ and $\sigma (\ell')=\sigma (\ell)-1$.  We need to prove
that  $F(\ell'_{can})=\Delta_{-}\, F(\ell_{can})$.
This equality follows from a \lq\lq local" equality involving only   $T$ and $T'$. To
formulate this local equality we first position $T$ as a trivial braid  in ${\bold
R}^2\times [0,1]$ with constant framing. The framed tangle $T'=T\cup t\subset
{\bold R}^2\times [0,1]$ 
 is obtained from $T$ as explained above. 
We  orient   $T$ and $t$. Note that  $C_T= ({\bold
R}^2\times [0,1])\backslash T$ is obtained from 
$C_{T'}= ({\bold R}^2\times
[0,1])\backslash {T'}$ by the surgery on  $t$. Therefore any group homomorphism $
g:\pi_1(C_T)\to \pi$ induces a group homomorphism $
g':\pi_1(C_{T'})\to \pi$ such that $g'$ maps the homotopy class of the
$(-1)$-longitude of $t$ into $1\in \pi$. Any coloring $u$ of $(T,g)$
induces a unique  coloring $u'$ of $(T',g')$ such that 
the sources and targets of $T$ and $T'$ coincide and the component $t$ of $T'$  
 has the canonical color.  The local equality mentioned above says  that for
any orientation  of $T\cup t$, any group homomorphism $
g:\pi_1(C_T)\to \pi$, and any coloring $u$ of $(T,g)$,
we have
$$F(T',g',u')= \Delta_{-}\, F(T,g,u).  \leqno (7.3.a)$$
Let us prove this formula. 
Let $(\varepsilon_1,\alpha_1,U_1),...,(\varepsilon_k,\alpha_k,U_k)$ be the   source 
of $T$. 
Using the standard technique of coupons colored with identity
morphisms we can reduce the general case to the case where $T$ consists
of only one string oriented from top to bottom and colored with  object 
$\otimes_{r=1}^k (U_r)^{\varepsilon_r}\in {\Cal C}$. Using  a decomposition of the
identity endomorphism of this object provided by axiom (6.1.3),
we can further reduce ourselves to the case where the string $T$ is 
colored with a  simple object  of ${\Cal C}$.
Thus we can assume that $T$ is a   string 
oriented from top to bottom and  the source  and target  of   $T$ (and $T'$) are a
$1$-term sequence $(+, \alpha , V)$ where $\alpha\in \pi$ and $V$ is a simple
object of $ {\Cal C}_\alpha$. 
By the argument used above, the invariant $F(T',g',u')$ does not change
if we invert the orientation of $t$.  Therefore we can assume that $t$ is
oriented so that its linking number with $T$ equals $-1$.  
Let us denote the   colored $\pi$-tangles $T,T'$ by $T_V, T'_V$. 
Clearly,  $F(T_V)=\id_V$. 
Since the object $V$ is simple,
$F(T'_V)=k\, \id_V$ with $k\in K$. We have to prove only that
$k=\Delta_{-} $.  To this end we 
 close $T'_V$ into   a colored 2-component $\pi$-link
$\hat {T}'_V$.  
As in the standard theory (cf. [Tu2, Corollary I.2.7.1]), 
$$F(\hat {T}'_V)
=\tr(F(T'_V))=\tr(k\,\id_V)=  k\, \tr(\id_V)=k\, \dim 
(V).$$ On the other hand, the  link  $\hat {T}'_V$  is obtained by doubling described in
Section 4.7.3 from an unknot $\ell^{-} \subset S^3$ endowed with framing $-1$,  
trivial   homomorphism  $\pi_1(C_{\ell^{-}})\to 1\subset \pi$ and
  color  $\omega_{\alpha^{-1}} V\in L_1$. By Lemma  6.6.1 and the results of
Section 6.7 the invariant $F$ of this colored $\pi$-unknot $(\ell^{-},\omega_{\alpha^{-1}}
V)$ is computed by  $$ F(\ell^{-},\omega_{\alpha^{-1}} V)=
F(\ell^{-}, \dim  (V) \,\omega_{1}  )= \dim  (V) \,  F(\ell^{-},
 \omega_1) =\dim  (V)\,
\Delta_{-}.$$ Since $F$ is preserved under the doubling, we have  $$F(\hat {T}'_V)=
 F(\ell^{-},\omega_{\alpha^{-1}} V)= \dim  (V)\, \Delta_{-}= 
\Delta_{-}\, \dim  (V).$$  Comparing these two computations of $F(\hat {T}'_V)$ and using 
Lemma 6.5, we obtain  $k=\Delta_{-}$.

 \skipaline  \noindent {\bf  7.4. Computations and remarks.} 
1.  The 3-sphere 
$S^3$ is simply connected and therefore admits a unique structure of a
$\pi$-manifold.
Presenting $S^3$ as the result of surgery  on $S^3$ along an empty
link we obtain $\tau_{\Cal C}(S^3)={\Cal D}^{-1}$. 

2. For each $\alpha \in \pi$, there is a unique $\pi$-structure $\xi_\alpha$ on
$S^1\times S^2$ whose monodromy    along  
$S^1\times pt$ equals $\alpha$. We  prove  in Section 7.7  that
$$\tau_{\Cal C} (S^1\times S^2,\xi_\alpha)=\cases
0,~ { \text {if \, the \, category} \, \,{\Cal C}_\alpha\,\, \text {is\, void} }, \\ 1,~ {
\text {otherwise}}.\endcases \leqno (7.4.a)$$

3.  Formula (7.2.a) can be rewritten in a 
 more symmetric form:   $$
\tau_{\Cal C}(M,\xi)= \Cal D^{-b_1(M)-1}\Delta_-^{-\sigma_-} \Delta_+^{-\sigma_+}
F(\ell_{can})$$  
where $b_1(M)=\#\ell -\sigma_+-\sigma_-$ is the first Betti number of $M$
and $\sigma_+$ (resp. $ \sigma_-$) is the number of positive (resp.
negative) squares in the diagonal decomposition of the intersection form
$H_2(W_\ell) \times H_2(W_\ell) \to \bold Z$. 
This   shows that the invariant
$$\tau'_{\Cal C}(M,\xi)=\Delta_-^{-\sigma_-} \Delta_+^{-\sigma_+}
F(\ell_{can})=\Cal D^{b_1(M)+1}
\tau_{\Cal C}(M,\xi)$$   does not depend on the choice of
$\Cal D$. Note that $\tau'_{\Cal C}(M,\xi)$ is defined  for  a wider class of   ribbon
crossed $\pi$-categories ${\Cal C}$ satisfying (6.1.1)-(6.1.4) and such that 
$\Delta_+, \Delta_-\in K^*$.  The latter condition is a weakened form of (6.1.5):
it follows from  (6.1.1)-(6.1.5) but in general  does not imply (6.1.5).
However,   the invertibility of the matrix $S$ is needed for the construction of a  
TQFT.

4. It is easy to deduce from definitions that the  invariant 
$\Cal D \, \tau_{\Cal C} $ is multiplicative with respect to the connected
sum of $\pi$-manifolds. In other words, for closed 
connected
oriented 3-dimensional $\pi$-manifolds  $M, N$
$$  \tau_{\Cal C} (M\#N)=
{\Cal D}\,\tau_{\Cal C}(M )\,\tau_{\Cal C}(N) $$
where the structure of a $\pi$-manifold on   $M\#N$
is defined so that its monodromy extends the monodromies of $M$ and $N$.

5.  We extend $\tau_{\Cal C}(M,\xi)$   to non-connected
closed oriented 3-dimensional $\pi$-manifolds by multiplicativity
so that $\tau_{\Cal C}(M \amalg  N,\xi) =\tau_{\Cal C}(M ,\xi\vert_M) 
\, \tau_{\Cal C}(N,\xi\vert_N)$.

6. Let $\Cal C= \Cal V(a)$ be the modular $(\bold Z/2\bold Z)$-category   
associated with a    primitive complex root of unity $a$ of order $4r$
 with odd $r\geq 3$
discussed in
Section 6.8. The  principle $\bold Z/2\bold Z$-bundles over
$M$ are numerated by elements   $\xi\in H^1(M;  \bold Z/2\bold Z)$. 
The corresponding  invariants $\tau_{\Cal C}(M ,\xi) $ were first introduced in 
[Bl], [KM], [Tu1]. The invariant $\tau_{\Cal C}(M
,\xi)$  
corresponding to   $\xi=0\in H^1(M;  \bold Z/2\bold Z)$   is called the quantum
$SO(3)$-invariant of $M$.   This example suggests that there should be similar
categorical structures yielding   invariants of spin 3-manifolds
or more generally 3-manifolds endowed with principal bundles over their spaces
of tangent frames.
The author hopes to consider this elsewhere.

 \skipaline  \noindent {\bf  7.5. Extended $\pi$-manifolds.}  The invariant
$\tau_{\Cal C}$ defined above can be generalized to    3-manifolds with partial
$\pi$-structure, i.e.,  with  a principal $\pi$-bundle on the complement of a
framed oriented link or more generally on the complement of a   ribbon graph.  This 
graph   should be colored over   ${\Cal C}$. 
We consider here   only  
3-manifolds   without
boundary and  ribbon graphs  without inputs or outputs,
 the  manifolds with boundary will be discussed in Section 10.
We proceed to precise definitions. 

Let $M$ be a closed connected oriented  3-manifold.    A {\it ribbon  graph}
  in   $M$ consists of a finite number of framed oriented
embedded arcs, circles and coupons which are disjoint, except that the endpoints of
the arcs lie on the bases of the coupons, and the framing satisfies the same
conditions as in Section 4.1.   Since $M$ is connected, the complement of a  ribbon
graph in $M$ is   connected.  A  {\it $\pi$-graph}  in  $M$ 
 is a ribbon graph $\Omega\subset M$ 
whose complement   is endowed with a base point
$z\in M\backslash \Omega$
 and a homomorphism $   \pi_1(M\backslash \Omega,z) \to \pi$.  This data
determines a principal $\pi$-bundle over $M\backslash \Omega$ which may extend
or not to $M$.   We shall use   the language of monodromies rather than
principal bundles; that is why we    pay attention to the base points. Note
that a   $\pi$-graph without coupons is a $\pi$-link as defined in Section  3.1.

Fix a  modular  crossed $\pi$-category $\Cal C$.  A  $\pi$-graph    in $M$
is {\it colored}  (over  $\Cal C$) if it is equipped with 
  two functions $u, v$ satisfying the same conditions as in Section
4.1. The notion of an ambient isotopy in $M$ applies to $\pi$-graphs and colored
$\pi$-graphs in $M$ in the obvious way. This allows us to  consider the (ambient)
isotopy classes of such graphs.
As in Section 3.1, we can transfer the structure
of a colored
$\pi$-graph $\Omega$   along 
paths   in $M\backslash \Omega$ relating various base points. 
The   transfers preserve  the
ambient isotopy class of  a colored $\pi$-graph.

A pair consisting of a closed connected  oriented  3-manifold
  $M$ and a colored
 $\pi$-graph in $  M$  is called  a {\it connected 3-dimensional   extended 
 $\pi$-manifold} (without boundary).
 Let $(M,\Omega)$,  $(M',\Omega')$ be two connected 3-dimensional extended 
  $\pi$-manifolds without boundary.
Here $\Omega=(\Omega \subset M,z, 
g:  \pi_1(M\backslash \Omega,z) \to \pi, u,v)$ is a 
colored $\pi$-graph in $M$ and   
 $\Omega'=(\Omega' \subset M', z', 
g':  \pi_1(M'\backslash \Omega',z') \to \pi, u', v')$ 
is a 
colored $\pi$-graph in $M'$.
By an {\it $e$-homeomorphism} $(M,\Omega)\to (M',\Omega')$
   we mean a degree $+1$ 
homeomorphism of triples $f: (M,\Omega,z)\to (M',\Omega',z')$
 preserving the framing, the orientation  and the splitting of $\Omega, \Omega'$
into strata and such that
  $g' f_{\#}=g:  \pi_1(M\backslash \Omega,z) \to \pi$ and for any 
path $\gamma$ in $M\backslash \Omega$ leading from $z$ to an arc or a circle of
$\tilde \Omega$  (resp. a coupon of $\tilde \Omega$) we have $u'_{f\circ
\gamma}=u_{\gamma}$ (resp. $v'_{f\circ \gamma}=v_{\gamma}$). 
In particular, if $M'=M$ and $\Omega'$ is obtained from $\Omega$ via an ambient
isotopy $\{f_t:M\to M\}_{t\in [0,1]}$ (where $f_0=\id_M$) then 
$f_1$ is an $e$-homeomorphism  $(M,\Omega)\to (M',\Omega')$.

The invariant $\tau_{\Cal C}$    generalizes to 
extended 
  $\pi$-manifolds   as follows.    Let
$\Omega=(\Omega  ,z,   g:  \pi_1(M\backslash \Omega,z) \to \pi, u,v)$ be a
colored $\pi$-graph  in a closed connected   oriented 3-manifold $M$.
Present $M$ as the result of surgery on $S^3$ along a framed
 link $\ell$.    As above, $M$ is obtained by gluing $\#\ell$
 solid tori to the exterior $E$ of $\ell$ in $S^3$.
 Applying   isotopy to  
  $\Omega\subset M$ we  can   deform it into  $E\subset M $.   
 Similarly, we can push the base point $z\in M\backslash \Omega$ into $E$. Thus, we may
assume that $\Omega\subset E$ and $z\in E$.
  The inclusion $E\subset M$
induces an epimorphism $\pi_1(S^3 \backslash (\ell \cup
\Omega),z)=\pi_1(E\backslash  
\Omega ,z) \to \pi_1(M \backslash  
\Omega,z)$.
Composing  this    with   $g$  we obtain a   homomorphism, $\tilde g:\pi_1(S^3
\backslash (\ell \cup \Omega),z) \to \pi$. Fix an arbitrary
orientation of  $\ell$.
 The triple $(\ell \cup \Omega,z, \tilde g)$ is thus  a $\pi$-graph in $S^3$.  
  We equip $\ell$ with the 
  canonical   color as in Section 7.2 and keep the given coloring of $\Omega$. 
(By the inclusion $E\subset M$,  any path in $E$ is also a
path in $M$.)  Denote
the resulting colored $\pi$-graph in $S^3$ by  $ \ell_{can} \cup \Omega$.
Set
$$\tau_{\Cal C}(M, \Omega )= \Delta_{-}^{\sigma(\ell)}
\, {\Cal D}^{ -\sigma(\ell)-\#\ell-1} F(\ell_{can} \cup \Omega)\in K. \leqno
(7.5.a)$$

\skipaline  \noindent {\bf  7.6. Theorem.}
 {\sl Let $\Omega $ be a colored $\pi$-graph  in a closed connected   oriented
3-manifold   $M$. Then $\tau_{\Cal C}(M, \Omega )$ does not depend on the
choices made in its definition and is  an $e$-homeomorphism invariant of  the pair  
$(M,\Omega)$.} \skipaline  

Theorem 7.6  includes Theorem 7.3 as a special case
$\Omega=\emptyset$.
 In this case the   homomorphism $g:  \pi_1(M\backslash \Omega,z) =\pi_1(M ,z)\to
\pi$  provides $M$ with a structure, $\xi$, of a $\pi$-manifold and  $
\tau_{\Cal C}(M, \Omega)=\tau_{\Cal C}(M,  \xi)$. Theorem 7.6 implies that 
$\tau_{\Cal C}(M, \Omega)$ is invariant under ambient isotopies of $\Omega$.
In particular, $\tau_{\Cal C}(M, \Omega)$ is invariant under transfers of the base
point.

The proof of Theorem 7.6 
reproduces the proof of Theorem 7.3 with natural changes (cf. [Tu2, Section II.3]).

Theorem 7.6   implies  that
  $\tau(M,\Omega)$ is an  isotopy invariant 
of     $ \Omega$. Restricting ourselves to colored
$\pi$-graphs consisting   of circles we
obtain   an isotopy invariant of    colored  
$\pi$-links in $M$.

     If $M=S^3$ then $\tau_{\Cal C}(M,\Omega) =\Cal D^{-1} F( \Omega)$ (we
may take $\ell=\emptyset$ to compute $\tau_{\Cal C}(S^3,\Omega)$). The invariant 
$\tau_{\Cal C}(M,\Omega)$ satisfies  the following multiplicativity law: $$
\tau_{\Cal C} (M_1\#M_2,\Omega_1\amalg  \Omega_2)=
{\Cal D}\,\tau_{\Cal C}(M_1,\Omega_1)\,\tau_{\Cal C}(M_2,\Omega_2)  
\leqno (7.6.a) $$
where $\Omega_1,\Omega_2$ are  colored $\pi$-graphs in closed 
connected
oriented 3-manifolds  $M_1,M_2$ respectively. 
The properties of the   invariant $F$  of colored $\pi$-graphs in $S^3$ established in
Section 4.7 generalize in the obvious way to  colored $\pi$-graphs in
closed 3-manifolds.  
 
  Taking   disjoint
unions of connected extended 
 $\pi$-manifolds
  we can obtain       non-connected   extended 
  $\pi$-manifolds.  
 The invariant $\tau_{\Cal C}$ extends to them by multiplicativity. 

\skipaline  \noindent {\bf  7.7. Proof of (7.4.a).} 
Set $$ d_\alpha=\dim (\omega_\alpha)=\sum_{i\in I_\alpha} (\dim
(V^\alpha_i))^2\in K $$  where $I_\alpha$ is the set of isomorphism classes of
simple objects in ${\Cal C}_\alpha$ and $V^\alpha_i$ is a simple  object of ${\Cal
C}_\alpha$ representing $i\in I_\alpha$.
We can obtain   $(S^1\times
S^2, \xi_\alpha)$   by surgery on $S^3$ along a $\pi$-unknot $\ell$ with zero framing and
with homomorphism $\pi_1(C_{\ell})\to \pi$ sending a meridian  of $\ell$ into
$\alpha\in \pi$.  Clearly, $\sigma(\ell)=0$.  By definition, 
    $$\tau_{\Cal C} (S^1\times S^2,\xi_\alpha)= {\Cal D}^{-2} F(\ell_{can} )= 
{\Cal D}^{-2}
d_\alpha. $$
 If the category ${\Cal C}_\alpha$ is void then $\tau_{\Cal C} (S^1\times S^2,\xi_\alpha)=
d_\alpha=0$. Assume that ${\Cal C}_\alpha$ is non-void. It suffices to prove that
 $ d_\alpha= {\Cal D}^{2}$.

By the assumption and axiom (6.1.3), the category 
${\Cal C}_\alpha$ has at
least one simple object. Choose a simple object $U\in {\Cal C}_\alpha$.  
Consider a link diagram in $\bold R^2$ obtained as disjoint union of two 
small embedded  circles $\ell_1,\ell_2$ equipped with  clockwise orientation.
We label $\ell_1,\ell_2$ with $(\alpha, U)$ and 
$(\beta , \omega_\beta)$, respectively, where $\beta\in \pi$.
This determines a special   colored $\pi$-link in $S^3$.
Apply  the surgery on $S^3$ along $\ell_2$. Then $\ell_1$ represents a 
colored $\pi$-knot (in fact an unknot),  $\Omega=\Omega(\alpha, \beta, U)$, in the
result of the surgery $S^1\times S^2$. It follows from definitions (using
multiplicativity of $F$ with respect to disjoint union) that
 $$\tau_{\Cal C} (S^1\times S^2,\Omega)={\Cal D}^{ -2} \dim (U)\, d_\beta .$$
This formula follows also from (7.6.a) 
 if we observe that $(S^1\times S^2,\Omega)$ is a
connected sum
of the $\pi$-manifold  $(S^1\times S^2,\xi_\beta)$ and 
the
extended 
  $\pi$-manifold  $(S^3,  \ell_1)$.
Observe now that the  
extended 
  $\pi$-manifolds 
$(S^1\times S^2,\Omega(\alpha, \beta, U))$ and 
$(S^1\times S^2,\Omega(\alpha, \beta \alpha , U))$ are $e$-homeomorphic. 
 To see this, we   identify   $\Omega$ with the
equatorial  circle 
  in the 2-sphere $pt\times
S^2\subset S^1\times S^2$. This circle splits $pt\times
S^2$ into two half-discs whose regular neighborhoods are 3-balls, containing
$\Omega$. These two 3-balls give rise to two     splittings of $(S^1\times
S^2,\Omega)$ into connected sum  
$$(S^1\times S^2,\Omega)=(S^1\times S^2,\xi_\beta)\# (S^3,  \ell_1)$$
and
$$(S^1\times S^2,\Omega)=(S^1\times S^2,\xi_{\beta \alpha})\# (S^3,  \ell_1).$$
Hence $(S^1\times S^2,\Omega(\alpha, \beta, U))$ and 
$(S^1\times S^2,\Omega(\alpha, \beta \alpha , U))$ are $e$-homeomorphic.
 This implies 
$${\Cal D}^{ -2}  \dim (U)\, d_\beta={\Cal D}^{ -2}  \dim (U)\, d_{\beta \alpha}.$$
Since $\dim (U) \in K^*$, we obtain 
$d_\beta = d_{\beta \alpha}$. For $\beta=1$, this gives 
$d_{\alpha}=d_1={\Cal D}^{ 2}$. 

\skipaline  \noindent {\bf  7.8. Remark.}   The argument in Section 7.7
shows that the set of $\alpha\in \pi$ such that ${\Cal C}_\alpha\neq \emptyset $
is a subgroup of $\pi$ (in fact a normal subgroup).

\skipaline \centerline {\bf 8.   A 2-dimensional homotopy modular functor}

\skipaline    A 
3-dimensional topological quantum field theory (TQFT) derived from a modular
category comprises two   ingredients:   a modular functor assigning $K$-modules 
(called  modules  of conformal blocks) to 
surfaces and an invariant of 3-dimensional 
cobordisms.  The surfaces in this theory have a certain additional structure
consisting of a finite family of marked points and a Lagrangian space in real
1-dimensional homology, 
  such surfaces are said to be {\it extended}. It turns out that a   modular
crossed $\pi$-category gives rise to a homotopy quantum field theory (HQFT) for 
$\pi$-surfaces and $\pi$-cobordisms. This has two   ingredients:   a homotopy 
modular functor assigning $K$-modules to    extended $\pi$-surfaces and an
invariant of 3-dimensional  $\pi$-cobordisms. In this section we   discuss the  
homotopy   modular functor. In the case $\pi=1$ we recover the standard   theory.

\skipaline \noindent {\bf 8.1.    Preliminaries.}  We shall use the language of  
pointed homotopy theory. A topological space is {\it pointed} if all its connected
components are provided with base points.  
A map  between pointed spaces   is a continuous map  sending base points
into base points and considered up to homotopy constant on
the base points.  

We fix an  Eilenberg-MacLane space $X=K(\pi,1)$ associated with   $\pi$
and a base point $x\in X$. We   assume that $X$ is a CW-space.
Note that maps from a pointed connected CW-complex $Y$   into  $X$ bijectively
correspond to group homomorphisms $\pi_1(Y)\to \pi=\pi_1(X,x)$. The
language of maps to $X$ is essentially equivalent to but slightly more convenient
than the language of group homomorphisms into $\pi$ because it allows to treat
connected and non-connected $Y$ on the same footing. 

 We fix a modular
crossed $\pi$-category $\Cal C$.

\skipaline \noindent {\bf 8.2.      Extended $\pi$-surfaces.} 
We first define   extended $\pi$-surfaces without marks. 
An {\it extended $\pi$-surface  without marks} is 
a pointed closed oriented surface $\Upsilon$ endowed with a map
  $\Upsilon\to X=K(\pi,1)$  and with 
a Lagrangian space
$\lambda \subset H_1(\Upsilon;\bold R)$. 
 According to our conventions, the map $\Upsilon\to X$   sends the base points
of all the components of $\Upsilon$ into $x\in X$ and is considered 
up to homotopy constant on
the base points.  
Recall that a  Lagrangian space in  
$  H_1(\Upsilon;\bold R)$ is a linear subspace of maximal dimension 
(equal to  $ \frac{1}{2} \dim H_1(\Upsilon;\bold R)$) on which   the homological
intersection form $ H_1(\Upsilon;\bold R)\times  H_1(\Upsilon;\bold R) \to \bold
R$ is zero.   
 
Now we  define more general extended $\pi$-surfaces with  marks. Let
$\Upsilon$ be a  pointed closed   oriented   surface. A point   $p\in \Upsilon$  is {\it
marked} if is   equipped with a sign $\varepsilon_p=\pm 1$ and a   tangent direction,
i.e.,  a ray $\bold {R_{+}} v$ where $v$ is a non-zero tangent vector at $p$.  A {\it
marking} of  $\Upsilon$ is a finite (possibly void) set of  distinct  marked points
$P\subset \Upsilon$ disjoint from the base points (of the components) of
$\Upsilon$.  Pushing slightly a marked point $p\in P$ in the given tangent direction
we obtain another point $\tilde p \in \Upsilon$ which in analogy with knot theory
can be viewed as a \lq\lq longitude" of $p$.  Set $\tilde P=\cup_{p\in P} \,\tilde p
\subset \Upsilon \backslash P$. Clearly, $\card( \tilde P)= \card (P)$.  For a path
$\gamma:[0,1] \to \Upsilon \backslash P$ connecting a point $z=\gamma(0) $ to   
$\tilde p \in \tilde P$, denote by $\mu_{\gamma}\in \pi_1(\Upsilon \backslash P,z)
$  the  homotopy class of the loop $(\gamma m_p 
\gamma^{-1})^{\varepsilon_p}=\gamma \,m_p^{\varepsilon_p} \,\gamma^{-1}$, where
$m_p $ is a small loop beginning and ending in $\tilde p$ and encircling the point
$p\in P$ in the clockwise direction. (The clockwise
 direction is opposite to the one induced by
the orientation of $\Upsilon$.)

  A {\it $\pi$-marking} of $\Upsilon$  is a marking
$P\subset \Upsilon$ endowed with a map $ g:\Upsilon \backslash P \to X=K(\pi,1)$ 
 sending the base points of  $\Upsilon$
into $x\in X$ and considered 
up to homotopy constant on
the base points. A $\pi$-marking   $P\subset \Upsilon$ is said to be {\it colored}
if it is equipped  with a function  $u$  which
assigns to every path $\gamma:[0,1] \to \Upsilon \backslash P$ leading from a
base point,     $z\in \Upsilon$,   to $\tilde P$   an object $u_\gamma\in
{\Cal C} $ such that

(i)  $u_\gamma$ is preserved under homotopies of $\gamma$ 
in $\Upsilon \backslash P$
fixing $\gamma(0) , \gamma(1)$;

(ii)  $u_\gamma\in
{\Cal C}_{{g_{\#}}(\mu_{\gamma})}$ where ${g_{\#}}:\pi_1(\Upsilon \backslash P,z)\to
\pi=\pi_1(X,x)$ is the homomorphism induced by  
 $ g $;

(iii) if $\beta $ is a  loop in $(\Upsilon \backslash P,z)$, then $u_{\beta
\gamma}=\varphi_{{g_{\#}}([\beta])}(u_\gamma)$.

An {\it extended $\pi$-surface} comprises a pointed closed oriented
surface $\Upsilon$, a colored $\pi$-marking   $P\subset \Upsilon$, and a
Lagrangian space $\lambda=\lambda (\Upsilon) \subset H_1(\Upsilon;\bold R)$.
In the case $P=\emptyset$ we obtain  
an extended $\pi$-surface  without marks as above.
A disjoint union of a finite number of
extended $\pi$-surfaces is an extended $\pi$-surface in the obvious way.
 The empty set is considered as 
an   empty extended $\pi$-surface.

A {\it weak $e$-homeomorphism}  of extended $\pi$-surfaces 
$(\Upsilon, P,g,u, \lambda)\to (\Upsilon', P',g',u'$, $\lambda')$ is a   degree $+1$
homeomorphism  of pairs    $f: (\Upsilon, P)\to (\Upsilon', P')$
such that

- $f$ preserves the base points, the signs of the marked points and their  tangent
directions;

- $g' f =g :\Upsilon \backslash P \to X$ and for any 
path $\gamma$ in $\Upsilon \backslash P  $ leading from a base point  to  
$\tilde P$   we have $u'_{f\circ \gamma}=u_{\gamma}$.

A weak $e$-homeomorphism $f$ as above is called an {\it $e$-homeomorphism} if
 the induced isomorphism $f_*:H_1(\Upsilon;\bold R)\to
H_1(\Upsilon';\bold R)$ maps $\lambda$ onto $\lambda'$.

For any extended $\pi$-surface $\Upsilon$, the opposite $e$-surface
 $-\Upsilon$ is obtained from
$\Upsilon$ by reversing the orientation of $\Upsilon$ and    multiplying the
signs of all the marked points   by $-1$   while keeping   the rest of the data.
Clearly, $-(-\Upsilon)=\Upsilon$. The transformation $\Upsilon \mapsto -\Upsilon$
is natural in the sense that any   (weak) $e$-homeomorphism   $f:\Upsilon\to
\Upsilon'$ gives rise to a (weak)
 $e$-homeomorphism   $-f:-\Upsilon\to -\Upsilon'$ which
coincides with $f$ as a mapping.

\skipaline \noindent {\bf 8.3.      Homotopy modular functor.}  
The modular crossed $\pi$-category $\Cal C$ gives rise to a 2-dimensional
homotopy modular functor $\Cal T=\Cal T_{\Cal C}$. This functor assigns 

- to each extended $\pi$-surface $\Upsilon$
a projective $K$-module  of finite type $\Cal T(\Upsilon)$;

- to each weak $e$-homeomorphism of  extended $\pi$-surfaces $f:\Upsilon  \to
\Upsilon'$ an isomorphism $f_\#:\Cal T(\Upsilon)\to \Cal T(\Upsilon')$.

A construction of $\Cal T$ will be outlined in Section 10. We  state here a  
few simple properties of $\Cal T$: 

(8.3.1) for  disjoint extended $\pi$-surfaces $\Upsilon, \Upsilon'$, there
is a natural  isomorphism $\Cal T (\Upsilon\amalg \Upsilon')=\Cal T(\Upsilon)
\otimes_K 
\Cal T (\Upsilon')$;

(8.3.2) $\Cal T (\emptyset)=K$;

(8.3.3) the isomorphism $f_\# $
associated to a weak $e$-homeomorphism $f$ is invariant under isotopy of $f$ in
the class of  weak $e$-homeomorphisms;

(8.3.4) for any weak 
$e$-homeomorphisms $f:\Upsilon\to \Upsilon'$  and $f':\Upsilon'\to \Upsilon''$,
 we
have    $$(f'f)_{\#}=({\Cal
D}\Delta_{-}^{-1})^{\mu  (f_*(\lambda(\Upsilon)), 
\lambda( \Upsilon'),(f')^{-1}_*(\lambda( \Upsilon'')))} f'_{\#}\, f_{\#}  \leqno
(8.3.a)$$ where $f_*, f'_*$ denote the action of $f,f'$ in the real 1-homology
and $\mu$ denotes the Maslov index of triples of Lagrangian subspaces of 
$H_1(\Upsilon';\bold R)$;

(8.3.5) for any  extended $\pi$-surface  $\Upsilon$,
there   is     a
non-degenerate bilinear 
  pairing $d_\Upsilon:\Cal T(\Upsilon) \otimes_K \Cal T(-\Upsilon) \to K$. The
 pairings   $\{d_\Upsilon\}_\Upsilon$ 
are natural with respect to
weak $e$-homeomorphisms, multiplicative with respect to disjoint unions,
  and
symmetric in the sense that $d_{-\Upsilon}$ is the composition of $d_\Upsilon$
with  the standard flip  $\Cal T(-\Upsilon)\otimes_K  \Cal T(\Upsilon) 
\to \Cal T(\Upsilon)\otimes_K  \Cal T(-\Upsilon) $.

The pairing $d_\Upsilon$
can be used to identify   $(\Cal T(\Upsilon)  )^*=\Hom_K(\Cal T(\Upsilon)
,K)$   with   $\Cal T(-\Upsilon) $. 

 If two extended $\pi$-surfaces
$\Upsilon,\Upsilon'$ 
differ only by the choice of a Lagrangian space in homology then the identity
homeomorphism $\id:\Upsilon\to \Upsilon'$ is a weak 
$e$-homeomorphism and defines thus an isomorphism 
$\Cal T(\Upsilon)\to \Cal T(\Upsilon')$. This   together with (8.3.a) shows that
the projective space associated with $\Cal T(\Upsilon)$ does not depend 
on the choice of $\lambda(\Upsilon)$.
Note also that 
the numerical factor in (8.3.a) equals to $1$ if  
$f$ or $f'$   is an $e$-homeomorphism: in this case    $\mu 
(f_*(\lambda(\Upsilon)),  \lambda( \Upsilon'),(f')^{-1}_*(\lambda( \Upsilon'')))=0$.

A connected extended
$\pi$-surface and the corresponding  module  
can be explicitly described (at least up to isomorphism) as follows. Let
$\Upsilon$ be a closed connected oriented surface of genus $n\geq 0$ with
base point $z$ and $ m\geq 0$ marked points $P=\{ p_1,...,p_m\}\subset  \Upsilon
\backslash \{z\}$. Let $\varepsilon_r=\pm 1$ be the sign of $p_r$. 
For $r=1,...,m$, choose an
 arc  $\gamma_r$ in $\Upsilon
\backslash P$ leading from  $z $ to $\tilde p_r$.  
We   assume that these  $m$  arcs
  are embedded and disjoint
except in  their common endpoint $z$.
Recall
the homotopy class   $\mu_{\gamma_r} \in \pi_1(\Upsilon
\backslash P,z) $  of the loop  encircling   $p_r$, see Section 8.2. 
The group  $\pi_1(\Upsilon \backslash
P,z)$ is known to be generated by  $\mu_{\gamma_1},...,\mu_{\gamma_m}$
and  $ 2n$ elements $ a_1,b_1,..., a_n,b_n$
subject to the only relation 
$$(\mu_{\gamma_1})^{\varepsilon_1}...\,(\mu_{\gamma_m})^{\varepsilon_m}\,
[a_1,b_1]\, ... \,[a_n,b_n]=1 \leqno (8.3.b)$$  where $[a
,b]= ab a^{-1}b^{-1}$.  
Choose    $m+2n$ elements $\mu_1,...,\mu_m,
\alpha_1,\beta_1,..., \alpha_n,\beta_n \in \pi$ satisfying  
$$(\mu_{ 1})^{\varepsilon_1}\,...\, (\mu_{
m})^{\varepsilon_m}\,[\alpha_1,\beta_1]\, ... \, [\alpha_n,\beta_n]=1.$$ Then the
formulas $\mu_{\gamma_r}\mapsto  \mu_r , a_s\mapsto \alpha_s,
b_s\mapsto \beta_s$ with $r=1,...,m; s=1,...,n$ define
 a group homomorphism $\pi_1(\Upsilon
\backslash P,z)\to \pi$ or equivalently a map, $g:\Upsilon
\backslash P\to X$. This makes $P$  a $\pi$-marking. 
For any  sequence of  objects $\{U_r\in {\Cal
C}_{ \mu_r }\}_{r=1}^m$, there is  a unique coloring $u$ of $P$ such
that $u_{\gamma_r}=U_r$ for all $r $. The linear subspace $\lambda$ of $
H_1(\Upsilon;\bold R)$  generated   by the homology classes of   $a_1,...,
a_n$ is a Lagrangian space. The tuple $(\Upsilon,  P,
g, u, \lambda)$ is an extended $\pi$-surface.  Recall the notation $U_r^{+}=U_r, 
U_r^{-}=U_r^*$. Then 
 $$\Cal T(\Upsilon,  P, g, u, \lambda) \leqno (8.3.c) $$
$$=\bigoplus_{i_1\in I_{ \alpha_1 },
...,i_n\in I_{ \alpha_n}}\Hom_{\Cal C} \left (\1, (U_1)^{\varepsilon_1}\otimes ...
\otimes (U_m)^{\varepsilon_m} \otimes \bigotimes_{s=1}^n 
\left  (V_{i_s}^{\alpha_s}\otimes (\varphi_{\beta_s} 
(V_{i_s}^{\alpha_s}))^*\right ) \right)$$ where for $\alpha\in \pi$ we denote by
$I_\alpha$ the set
of isomorphism classes of simple  objects in the category
${\Cal C}_\alpha$ and denote by $\{V^\alpha_i\}_{i\in
I_\alpha}$  certain representatives of these  classes. Note that the group of the
isotopy classes of (weak) $e$-homeomorphisms of 
$(\Upsilon,  P, g, u, \lambda)$ onto itself acts (projectively) on  
$\Cal T(\Upsilon,  P, g, u, \lambda) $.

 It is clear that any extended $\pi$-surface $\Upsilon$ is weakly
$e$-homeomorphic to an extended $\pi$-surface 
$(\Upsilon,  P, g, u, \lambda)$ of the type   described in the previous paragraph. 
Formula
(8.3.c) allows  thus  to compute  
$\Cal T (\Upsilon)$.

Most properties of the   modular functors  known for $\pi=1$ extend
to the general case.  We leave a detailed discussion to another place. Note only
that the splitting formula for the  modules  of conformal blocks   along simple
closed curves on surfaces extends to our setting. Instead of giving a detailed
statement we formulate the key algebraic fact underlying this formula.
 
 \skipaline \noindent {\bf 8.4.  Lemma.  } {\sl  
Let $\alpha\in \pi$ and  $\{V^\alpha_i\}_{i\in I_\alpha}$ be representatives of
the isomorphism classes of simple  objects in the category ${\Cal C}_\alpha$. 
For any objects  $V\in {\Cal C}_\alpha, W\in {\Cal C}_{\alpha^{-1}}$, 
there is a canonical isomorphism}
$$\Hom_{\Cal C} (\1,V\otimes W) =\bigoplus_{i\in I_\alpha} 
\left (\Hom_{\Cal C} (\1,V\otimes  (V^\alpha_i)^*) \otimes_K 
\Hom_{\Cal C} (\1, V^\alpha_i\otimes W) \right ).$$
 
\skipaline {\sl Proof.} This   is equivalent
to  
$$\Hom_{\Cal C} (W^*,V) =\bigoplus_{i\in I_\alpha} 
\left ( \Hom_{\Cal C} (V^\alpha_i,V ) \otimes_K 
\Hom_{\Cal C} (W^*, V^\alpha_i )  \right ).\leqno (8.4.a)$$
This isomorphism   sends $h\otimes h'$ with $h\in 
\Hom_{\Cal C} (V^\alpha_i,V ), h'\in \Hom_{\Cal C} (W^*, V^\alpha_i )$ into
$hh'$. Equality (8.4.a) involves only   morphisms in  the category ${\Cal
C}_\alpha$
and follows from axiom (6.1.3),  cf. [Tu2, Lemma II.4.2.2].

\skipaline \noindent {\bf 8.5.     Action of $\pi$.}  
There is a canonical left action  of   $\pi $ 
on   extended $\pi$-surfaces and (weak) $e$-homeomorphisms.  The
action of 
  $\alpha\in \pi$ transforms an extended $\pi$-surface 
$\Upsilon=(\Upsilon, P,g,u, \lambda)$ into ${}^\alpha\!\Upsilon=(\Upsilon, P,
\alpha_* g, \varphi_\alpha u, \lambda)$ where $\alpha_*:(X,x)\to (X,x)$ is the
  map inducing the conjugation  $\beta\mapsto \alpha\beta\alpha^{-1}$ in
$\pi_1(X,x)=\pi$ and $\varphi_\alpha$ is the given action of $\alpha$ on $\Cal C$.
For a  (weak) $e$-homeomorphism   of extended $\pi$-surfaces 
$f:\Upsilon\to \Upsilon'$, we define ${}^\alpha\! f$ as 
the same map $f$  viewed as a  (weak) $e$-homeomorphism  
${}^\alpha\!\Upsilon\to {}^\alpha\!\Upsilon'$. 

The   modular functor $\Cal T$ can be enriched as follows:
for every extended $\pi$-surface $\Upsilon$ and each $\alpha\in \pi$
there is a   canonical  isomorphism $\alpha_*: \Cal T(\Upsilon)\to \Cal
T({}^\alpha\!\Upsilon)$, see Section
10.3. For  $\alpha, \beta\in \pi$ we have $(\alpha \beta)_*=\alpha_*
\beta_*$. For a   weak  $e$-homeomorphism   of extended $\pi$-surfaces 
$f:\Upsilon\to \Upsilon'$,  we have a  commutative diagram
$$\CD \Cal T(\Upsilon) @>f_\#>>  \Cal T(\Upsilon') \\ 
 @V\alpha_*VV    @VV\alpha_*V\\
  \Cal T({}^\alpha\!\Upsilon) @>({}^\alpha\! f)_\#>>  \Cal T({}^\alpha\!\Upsilon') .
\endCD  $$
  
\skipaline \noindent {\bf 8.6.    Computations on  the torus.}  
Consider formula (8.3.c) in the case where $\Upsilon$ is a   torus without marked
points, i.e.,  $m=0, n=1$.   The  group  $\pi_1(\Upsilon ,z)$ is  generated by  two
elements $ a,b$
subject to  
 $
[a ,b ] =1$. A map $g:\Upsilon
\to X$ is determined by   $ g_{\#} (a)=\alpha, g_{\#} (b)=\beta$ where
$ \alpha, \beta $ are   commuting elements of $\pi$.
 The 1-dimensional   subspace $\lambda$ of $
H_1(\Upsilon;\bold R)$  generated   by the homology class  of   $a $ is a Lagrangian
space. By (8.3.c),
 $$\Cal T(\Upsilon,    g,  \lambda)  
 =\bigoplus_{i \in I_{ \alpha } }\Hom_{\Cal C}  (\1,
V_{i}^{\alpha}\otimes (\varphi_{\beta} 
(V_{i}^{\alpha}))^*)
=\bigoplus_{i \in I_{ \alpha } }\Hom_{\Cal C}  (\varphi_{\beta} 
(V_{i}^{\alpha}),
V_{i}^{\alpha}).  $$
Observe that for any simple objects $U,U'$ of
$\Cal C$, $$ \Hom_{\Cal C} (U, U' )=\cases K,~ {
\text {if} \, \,U \,\, \text {is\, isomorphic \, to} }\,\, U', \\ 0,~ { \text
{otherwise}}.\endcases $$ Note also that $\varphi_{\beta}$ maps ${\Cal C}_\alpha$
into itself and induces a permutation on the set $I_\alpha$. Therefore
$\Cal T(\Upsilon,    g,  \lambda)=K^N$ where $N=\card
\{i \in I_{ \alpha } \, \vert \, \varphi_{\beta}(i)=i\}$.
For instance if $\beta=1$ then $N=\card
 (I_{ \alpha })$.  More generally, if $(\Upsilon,    g, 
\lambda)$ is an extended $\pi$-torus without marks such that the image of 
$g_\#:\pi_1(\Upsilon) \to \pi$ is  a cyclic  group generated by $\alpha\in \pi$, then 
$$\Cal T(\Upsilon,    g,  \lambda)  =K^{\card \,(I_{ \alpha })}.$$ Indeed, we can  
choose the generators $a,b\in \pi_1(\Upsilon)$ as above so that
$\alpha=g_\#(a), \beta=g_\#(b)=1$.

Consider now the case where $\Upsilon$ is a   torus with one marked
point $p$, i.e.,  $m=n=1$. We use the notation of  Section 8.3
 but omit the  index 1 and write $p,
\varepsilon, \gamma, a,b, \mu,\alpha, \beta, U$ for $p_1,\varepsilon_1, \gamma_1, 
a_1,b_1, \mu_1,\alpha_1,
 \beta_1, U_1$.  Thus $P=\{p\}$  and  the  group  $\pi_1(\Upsilon \backslash
P,z)$ is  generated by  three elements
$\mu_{\gamma}, a,b$
subject to  
 $(\mu_{\gamma})^{\varepsilon} \,
[a ,b ] =1$. The map $g:\Upsilon
\backslash P\to X$ is determined by   $g_{\#} (\mu_{\gamma})=
 \mu, g_{\#} (a)=\alpha, g_{\#} (b)=\beta$ where
$\mu, \alpha, \beta\in \pi$ satisfy 
$\mu^{\varepsilon}\,[ \alpha, \beta]=1 $.
The coloring $u$ of $P$ is determined by 
$u_{\gamma}=U$. The 1-dimensional   subspace $\lambda$ of $
H_1(\Upsilon;\bold R)$  generated   by the homology class  of   $a $ is a Lagrangian
space. Then
 $$\Cal T(\Upsilon,  P, g, u, \lambda)  =\bigoplus_{i \in I_{ \alpha } }\Hom_{\Cal C}
\left (\1, U^{\varepsilon}\otimes 
V_{i}^{\alpha}\otimes (\varphi_{\beta} 
(V_{i}^{\alpha}))^* \right).  $$

In general, taking different systems of generators (8.3.b)
of  $\pi_1(\Upsilon \backslash
P,z)$ we can obtain different descriptions of one and the same 
extended $\pi$-surface.
The modules appearing on  the
right-hand side of (8.3.c) are then isomorphic. For instance, in the case $m=n=1$ the
identity $[a ,b ]=  
[aba^{-1} ,a^{-1} ]$ implies that 
the  group  $\pi_1(\Upsilon \backslash
P,z)$ is  generated by  three elements
$\mu_{\gamma}, aba^{-1}, a^{-1}$
subject to  
$ (\mu_{\gamma})^{\varepsilon} \,
[aba^{-1} ,a^{-1} ]=1$. The same map $g:\Upsilon
\backslash P\to X$ as above is determined by   $g_{\#}
(\mu_{\gamma})=
 \mu , g_{\#} (aba^{-1})=\alpha\beta \alpha^{-1}, 
g_{\#} (a^{-1})=\alpha^{-1}$. The  same  coloring $u$  is determined by 
$u_{\gamma}=U$.   Thus under these  two choices of generators of 
$\pi_1(\Upsilon \backslash
P,z)$
we obtain  two extended $\pi$-surfaces which differ only in the  
Lagrangian spaces generated   by the homology classes  of  $a$ and  $b $,
respectively. The corresponding $K$-modules are then isomorphic: 
$$\bigoplus_{i \in I_{ \alpha } }\Hom_{\Cal C} \left (\1, U^{\varepsilon}\otimes 
V_{i}^{\alpha}\otimes (\varphi_{\beta} 
(V_{i}^{\alpha}))^*  \right) \leqno (8.6.a)$$
$$
=\bigoplus_{j \in I_{ \alpha\beta \alpha^{-1} } }\Hom_{\Cal C} \left (\1,
U^{\varepsilon}\otimes 
V_{j}^{\alpha\beta \alpha^{-1}}\otimes (\varphi_{\alpha^{-1}} 
(V_{j}^{\alpha\beta \alpha^{-1}}))^*  \right) $$
$$=\bigoplus_{k \in I_{ \beta  }
}\Hom_{\Cal C} \left (\1, U^{\varepsilon}\otimes 
 \varphi_{\alpha}(V_{k}^{\beta})\otimes
 (V_{k}^{\beta })^*  \right)$$
where the last equality follows from the fact that 
$\varphi_{\alpha}$  induces a bijection $I_{\beta}\to I_{\alpha\beta \alpha^{-1}}$.

\skipaline \centerline {\bf 9.   A 2-dimensional HQFT}

\skipaline    We discuss the 2-dimensional homotopy quantum field theory
(HQFT)  underlying the homotopy 
modular functor of Section 8.

\skipaline \noindent {\bf 9.1.    Preliminaries on   2-dimensional HQFT's.} 
We recall briefly   the notion of a       2-dimensional HQFT  referring
for details to [Tu3]. A {\it 2-dimensional HQFT with target} $X=K(\pi,1)$ has two
ingredients: a family  $\{L_\alpha\}_{\alpha\in \pi}$ of projective $K$-modules of
finite type numerated by the elements of $\pi$ 
 and  a function $\tau$ assigning to surfaces equipped with a map to $X$  
certain   $K$-linear homomorphisms.  More precisely, let 
    $W$ be a compact oriented surface 
 with pointed oriented boundary.  A component $c$ of $\partial W$ is {\it positive}
(resp. {\it negative}) if its given orientation is opposite to the one induced
from $ W$  (resp. coincides with the one  
  induced
from $ W$).  We write $\varepsilon_c=+$ and $\varepsilon_c=-$, respectively. 
 We   view $W$ 
as a  cobordism between    $\cup_{c, \varepsilon_c=-} c  $ and
  $\cup_{c, \varepsilon_c=+} c  $. 
Let $g:W\to X$ be a map
 sending the base points of all the components
of $\partial W$ into a base point $x\in X$. 
Each component  $c\subset \partial W$ is pointed and
oriented so that the map   $g\vert_{c }:c \to X$  is a loop  in $(X,x)$. 
Denote its homotopy class by $\alpha_c\in \pi$.
 The function $\tau$ assigns to each such pair $(W,g)$ a $K$-homomorphism 
$$\tau(W,g): \bigotimes_{c, \varepsilon_c=-} L_{\alpha_c}\to  
\bigotimes_{c, \varepsilon_c=+} L_{\alpha_c}. \leqno (9.1.a)$$ The homomorphism
$\tau(W,g)$ should satisfy a few   axioms. It should be   preserved  under   homotopy
of   $g $ relative to the base points on $\partial W$. It should be  multiplicative
under disjoint union of cobordisms. The gluing of cobordisms should correspond to
composition of homomorphisms. The value of $\tau$ on a cylinder should be the
identity homomorphism.  

In particular, if $W$ is a closed oriented surface   endowed with a map $ g:W\to
X$ then the  
homomorphism   $\tau(W,g):K\to K$ is multiplication by a certain  element of $K$.
This element,  denoted also by  $\tau(W,g) $,  is a homotopy
invariant of   $g$. 

The algebraic counterparts of  2-dimensional HQFT's are crossed $\pi$-algebras over
$K$.  A  $\pi$-graded algebra   or, briefly, a    $\pi$-algebra  over  
$K$ is an associative algebra   $L$ over   $K$
endowed with a splitting     
$L=\bigoplus_{\alpha\in \pi} L_{\alpha}$ such that 
each $L_{\alpha}$ is a projective $K$-module  of finite type, $L_{\alpha}   L_{\beta} \subset
L_{\alpha \beta}$ for any $\alpha, \beta\in \pi$, and
 $L$ has a (right and left) unit $1_L\in L_1$ where     $1$ is the  neutral element of
$\pi$.    
A   
{\it crossed $\pi$-algebra} over    $K$ 
  is a $\pi$-algebra   $L$ over   $K$
endowed with 
a symmetric $K$-bilinear 
  form  (inner product) $\eta : L \times L \to K $ and  a group homomorphism $
\varphi : \pi \to \Aut  (L)$
satisfying the following conditions:

(9.1.1)   $\eta(L_{\alpha}\otimes
L_{\beta})=0$ if $\alpha \beta\neq 1$;     the restriction of  $\eta$ to
$L_{\alpha}\otimes L_{\alpha^{-1}}$ is     non-degenerate for all $\alpha \in \pi$,
and
 $\eta(ab,c)=\eta(a,bc)$ for any $a,b,c\in L$;

(9.1.2) for all $\alpha\in \pi$, $\varphi_{\alpha}=\varphi (\alpha)$ is an algebra
automorphism of $L$ preserving  $\eta$ and  such that $\varphi_{\alpha}
(L_{\beta})\subset L_{ \alpha\beta\alpha^{-1}}$ for all $\beta \in \pi$;

(9.1.3) $\varphi_{\alpha}\vert_{L_{\alpha}}=\id$, for all $\alpha\in \pi$ and 
$\varphi_{\alpha}(b)a=ab$ for any $a\in L_{\alpha}, b\in L_{\beta}$;

(9.1.4) for any $\alpha, \beta \in \pi$ and any $c\in   L_{\alpha\beta \alpha^{-1}
\beta^{-1}}$,
$$\Tr\, (c\, \varphi_{\beta}:L_{\alpha}\to L_{\alpha})=
\Tr\, ( \varphi_{\alpha^{-1}} c:L_{\beta}\to L_{\beta}). $$
Here the homomorphism  on the left-hand side  sends any
$a\in L_{\alpha}$ into $c\, \varphi_{\beta}(a)\in L_{\alpha}$
and the homomorphism  on the right-hand side sends
any
$b\in L_{\beta}$ into $\varphi_{\alpha^{-1}} (cb) \in L_{\beta}$.

According to [Tu3],  for any  2-dimensional HQFT 
$(\{L_\alpha\}_{\alpha\in \pi}, \tau) $ with target $X=K(\pi,1)$
the direct sum  $L=\bigoplus_{\alpha\in \pi} L_{\alpha}$ has the structure
of a crossed $\pi$-algebra. Moreover, this  establishes a 
bijective correspondence between   the isomorphism classes of 
$2$-dimensional HQFT's with target   $K(\pi,1)$ and  the isomorphism classes of
crossed $\pi$-algebras. 

\skipaline \noindent {\bf 9.2.    Underlying 2-dimensional HQFT.} 
The 2-dimensional homotopy modular functor $\Cal T=\Cal T_{\Cal C}$
derived from a modular crossed $\pi$-category  $\Cal C$ has an \lq\lq  underlying"
 2-dimensional HQFT $(\{L_\alpha\}_{\alpha\in \pi}, \tau) $. 
Here $L=\bigoplus_{\alpha\in \pi} L_{\alpha}$ is   the algebra of colors of $\Cal
C$ defined in Section  5.3.  
The values of $\tau$ are obtained, roughly
speaking,   by taking the
dimensions of the modules assigned by  $\Cal T$ to the surfaces.
In particular, for a closed oriented surface  $W$ and  a map $ g:W\to
X$ we have   $\tau(W,g)= \Dim\, \Cal T (\Upsilon_{W,g})\in K$
where $\Upsilon_{W,g}$ is an extended $\pi$-surface 
without marks obtained from $W$
by choosing an arbitrary Lagrangian space in $H_1(W;\bold R)$, choosing  
arbitrary  base points on the components of $W$   and deforming $g$ so that  
it maps these base points into $ x$.  According to the results of Section 8, the
isomorphism class of the module $\Cal T (\Upsilon_{W,g})$ and therefore
its dimension do  not depend on the
choices in the construction of $\Upsilon_{W,g}$.

To describe the homomorphism (9.1.a) for any  pair $(W,g:W\to X)$ as in Section 9.1
we proceed as follows. For each component  $c$ of $\partial W$, the $K$-module
$L_{\alpha_c}$ is free with   basis given by Lemma 6.2. 
(Recall that $\alpha_c\in \pi$ is represented by the loop $g\vert_c$ in $(X,x)$.) 
  We present $\tau(W,g)$ by a matrix with respect to 
the tensor products of these bases. Fix for all  
  $c$   a simple object   
$V_c\in {\Cal C}_{\alpha_c}$. Consider the 
 basis
elements 
$$\bigotimes_{c, \varepsilon_c=-}\langle V_c \rangle\in \bigotimes_{c,
\varepsilon_c=-} L_{\alpha_c}\,\,\, {\text {and }}\,\,\,  
 \bigotimes_{c, \varepsilon_c=+}\langle V_c \rangle\in \bigotimes_{c,
\varepsilon_c=+} L_{\alpha_c}.$$ It suffices to describe the corresponding
matrix term  of $\tau(W,g)$. We denote it by 
$\langle \tau(W,g)\, \vert \, \{V_c\}_c \rangle $.

We upgrade $W$ to  an extended
$\pi$-surface  with  marks as follows. Assume for simplicity that $W$ is connected
and choose a base point $z\in W\backslash \partial W$.  We deform $g:W\to X$
relative to $\partial W$ so that $g(z)=x$. We cap $W$ with 2-discs by gluing to each
component  $c$ of $\partial W$ a copy of the unit complex 2-disc
$D=\{a\in \bold C\,\vert\, \vert a \vert \leq 1\}$. The gluing is
effected so that the point $1\in \partial D$ is identified with  the base point  of
$c$. This results in a closed orientable surface
$\Upsilon$ which we provide with orientation extending the one in $W\subset
\Upsilon$.  The centers (corresponding to $a=0\in D$) of the glued  2-discs form a
finite set $P\subset \Upsilon$. We provide each    $p\in P$ with sign $+$ and 
tangent direction corresponding to $ \bold  {R_+}\subset \bold C$.
  We agree  that 
pushing   $p\in P$ along its tangent direction  we obtain the base point $\tilde p$ of 
the corresponding component,   $c_p$,  of $\partial W$. 
 
Clearly,  $W$ is a deformation retract of 
$\Upsilon \backslash P$. Therefore the map $g:W\to X$ extends to 
a map $\Upsilon \backslash P$ denoted also by $g$. This makes $P$ a $\pi$-marking
on $\Upsilon$. We color $P$ as follows.   For a path $\gamma:[0,1] \to
\Upsilon \backslash P$ leading from       $z $    to $\tilde p$,  the path
$g\circ \gamma$ is a {\it loop} in $ (X,x)$ so that we can consider its homotopy
class $[g\circ \gamma]\in \pi$. Set $$
 u_\gamma
=\cases
\varphi_{[g\circ \gamma]}
(V_{c_p}) ,~ { \text {if} \, \,\varepsilon_{c_p}=+1}, \\ 
\varphi_{[g\circ \gamma]}
(V^{*}_{c_p}) ,~ { \text {if} \, \,\varepsilon_{c_p}=-1 }.\endcases
\leqno (9.2.a)$$ 
 Conditions (i) - (iii) of Section 8.2 are straightforward. In particular, Condition (ii)
follows from the equality
$${g_{\#}}(\mu_{\gamma})= [g\circ \gamma] (\alpha_{c_p})^{\varepsilon_{c_p}} 
[g\circ \gamma]^{-1}.$$
Choosing  an arbitrary Lagrangian space   $\lambda \subset H_1(\Upsilon;\bold R)$
we obtain an extended $\pi$-surface $(\Upsilon, g, P, u, \lambda)$.
Set
$$\langle \tau(W,g)\, \vert \, \{V_c\}_c \rangle= \Dim \, \Cal T (
\Upsilon, g, P, u, \lambda).$$

One can check that  $(\{L_\alpha\}_{\alpha\in \pi}, \tau) $ 
is a  2-dimensional HQFT with target $X$. We  
describe here the structure of a crossed $\pi$-algebra in $L=\bigoplus_{\alpha\in
\pi} L_{\alpha}$ underlying this  2-dimensional HQFT. First,  define a pairing $\eta :
L \times L \to K $ by 
 $\eta( \langle U,f\rangle, \langle U',f'\rangle)=
\Tr \,  (f\otimes f')_* $ where $\langle U,f\rangle, \langle U',f'\rangle$ are
arbitrary additive generators of $L$ as in Section 5.3 and $(f\otimes f')_*$ denotes 
the endomorphism of $\Hom_{\Cal C} (\1, U\otimes U')$ sending each 
$h\in \Hom_{\Cal C} (\1, U\otimes U')$ into $(f\otimes f')h$.  It is easy to check that
$\eta$ is  a well-defined bilinear pairing. 

\skipaline  \noindent {\bf  9.3. Theorem.} {\sl The $\pi$-graded algebra 
$L=\bigoplus_{\alpha\in \pi} L_{\alpha}$ with
form $\eta$ and the action $\varphi$ of $\pi$ defined in Section 
5.3  is a  crossed $\pi$-algebra.} \skipaline  

\skipaline {\sl Proof.}  
For each $\alpha\in \pi$, fix representatives $\{V^\alpha_i\}_{i\in I_\alpha}$  of
the isomorphism classes of simple  objects in the category ${\Cal C}_\alpha$. 
By Lemma 6.2,    $L_\alpha$ is a free
$K$-module with basis $\{\langle V^\alpha_i \rangle   \}_{i\in I_\alpha}$.

Let us check (9.1.1).  Observe that for any objects $U,U'$ of
$\Cal C$, $$\eta ( \langle U \rangle, \langle U' \rangle)= \Dim (\Hom_{\Cal C} (\1, U\otimes
U'))=\Dim (\Hom_{\Cal C} (U^*, U')).$$ In particular if the objects $U,U'$ are simple
then
$$\eta ( \langle U \rangle, \langle U' \rangle)
=\cases
1,~ { \text {if} \, \,U' \,\, \text {is\, isomorphic \, to} }\,\, U^*, \\ 0,~ {
\text {otherwise}}.\endcases \leqno (9.3.a)$$
This implies that
 $\eta(L_{\alpha}\otimes
L_{\beta})=0$ if $\alpha \beta\neq 1$ and the bases   
$\{\langle V^\alpha_i \rangle   \}_{i\in I_\alpha}$
and
$\{\langle V^{\alpha^{-1}}_i \rangle   \}_{i\in I_{\alpha^{-1}}}$   are dual to each
other with respect to $\eta$. Hence the    non-degeneracy of $\eta$.
Formula $\eta(ab,c)=\eta(a,bc)$ follows from the definition of $\eta$.

By (9.3.a) and Corollary 4.6, the form $\eta$ is symmetric.  Formula (9.3.a) implies
that $\eta$ is invariant under the action of $\pi$. All other  conditions in (9.1.2),
(9.1.3) were checked in Section 5.3.
 
 To check (9.1.4) let $c=\langle U \rangle$ where $U$ is an object of
${\Cal C}_{\alpha\beta \alpha^{-1}
\beta^{-1}}$. The homomorphism
$c\, \varphi_{\beta}:L_{\alpha}\to L_{\alpha}$ sends  
$\langle V^\alpha_i \rangle $ into $$\langle U\otimes
\varphi_\beta(V^\alpha_i)\rangle = \sum_{j\in I_\alpha} \mu_{V^\alpha_j, 
U\otimes
\varphi_\beta(V^\alpha_i)} \langle V^\alpha_j \rangle.$$
Therefore
$$\Tr\, (c\, \varphi_{\beta})= \sum_{i\in I_\alpha}
\mu_{V^\alpha_i, 
U\otimes
\varphi_\beta(V^\alpha_i)}\leqno (9.3.b)$$
$$=\sum_{i\in I_\alpha}
\mu_{ 
U\otimes
\varphi_\beta(V^\alpha_i), V^\alpha_i}= 
\sum_{i\in I_\alpha}
\Dim (\Hom_{\Cal C} ( \1,
U^*\otimes V^\alpha_i \otimes
(\varphi_\beta (V^\alpha_i))^* ) ).$$
A similar computation shows that
$$\Tr\, ( \varphi_{\alpha^{-1}} c:L_{\beta}\to L_{\beta})=
\sum_{k\in I_\beta}
\Dim (\Hom_{\Cal C} ( \1,
U^*\otimes \varphi_\alpha(V^\beta_k) \otimes
 (V^\beta_k)^*)). \leqno (9.3.c)$$
By (8.6.a) (where $U$ should be replaced with $U^*$), the right-hand sides of (9.3.b)
and (9.3.c) are the dimensions of isomorphic $K$-modules and therefore they
are equal.  

\skipaline  \noindent {\bf  9.4. Remarks.} 1. The proof of   (9.1.4) in Theorem 9.3
uses the existence of the homotopy modular functor associated with $\Cal C$. It
would be useful to have a direct algebraic proof.

2. The  neutral component  $L_1$ of the algebra of colors of $\Cal C$ is known to be
a semisimple commutative $K$-algebra (cf.  [Tu2, Section IV.12.4]). More precisely,
the algebra
$L_1$ splits as a direct sum of $\card (I_1)$ copies of $K$. This and the results of
[Tu3] imply that if $K$ is a field of characteristic 0 then   the HQFT
$(\{L_\alpha\}_{\alpha\in \pi}, \tau) $ is semi-cohomological, i.e., is determined by
a finite family of $K$-valued weights and 2-dimensional cohomology classes of
subgroups of $\pi$ of finite index.

\skipaline \centerline {\bf 10.   A 3-dimensional HQFT}

 \skipaline 

Each
modular crossed
 $\pi$-category $\Cal C$ gives rise to a 3-dimensional
homotopy quantum field theory (HQFT) with target space $K(\pi,1)$. As  in the
standard case $\pi=1$, this HQFT has several   equivalent versions. We  
describe  here one of these versions formulated in terms of  extended
 $\pi$-manifolds.

Throughout Section 10 we fix
a
   modular  crossed $\pi$-category $\Cal C $ and 
   an  Eilenberg-MacLane CW-complex $X=K(\pi,1)$   with base point $x\in X$.

 \skipaline \noindent {\bf 10.1.   Extended   $\pi$-manifolds with
boundary.}  Recall the notion of a connected  extended 
  $\pi$-manifold  without boundary introduced in Section 7.5.
  In the present context it is more convenient to use
  (pointed) maps to $X=K(\pi,1)$ rather than   group
homomorphisms to $\pi$. 
Thus, instead of a group homomorphism   $\pi_1(M\backslash \Omega,z) \to \pi$
as in Section 7.5, we shall  consider a (pointed) map $g: M\backslash
\Omega   \to X$ inducing   this homomorphism. 
(In the definition of a coloring   we should then 
   replace $g$ by $g_{\#}$.) After this minor modification of the definitions of
Section 7.5 we can  extend them to manifolds with boundary.
 
Let $M$ be a compact connected oriented 3-manifold  with non-void pointed
boundary.   A  {\it ribbon  graph}
    $\Omega\subset M$ consists of a finite family of framed oriented
embedded arcs, circles and coupons (the {\it strata} of $\Omega$) such that

(i) the  strata of $\Omega$    are disjoint except that  
some endpoints of the arcs lie on the bases of the coupons;

(ii) all   other endpoints of the arcs  of $\Omega$ lie on $\partial M$ and form
the (finite) set   $\Omega \cap \partial M$; this set
 does not contain the base points of
$\partial M$; 

(iii) the framings of the strata form a framing of $\Omega$, i.e., a continuous
non-singular vector field on $\Omega$ transversal to $\Omega$; at  
$\Omega\cap \partial M$ the framing is tangent to $\partial M$;

(iv) on each coupon of $\Omega$ the framing is transversal to the coupon and yields
together with the orientation of the coupon the given orientation of $M$. 

Slightly pushing  $\Omega$ along its framing we obtain a
disjoint   copy $\tilde \Omega$ of $\Omega$.
Pushing a stratum  $t$  of $\Omega$  along the   framing we obtain a stratum 
  $\tilde
t$ of $ \tilde \Omega$. 

  A  {\it $\pi$-graph}  in  $M$ 
 is a ribbon graph $\Omega\subset M$ 
 endowed with a map
$g: M\backslash \Omega\to X$ which sends   the base points of all the
components of $\partial M$ into $x\in X$. 
 We consider $g$ up to homotopy constant on the base points
of $\partial M$. 

Let $ (\Omega ,g)$ be 
 a $\pi$-graph in $M$. A {\it coloring of $\Omega$ with respect to a base point
$z$ of a component of $ \partial M$} 
 comprises  two functions $u$ and $v$.  
 The     function  $u$ assigns to
every arc or circle $t$ of $\Omega$ and to every path $\gamma:[0,1] \to
M\backslash \Omega$ connecting $z$ to       $\gamma(1) \in \tilde t\subset \tilde
\Omega$ a certain  object $u_\gamma\in {\Cal C}_{g_{\#}(\mu_{\gamma})}$  where
$\mu_{\gamma}\in \pi_1(M\backslash \Omega, z)$ is the   meridian of $t$
determined by $\gamma$ and $g_\#:\pi_1(M\backslash \Omega,z)\to \pi$ 
is the homomorphism   induced
by $g$.   The   function  $v$ assigns to every coupon  $Q$ of $ \Omega$ and to every 
path $\gamma$ in $M\backslash \Omega$ connecting $z$
to      $\gamma(1) \in \tilde Q \subset \tilde
\Omega$  a certain  morphism $v_\gamma$ in $
{\Cal C}_{g_{\#}(\mu_{\gamma})}$.
The functions $u,v$ should satisfy the same conditions   as in Section
4.1 (where $g$ should be replaced with $g_{\#}$).

A coloring $(u,v)$ of $\Omega$ with respect to the base point
$z$ of a  component  of $  \partial M$ can be canonically transported  to a 
coloring   of $\Omega$   with respect to  the base point, $z'$,  of any other 
component  of $  \partial M$. (Here it is essential that $M$ is connected). 
Namely, choose
a path $\rho:[0,1]\to M\backslash \Omega$ with $\rho(0)=z, \rho(1)=z'$. 
For an arc or
circle $t$ of $\Omega$ and a path $\gamma:[0,1] \to M\backslash \Omega$
connecting $z'$ to      $\gamma(1) \in \tilde t $ set 
$$u_\gamma= (\varphi_{[g\circ \rho]})^{-1} (u_{\rho \gamma}).$$
Here the path $g\circ \rho $ is a {\it loop} in $(X,x)$
so that we can consider its homotopy class $[g\circ \rho]\in \pi$.
The path $\rho \gamma$ connects $z$ to  $\gamma(1) \in \tilde t $
so that $u_{\rho \gamma}$ is an  object of $ {\Cal
C}_{g_{\#}(\mu_{\rho \gamma})}$.
We have $\mu_{\rho \gamma}=\rho \mu_{\gamma} \rho^{-1}$ and so 
$g_{\#}(\mu_{\rho \gamma})= [g\circ \rho] g_{\#}(\mu_{  \gamma})
[g\circ \rho]^{-1}$. Hence  the functor
$(\varphi_{[g\circ \rho]})^{-1}$ maps  ${\Cal C}_{g_{\#}(\mu_{\rho \gamma})}$
into  ${\Cal
C}_{g_{\#}(\mu_{\gamma})}$ and  $u_\gamma$ is an object of 
${\Cal
C}_{g_{\#}(\mu_{\gamma})}$. Similarly, for a coupon
$Q$ of $\Omega$ and   a path
$\gamma:[0,1] \to M\backslash \Omega$ connecting $z'$ to      $\gamma(1) \in
\tilde Q $ set  $$v_\gamma= (\varphi_{[g\circ \rho]})^{-1} (v_{\rho \gamma}).$$
It is easy to check that these formulas define 
a  coloring  of $\Omega$ with respect to  
$z'$. It follows from definitions that this coloring does not depend on the
choice of $\rho$. 
Thus, a coloring   of $\Omega$ with respect to the base point
  of one component of $  \partial M$  canonically extends to  a system of colorings
of $\Omega$ with respect to  the base points of all components of $  \partial M$.
 This system is called a 
  {\it coloring} of  $\Omega$.  To specify a
coloring of $\Omega$ it is enough to specify 
a
coloring of $\Omega$  with  respect to one base
point, the colorings of $\Omega$   with  respect to all the other base points are
obtained by the transport as above.

  A tuple
  $(M,\Omega,g,u,v,\lambda)$   consisting of a  colored $\pi$-graph 
 $(\Omega,g,u,v)$ in $M$ and  a Lagrangian space  
$\lambda\subset H_1(\partial M; \bold R)$  is called 
a  {\it  connected extended 
  $\pi$-manifold with boundary}.

 Given a connected extended $\pi$-manifold   $(M,\Omega,g,u,v, \lambda)$, the
surface $\partial M$  becomes an extended $\pi$-surface as follows. By definition, 
$\partial M$ is   pointed. We provide $\partial M$  with orientation induced by the
one in $M$ (see Section 7.2 for our orientation convention). Set $P=\Omega \cap
\partial M$. Clearly, $P$ is a finite subset of $\partial M$ disjoint from the base
points. A point $p\in P$ is provided with sign $-$ if the adjacent arc  of $\Omega$ 
is oriented towards $p$ and with sign $+$ otherwise.  We provide $p$ with tangent
direction generated by the framing vector  at $p$ (by (iii) this vector is tangent to
$\partial M$). We can assume that $\tilde P=\tilde \Omega\cap \partial M$.
Denote by $\partial g$ the restriction of  
$g: M\backslash \Omega\to X$ to $\partial M\backslash P $. The
pair
$(P, \partial g)$ is    a $\pi$-marking on $\partial M$. We  
define its coloring $\partial u$  as follows.  Every path
$\gamma:[0,1] \to \partial M \backslash P$ leading from a base point   to $\tilde P$
is a path in $M \backslash \Omega$ leading to $\tilde \Omega$. 
Set $(\partial u)_\gamma=u_\gamma\in {\Cal C} $.
It is clear that the tuple  $(\partial M, P, \partial g, \partial u, \lambda)$
is an extended $\pi$-surface.  By definition, the boundary of an 
extended $\pi$-manifold without
boundary  is an empty surface.

Let $(M,\Omega,g,u,v,\lambda)$ and $(M',\Omega',g',u',v',\lambda')$ be two  
connected extended 
  $\pi$-manifolds with boundary.
A {\it weak  $e$-homeomorphism} $$
(M,\Omega,g,u,v,\lambda)\to (M',\Omega',g',u',v',\lambda')$$
   is a  degree $+1$  pointed 
homeomorphism of pairs $  f:(M,\Omega )\to (M',\Omega' )$
 preserving the framing, the orientation  and the splitting of $\Omega, \Omega'$
into strata and such that
  $g' f =g:   M\backslash \Omega  \to X$ and for any 
path $\gamma$ in $M\backslash \Omega$ leading from a base point of $\partial M$
to an arc or a circle of $\tilde \Omega$  (resp. a coupon of $\tilde \Omega$) we have
$u'_{f\circ \gamma}=u_{\gamma}$ (resp. $v'_{f\circ \gamma}=v_{\gamma}$). 
A weak  $e$-homeomorphism  is an {\it $e$-homeomorphism} if the induced
isomorphism $H_1(\partial M; \bold R)\to H_1(\partial M'; \bold R)$ maps
$\lambda$ onto $\lambda'$. It is clear that a (weak) 
$e$-homeomorphism 
of extended 
  $\pi$-manifolds induces a (weak) 
$e$-homeomorphism 
of their boundaries.

Taking disjoint unions of connected  extended $\pi$-manifolds with or without
boundary we obtain general   extended $\pi$-manifolds.
The notions of boundary, $e$-homeomorphism and weak $e$-homeomorphism
generalize to them in the obvious way. 

An  example of an extended $\pi$-manifold is provided
by the cylinder $\Upsilon\times
[0,1]$ over an extended $\pi$-surface  $\Upsilon=(\Upsilon, P\subset
\Upsilon,g:\Upsilon\backslash P\to X,u, \lambda)$. The cylinder   has  the 
following  structure of an extended $\pi$-manifold. 
We provide $\Upsilon\times
[0,1]$ with the product orientation where the interval $[0,1]$ is
oriented from left to right. We provide
$\partial (\Upsilon\times
[0,1])= (\Upsilon\times 0) \amalg (\Upsilon\times1)$ with base points
$z\times 0, z\times 1$ where   $z\in \Upsilon\backslash P$ 
runs over the base points of the components of $\Upsilon$. Set  
$\Omega=\cup_{p\in P}  \, (p\times [0,1])$ and provide each arc $p\times [0,1]$ with
constant framing determined by a  vector representing the given tangent direction
at $p$. We orient the arc $p\times [0,1]$ towards $p\times 1$ if the sign of $p$ is
$-$ and towards $p\times 0$ otherwise.  Let $\pr$ be  the projection $ 
(\Upsilon\times [0,1]) \backslash \Omega \to  \Upsilon\backslash P$.
Then the map $g\circ \pr$ makes $\Omega$ a $\pi$-graph and the formula
$\gamma\mapsto u_{\pr\circ \gamma}$  defines   its
coloring. Finally, 
we   provide
$ \Upsilon\times
[0,1] $ with  Lagrangian space in real 1-homology of the boundary equal to the direct
sum of
  copies of $\lambda  $ in $H_1(\Upsilon\times 0;\bold R)$ and $ 
H_1(\Upsilon\times 1;\bold R)$. 
This makes $\Upsilon\times
[0,1]$ an  extended $\pi$-manifold. Clearly, $\partial (\Upsilon\times
[0,1])= ((-\Upsilon)\times 0) \amalg (\Upsilon\times1)=
(-\Upsilon)  \amalg  \Upsilon$ where $=$ means the equality   of extended
$\pi$-surfaces.

 \skipaline \noindent {\bf 10.2.    The HQFT  $(\Cal T_{\Cal C},\tau_{\Cal C})$.}
The HQFT derived from a modular crossed
 $\pi$-category $\Cal C$ comprises
the 2-dimensional
homotopy modular functor $\Cal T=\Cal T_{\Cal C}$ discussed in Section 8 
and a function, $\tau=\tau_{\Cal C}$,   assigning to every 
extended 
 3-dimensional $\pi$-manifold $M $ 
a vector $\tau (M)\in \Cal T(\partial M)$. 
A construction of $\tau$ will be outlined in Section 10.4. We  state here main   
properties of $\tau$: 

(10.2.1)  if $f:M\to M'$ is an $e$-homeomorphism of  extended $\pi$-manifolds
then $(f\vert_{\partial M})_{\#}(\tau(M))= \tau(M')$ where  $(f\vert_{\partial
M})_{\#}:\Cal T(\partial M)\to \Cal T(\partial M')$ is the  isomorphism 
induced by the $e$-homeomorphism
$f\vert_{\partial M}:\partial M \to  \partial M'$;

(10.2.2)  for 
disjoint extended $\pi$-manifolds $M, M'$, we have $$\tau (M\amalg M')
=\tau(M)\otimes \tau(M') \in \Cal T(\partial M) \otimes_K  \Cal T (\partial M')=\Cal
T (\partial M\amalg \partial M');$$

(10.2.3) for an     extended $\pi$-manifold 
$M$ without boundary, the invariant 
$\tau (M )\in \Cal T(\emptyset)=K$ is the invariant  introduced
in Section 7.5;

(10.2.4) for   any extended $\pi$-surface $\Upsilon$,  the vector
 $ \tau(\Upsilon\times
[0,1]) \in \Cal T(-\Upsilon)\otimes_K  \Cal T( \Upsilon)$ 
 is the image of 
$d_\Upsilon\in  (\Cal T(\Upsilon)\otimes_K  \Cal T(-\Upsilon))^*$ under the
identifications  $$  (\Cal T(\Upsilon) \otimes_K \Cal
T(-\Upsilon))^*= (\Cal T(\Upsilon)  )^*
\otimes_K ( \Cal T(-\Upsilon))^*= \Cal T(-\Upsilon) \otimes_K \Cal
T( \Upsilon) $$
induced by   $d_\Upsilon$. 
 
A fundamental   property of $\tau$ is a  formula
which   computes $\tau(M) $ for an extended $\pi$-manifold $M$ obtained by gluing
two extended $\pi$-manifolds $M_1,M_2$ along   their
boundaries or  along several components of their
boundaries. In other words, if a closed surface in $M$ splits   $M$ into two pieces
$M_1,M_2$ then  the gluing formula computes $\tau(M)$  from
$\tau(M_1),\tau(M_2)$. By (10.2.2), it is enough to state this formula for
connected $M$.

We first consider  the case   $\partial M=\emptyset$. 
Let $M$ be  a closed connected  oriented  3-manifold
and $\Omega=(\Omega \subset M,z\in M\backslash \Omega, 
g:  M\backslash \Omega \to X, u,v)$ be a 
colored $\pi$-graph in $M$. Let
$\Sigma\subset M $ be a   closed non-void surface in $M$
splitting  $M$ into two compact submanifolds $M_1,M_2$ with 
$\partial M_1=\partial M_2=\Sigma$. We provide each $M_r$ with
orientation induced by the one in $M$.
Assume (deforming if necessary $\Sigma$ in
$M$) that 
   $\Sigma$ meets $\Omega$ transversally in a finite set of points lying on
the arcs and circles of $\Omega$. We fix a base point on each component of
$\Sigma$ so that these base points do not belong to   $\Sigma\cap \Omega$.   It is
clear that     $\Omega_r=\Omega\cap M_r$ is a ribbon graph
 in $M_r$ for $r=1,2$.  
Here the framing of $\Omega_r $ is the restriction of the framing of
$\Omega$.
We can thus assume that $\tilde \Omega_r=\tilde \Omega\cap M_r$. We upgrade 
$ \Omega_r$ to a colored $\pi$-graph  as follows.  
Deforming if necessary  $g: M\backslash \Omega \to X$ (keeping $g(z)=x)$ we
can assume that $g$  sends all the base points of $\Sigma$ into $x\in X$.
 Denote by
$g_r$ the restriction of $g$ to $M_r\backslash \Omega_r\subset  M\backslash
\Omega$. The pair $(\Omega_r, g_r)$ is   a $\pi$-graph in $M_r$. 
Transporting (as in Section 10.1) the given coloring $(u,v)$ of $\Omega$   along  
arbitrary paths in $M$ connecting $z$ to  the  base points of $\Sigma$  we obtain
colorings of $\Omega$ with respect to  these     points.
(These colorings do not depend on the choice of the paths). 
Now restricting these colorings to the
paths $\gamma$  in  $M_r\backslash \Omega_r $
 connecting the   base points  of $\Sigma=\partial M_r$ to $\tilde \Omega_r$ we
obtain a coloring  $(u_r,v_r)$  of $(\Omega_r, g_r)$. 
 Fix 
 a Lagrangian space  
$\lambda_0\subset H_1(\Sigma; \bold R)$. Then the tuple $(M_r, \Omega_r, g_r,
u_r,v_r,\lambda_0)$ is  an extended 
  $\pi$-manifold briefly denoted  $M_r$.
Its boundary  $\partial M_r$   is an extended $\pi$-surface   with 
  underlying  surface  
 $\Sigma$. It follows from definitions that
$\partial M_2=- \partial M_1$ where minus denotes the negation of extended
$\pi$-surfaces. We have the vectors 
$\tau(M_1 )\in \Cal T(\partial M_1)$
and $\tau(M_2 )\in \Cal T(\partial M_2)=\Cal T(-\partial M_1)$.
Then $$\tau(M )= ({\Cal D} \Delta_{-}^{-1})^n\, d_{\partial M_1}
(\tau(M_1 )\otimes \tau(M_2 )) \in K\leqno (10.2.a) $$
 where  $M$ denotes the extended $\pi$-manifold 
$(M,\Omega ,z,  g, u,v)$,
  $d_{\partial M_1}$ is the pairing 
$\Cal T(\partial M_1) \otimes_K \Cal T(-\partial M_1)\to K$
provided by (8.3.5) and 
 $n$ is an integer computed as follows.
Consider the inclusion homomorphism $H_1(\partial M_r;\bold R)
\to
H_1(  M_r;\bold R)$   and observe that its kernel,  $\lambda_r$, is 
a  Lagrangian subspace in $H_1(\partial M_1;\bold R)=H_1(\partial M_2;\bold R)
=H_1(\Sigma;\bold R)$. Then $n=\mu ( \lambda_1,\lambda_0, \lambda_2) \in
\bold Z $
 is the Maslov index  of the Lagrangian spaces $\lambda_1,\lambda_0, \lambda_2$
in  $H_1(\partial M_1;\bold R)$.

 Let now   $M=(M,\Omega,g:M\backslash \Omega \to X,u,v,\lambda)$  be a connected  
extended 
  $\pi$-manifold with boundary. Let
$\Sigma\subset M\backslash \partial M $ be a closed non-void surface
  splitting  $M$ into two compact submanifolds $M_1,M_2$. For   $r=1,2$, we have 
  $\partial M_r=\Sigma \amalg   (\partial M \cap   M_r)$   
where  $\partial M \cap   M_r=\partial M  \cap \partial M_r$ is a disjoint union
(possibly void) of certain components of $\partial M$ such that $\partial M= 
(\partial M \cap   M_1) \amalg (\partial M \cap   M_2)$.  
Assume that $\lambda=\lambda^1\oplus \lambda^2$ where 
$\lambda^r$ is 
a  Lagrangian  space in $H_1( \partial M \cap   M_r; \bold R)$.
Then the extended $\pi$-surface $\partial M$ splits as the disjoint union 
of  extended $\pi$-surfaces 
$\Upsilon^1 =(\partial M \cap   M_1,\lambda^1)$ and 
$\Upsilon^2 =(\partial M \cap   M_2,\lambda^2)$.

Assume    that 
   $\Sigma$ meets $\Omega$ transversally in a finite set of points lying on
the arcs and circles of $\Omega$.
 We fix a base point on each component of
$\Sigma$ so that these base points do not belong to   $\Sigma\cap \Omega$.
Deforming if necessary  $g: M\backslash \Omega
\to X$ (keeping $g$ on $\partial M$) we can assume that $g$  sends all the base
points of $\Sigma$ into $x\in X$.
As above,  $\Omega_r=\Omega\cap M_r$ is a ribbon graph
 in $M_r$ for $r=1,2$ where the orientation in
$M_r$ is induced by the one in $M$ and the framing on $\Omega_r$ is the
restriction of the one on $\Omega$.  
 Denote by
$g_r$ the restriction of $g$ to $M_r\backslash \Omega_r\subset  M\backslash
\Omega$. The pair $(\Omega_r, g_r)$ is   a $\pi$-graph in $M_r$. 
Choose a    component of $ \partial M$ and transport   the
given coloring $(u,v)$ of $\Omega$  at its  base point    along  arbitrary paths in
$M$
  to  the  base points of $\Sigma$. This gives   colorings of $\Omega$
with respect to  the  base points of $\Sigma$. These colorings  depends neither on 
 the choice of the initial 
component of $ \partial M$ nor on the choice of the paths in $M$. Now restricting
the   colorings of $\Omega$ (with respect to the base points of
$\Sigma$ and $   \partial M \cap   M_r$) to the paths
$\gamma$  in  $M_r\backslash \Omega_r $
  we obtain a coloring $(u_r,v_r)$  of $(\Omega_r, g_r)$. 
Fix 
 a Lagrangian space  
$\lambda_0\subset H_1(\Sigma; \bold R)$. 
  Then the tuple $(M_r, \Omega_r, g_r,
u_r,v_r,\lambda_0\oplus \lambda^r)$ is  an extended 
  $\pi$-manifold   briefly denoted  $M_r$.
Its boundary $\partial M_r$ is a disjoint
union of $\Upsilon^r$ and an extended $\pi$-surface  
$\Upsilon_r  $ with
  underlying   surface   $\Sigma$. Clearly, $\Upsilon_2=-\Upsilon_1$. 
For $r=1,2$, we have  a vector 
$$\tau(M_r )\in \Cal T(\partial M_r)= \Cal T (\Upsilon^r)\otimes_K  
\Cal T (\Upsilon_r).$$ 
Then $$\tau(M )= ({\Cal D} \Delta_{-}^{-1})^n \,d_{\Upsilon_1 }
(\tau(M_1 )\otimes \tau(M_2 ))
\in \Cal T (\Upsilon^1)\otimes_K 
\Cal T (\Upsilon^2) \leqno (10.2.b) $$
where $d_{\Upsilon_1}$ is the pairing 
$\Cal T(\Upsilon_1) \otimes_K \Cal T(\Upsilon_2)\to K$
provided by (8.3.5)
and $n$ is an integer computed as follows.
Let   $\lambda_r$ be the subset of $H_1(\Upsilon_r;\bold R)
=H_1(\Sigma;\bold R)$
consisting of   elements  homological in $M_r$ to 
elements of $\lambda^r\subset H_1(\partial M \cap   M_r;\bold R)$. 
 One can check that
$\lambda_r$ is a Lagrangian  space.  
Then
 $n=\mu ( \lambda_1,\lambda_0, \lambda_2) \in \bold Z $
is the Maslov index  of     $\lambda_1,\lambda_0, \lambda_2$
in  $H_1(\Upsilon^1;\bold R)$.

 \skipaline \noindent {\bf 10.3.   Extended    $\pi$-cobordisms.}
As   in the study of TQFT's, we can reformulate the vector $\tau_{\Cal
C}$ as an operator invariant  of   $\pi$-cobordisms. An {\it extended 
(3-dimensional) $\pi$-cobordism} is a triple $(M, \partial_- M, \partial_+ M)$
where $M$ is a  3-dimensional extended  $\pi$-manifold  and  $\partial_- M,
\partial_+ M$ are extended $\pi$-surfaces  such that $\partial M=(-\partial_- M)
\amalg  \partial_+ M$. We call  $\partial_- M, \partial_+ M$
  the {\it bottom base} and {\it top base} of the cobordism,
respectively. Clearly, $$\Cal T( \partial M)=  \Cal T(
-\partial_- M)\otimes_K \Cal T( \partial_+ M) =  (\Cal T( \partial_-
M))^*\otimes_K \Cal T( \partial_+ M) $$
$$
=\Hom_K (\Cal T( \partial_- M), \Cal T( \partial_+ M)).$$
 By these identifications,  the vector $\tau(M)\in \Cal T( \partial M)$ determines 
  a
homomorphism $ \Cal T( \partial_- M)\to \Cal T( \partial_+ M)$
denoted $\tau (M, \partial_-
M, \partial_+ M)$. Note that this definition   implicitly involves the pairing
$d_{\partial_{-}  M} $. 
 By axioms (10.2.1) and (10.2.2), the homomorphism $\tau (M, \partial_-
M, \partial_+ M)$ is natural with respect to $e$-homeomorphisms of cobordisms
and $\otimes$-multiplicative with respect to disjoint unions.

For example, if $\partial_-
M=\emptyset$ then $\partial_+ M=\partial  M$ and 
$\tau (M, \partial_-
M, \partial_+ M):\Cal T( \emptyset)=K\to \Cal T( \partial M)$
sends any $k\in K$ into $k \,\tau(M)$. 
If $\partial_+ M=\emptyset$ then
$\partial_- M=-\partial  M$ and 
$\tau (M, \partial_-
M, \partial_+ M):\Cal T(-\partial  M) \to \Cal T( \emptyset)=K$
sends any $h\in \Cal T(-\partial  M)$ into $d_{\partial  M} (\tau(M),h)$.
By
  (10.2.3),  if $\partial  M=\emptyset$ then  
$\tau (M, \partial_-
M, \partial_+ M):K\to K$
is multiplication by  the invariant 
$\tau (M )\in  K$   introduced
in Section 7.5. 

We can reformulate  formulas (10.2.a) and (10.2.b)
as $$\tau (M, \partial_-
M, \partial_+ M)= ({\Cal D} \Delta_{-}^{-1})^n \tau (M_2, \partial_-
(M_2), \partial_+ (M_2))
\circ \tau (M_1, \partial_-
(M_1), \partial_+ (M_1))$$
where we view $M, M_1,M_2$  as extended 
$\pi$-cobordisms with  bottom bases $-(\partial M \cap   M_1),
-(\partial M \cap   M_1)  ,-\Upsilon_2 =\Upsilon_1$
and top bases $\partial M \cap   M_2,
\Upsilon_1 , \partial M \cap   M_2$, respectively.

For an extended $\pi$-surface $\Upsilon=(\Upsilon, P\subset
\Upsilon,g:\Upsilon\backslash P\to X,u, \lambda)$,  
   the triple  $ 
(\Upsilon\times [0,1], \Upsilon\times 0, \Upsilon\times 1)$ is 
an extended   $\pi$-cobordism whose bases are copies of $\Upsilon$. 
By   (10.2.4), 
$$\tau(\Upsilon\times
[0,1], \Upsilon\times 0, \Upsilon\times 1)=\id_{\Cal T( \Upsilon)}
: T(\Upsilon)\to  \Cal T( \Upsilon).$$
For  every $\alpha\in \pi$, we   define a twisted cylinder 
$(\Upsilon\times [0,1])^\alpha$ as follows. It differs from 
$\Upsilon\times [0,1]$ only by the choice of the map 
$(\Upsilon\times [0,1]) \backslash \Omega \to X$ and the choice of the coloring
of the ribbon graph $\Omega=P\times [0,1]$.
The map in question is chosen so that its restriction to the bottom base equals $g$
and its restrictions to   arcs $z\times [0,1]$  are loops representing 
$\alpha^{-1}$ for all the base points   $z\in \Upsilon\backslash P$ 
  of the components of $\Upsilon$.
The coloring of  $\Omega$ is chosen so that its restriction to the bottom base
coincides with $u$. It is clear that $(\Upsilon\times [0,1])^\alpha$ 
is 
an extended   $\pi$-cobordism with bottom base  
$\Upsilon $ and top base  ${}^\alpha\!\Upsilon$. The  operator invariant of
this cobordism   defines the action of $\alpha$
discussed in Section  8.5:
$$\alpha_*=\tau ((\Upsilon\times [0,1])^\alpha, \Upsilon, {}^\alpha\!\Upsilon): \Cal T(\Upsilon)\to \Cal
T({}^\alpha\!\Upsilon).$$

 \skipaline \noindent {\bf 10.4.   Construction of   $(\Cal T_{\Cal C},\tau_{\Cal
C})$.} The construction closely follows the known construction of   3-dimensional
TQFT's from modular categories, see [Tu2]. First, one defines the operator invariant
for  extended $\pi$-cobordisms 
whose boundary components are parametrized, i.e., identified with
  standard $\pi$-surfaces in $S^3$. (Only the geometric position of the surfaces in
$S^3$ is   standard, the map to $X=K(\pi,1)$ is arbitrary). Then one uses these
operators to define the action of weak $e$-homeomorphisms. Finally, one
replaces the parametrizations with Lagrangian spaces in homology. 
We skip the details and  give here only one of the key  lemmas whose proof is
somewhat different from the standard case.

 \skipaline \noindent {\bf 10.5.  Lemma.  }
{\sl Let $T \subset {\bold R}^2\times [0,1]$ be a tangle 
consisting of a vertical interval $t$ oriented downwards and its meridian $m$,
both with zero framing.  Sending all meridians of $t$ into $1\in \pi$ and all
meridians of $m$ into  $\alpha\in \pi$,
we make $T$ a $\pi$-tangle, $T_\alpha$. 
 Let   $V\in {\Cal C}_1$ be a simple object in the neutral
component 
  of  ${\Cal C}$. Let   $T_\alpha (V)$ be
   $T_\alpha$ colored so that
$m$ acquires the canonical color as in Section 7.2 and the  
target  of   $T_\alpha$ is the triple 
$(+1, 1\in \pi, V)$.
Then the
source   of   $T_\alpha (V)$ is the triple  $(+1, 1\in \pi, \varphi_{\alpha^{-1}}(V))$
and} $ F(T_\alpha (V))\in \Hom_{\Cal C} (\varphi_{\alpha^{-1}}(V), V)$ {\sl  is
computed by} $$ F(T_\alpha (V))=\cases {\Cal D}^2\, \id_{\1},~ { {if}}\,\,\, 
$V=\1$, \\ 0,~  { {if}}\,\,\, V \,\, { \text  {is not isomorphic to}} \,\,\1.\endcases
$$ 

\skipaline {\sl Proof.}  We begin with another  useful  identity. 
Fix $\beta\in
\pi$ and  an object $W\in {\Cal C}_\beta$.  Consider the $\pi$-tangle
$T_\beta$ as in the statement of the  lemma  with  $\alpha$ replaced by $\beta$.
We    present $T_\beta$ by a plane diagram with   two crossings. We 
 attach $W$ to the
  arc  of the diagram representing $m$ and attach $V$
 (resp. $\varphi_{\beta^{-1}} (V)$) to the arc
incident to the   output (resp. intput).  
Denote the resulting  
colored $\pi$-tangle by
$T^V_W$. Its   target  and source are the triples  $(+1,
1\in \pi, V)$ and  $(+1,
1\in \pi, \varphi_{\beta^{-1}} (V))$, respectively. 
We claim that
$$F(T^V_W) \, F(T_\alpha (\varphi_{\beta^{-1}} (V)))=
\dim (W) \, F(T_{\beta\alpha} (V )):\varphi_{\alpha^{-1} \beta^{-1} } (V) \to
V.\leqno (10.5.a)$$
We first prove this equality   and then deduce from it the claim of
the lemma.  
Glueing $T^V_W$ on the top of 
$T_\alpha (\varphi_{\beta^{-1}} (V))$ we obtain a colored $\pi$-tangle  $$\tilde
T= T^V_W \circ T_\alpha (\varphi_{\beta^{-1}} (V)).$$ Geometrically, 
$\tilde T$
consists of a vertical interval with two   meridians and all framings
zero. Clearly,
$$F(T^V_W) \, F(T_\alpha (\varphi_{\beta^{-1}} (V)))=
F(\tilde T).$$ 
Denote by $\ell^{\beta}_W  $  a colored $\pi$-knot  (or rather unknot) in 
  ${\bold
R}^2\times [0,1]$ represented by a  plane circle labeled with $(\beta, W)$.
Sliding the $W$-colored circle of $\tilde T$ along the $can$-colored
circle we can transform $\tilde T$  into a
disjoint union  $\ell^{\beta}_{W} \amalg T_{\beta\alpha} (V)$. (An explicit splitting
of this handle sliding into a 
 composition of Kirby-Fenn-Rourke moves is given in [Tu2, p. 93]). 
As in the proof of Theorem 7.3, we obtain 
$$F(\tilde T)=F(\ell^{\beta}_{W} \amalg T_{\beta\alpha} (V))= F(\ell^{\beta}_{W})\,
F(T_{\beta\alpha} (V)) =\dim (W) \, F(T_{\beta\alpha} (V)).$$
This implies (10.5.a).

Now we can prove the claim of the lemma. If $V=\1$ then  
 formula (2.2.f)   implies that we can push the circle stratum   of $T$
across the  arc  stratum without changing the operator invariant.
This yields our claim in this case.
Assume that $V$ is not isomorphic to $\1$. 
 Let $I$ be the set of   isomorphism 
classes of simple objects in ${\Cal C}_1$ and let $\{V_i\in {\Cal C}_1\}_{i\in I}$ be
representatives of these classes. We can assume that
$V=V_i$  for a certain $i\in I$.  For  $W=V_j$ with $j\in I$ we
can compute $F(T^V_W)\in \End(V)$ 
explicitly. Since $V$ is simple,  $F(T^V_W)=k\, \id_V$ with $k\in K$. The
closure of $T^V_W$  is the Hopf link whose components are colored with
$V=V_i, W=V_j$. Therefore $k \,\dim   (i)=S_{i,j}$ where 
$\dim   (i)=\dim  (V)= \dim (V_i)$.
Hence $k=  (\dim   (i))^{-1} S_{i,j}$ where we use the invertibility of $\dim (i)$ (Lemma
6.5). Substituting this in (10.5.a)  (for $\beta=1$) we
obtain $$ (\dim   (i))^{-1} S_{i,j}   F(T_{\alpha} (V))=   \dim (j )\,
F(T_{\alpha} (V)). $$
 By [Tu2, formula (3.8.b)],  $\sum_{j\in I} \dim (j) S_{i,j}  = 0$.
Hence
$$   F(T_{\alpha} (V)) = {\Cal D}^{-2} 
\sum_{j\in I} \, \dim (j) \, \dim (j)  \, F(T_{\alpha} (V))$$
$$
=  {\Cal D}^{-2} \,(\dim   (i))^{-1} \sum_{j\in I}  \dim (j)\, S_{i,j} \,  F(T_{\alpha}
(V)) = 0.$$

\skipaline \centerline {\bf 11.   Hopf group-coalgebras}

\skipaline \noindent {\bf 11.1.     Hopf  algebras.} For convenience of the
reader we   recall the standard definitions of quasitriangular and
ribbon Hopf algebras,   see for instance [KRT], [Tu2]. A  {\it Hopf algebra}
over $K$ is a tuple $(A, \Delta, \varepsilon, s)$
where $A$ is an associative unital algebra over $K$, 
$\Delta :  A  \to A^{\otimes 2}=A \otimes A$  and 
$\varepsilon:A\to K$ are algebra homomorphisms, $s:A\to A$ is an algebra
anti-homomorphism such that $\Delta$ is coassociative and 
$$( \id_{A} \otimes \varepsilon)
\Delta = ( \varepsilon \otimes \id_{A})
\Delta =\id_{A} :A\to A,$$ $$ \mu 
(s  \otimes \id_{A}) \Delta 
=\mu  ( \id_{A} \otimes s ) \Delta  = 1_A\, \varepsilon:A \to
A $$ where $ \mu $  is
multiplication in $A$. 

Let $A $ be a Hopf algebra over $K$. Denote  the flip (permutation) $A^{\otimes 2} \to
A^{\otimes 2}$   by $\Perm$. A
pair $(A,R\in A^{\otimes 2})$ is called a {\it quasitriangular Hopf algebra} if $R$ is
invertible in $ A^{\otimes 2}$ and  satisfies the following identities: $$\Perm
(\Delta(a))= R \Delta(a) R^{-1},$$ for all $a\in A$  and 
$$(\id_A \otimes \Delta)(R)=R_{13} R_{12}, \,\,\,\, (\Delta\otimes \id_A)
(R)=R_{13} R_{23}$$ where $R_{12}=R\otimes 1_A\in A^{\otimes 3} $,
$R_{23}=1_A\otimes R\in A^{\otimes 3} $
and 
$R_{13}=(\id_A \otimes \Perm)(R_{12})\in  A^{\otimes 3} $.
These identities imply the Yang-Baxter equality
$R_{12} R_{13} R_{23}=  R_{23} R_{13} R_{12}$.

Let $(A,R)$ be a quasitriangular Hopf algebra.  A triple
$(A,R, v\in A)$ is called a {\it ribbon Hopf algebra} if $v$ is an invertible element
of the center of $A$ such that
$s(v)=v$ and $\Delta(v)= (v\otimes v) \Perm(R) R$. 
If $A=(A,R, v)$ is   a {  ribbon Hopf algebra} then the tuple $(A,(\Perm
(R))^{-1},v^{-1})$ is also a ribbon Hopf algebra.  

It is well-known that the category
of representations of a quasitriangular (resp. ribbon) Hopf algebra is a  braided
(resp. ribbon) 
monoidal category.

\skipaline \noindent {\bf 11.2.    Hopf $\pi$-coalgebras.} Let $\pi$ be a
group.  The notion of a
$\pi$-coalgebra is dual to the notion of a (unital) $\pi$-graded algebra. 
 By a {\it    $\pi$-coalgebra} over   $K$, we mean the following data:

- a family of  $K$-modules $\{A_{\alpha}\}_{\alpha\in \pi}$; 

- a family of  $K$-linear  homomorphisms (called the comultiplication)
$$\Delta=\{\Delta_{\alpha,\beta
}:  A_{\alpha \beta} \to A_{\alpha}\otimes A_{\beta}\}_{\alpha,\beta \in \pi} ;
$$

- a $K$-homomorphism $\varepsilon_1:A_1\to K$ (the counit).

This data should satisfy the following two axioms:

(11.2.1) $ \Delta $ is coassociative in the sense that for any $\alpha,
\beta,\gamma \in \pi$, $$ (\Delta_{\alpha,\beta }\otimes \id_{A_\gamma})
\Delta_{\alpha \beta,\gamma } = (\id_{A_\alpha} \otimes \Delta_{ \beta,\gamma 
})
\Delta_{\alpha, \beta \gamma }:  A_{\alpha \beta\gamma }\to A_{\alpha }\otimes
A_{\beta } \otimes A_{\gamma } ;$$

(11.2.2) for any $\alpha \in \pi$, $$( \id_{A_\alpha} \otimes \varepsilon_1)
\Delta_{\alpha,1}= ( \varepsilon_1 \otimes \id_{A_\alpha})
\Delta_{1,\alpha}=\id_{A_\alpha} :A_\alpha\to A_\alpha$$

A {\it  Hopf $\pi$-coalgebra} over   $K$ is a $\pi$-coalgebra
$(A,\Delta,\varepsilon)$ where  each $A_\alpha$ 
is an associative $K$-algebra with multiplication $ \mu_{\alpha}$ and (left and
right) unit  $1_\alpha$ endowed with  
  algebra anti-isomorphisms $s= \{s_{ \alpha }:  A_{\alpha } \to 
A_{\alpha^{-1}}\}_{\alpha \in \pi}$ (the antipode) such that

(11.2.3) for all $\alpha,\beta\in \pi$ the comultiplication $\Delta_{\alpha,\beta
}$ is an algebra homomorphism and  $\Delta_{\alpha,\beta
}(1_{\alpha\beta})=1_\alpha\otimes 1_\beta$;

(11.2.4) the counit 
$\varepsilon_1:A_1\to K$ is an algebra homomorphism and 
$\varepsilon_1(1_1)=1_K$;

(11.2.5) for any $\alpha \in
\pi$,   $$ \mu_{\alpha} (s_{\alpha^{-1}} \otimes \id_{A_\alpha})
\Delta_{\alpha^{-1},\alpha} =\mu_{\alpha} ( \id_{A_\alpha} \otimes
s_{\alpha^{-1}}) \Delta_{\alpha, \alpha^{-1}} = 1_\alpha\, \varepsilon_1:A_1\to
A_\alpha. $$

It follows from the properties of   algebra  anti-isomorphisms  that
$s_\alpha(1_\alpha)=1_{\alpha^{-1}}$  for all $\alpha \in \pi$.  
Note
also that the tuple $(A_1,\Delta_{1,1}, \varepsilon_1, s_1)$ is a Hopf algebra in
the usual sense of the word. We call it the {\sl neutral component} of $A$. 

Warning:  
the notion of a Hopf $\pi$-coalgebra  is not self-dual.
It is quite interesting to study the dual notion of a Hopf $\pi$-algebra but
we shall not need it in this paper.  

 A {\it 
crossed
Hopf $\pi$-coalgebra} over   $K$ is  a
Hopf $\pi$-coalgebra  $(\{A_{\alpha}\}_{\alpha\in \pi},
\Delta , \varepsilon_1, s)$  
endowed with 
 a family of algebra isomorphisms $\varphi=\{\varphi_{ \alpha }:  A_{\beta } 
\to
 A_{\alpha \beta \alpha^{-1}}\}_{\alpha, \beta \in \pi}$ such that

(11.2.6) each $\varphi_{ \alpha }$ preserves the counit, the antipode, and the
comultiplication, i.e., for any $\alpha, \beta, \gamma\in \pi$, we have $$
\varepsilon_1 \varphi_{ \alpha }\vert_{A_1}=\varepsilon_1,$$ $$\varphi_{ \alpha }
s_{\beta}=s_{\alpha\beta\alpha^{-1}}\varphi_{\alpha }:A_\beta\to
A_{\alpha\beta^{-1}\alpha^{-1}},$$ $$(\varphi_{ \alpha }\otimes \varphi_{ \alpha 
})
\Delta_{\beta,\gamma} =\Delta_{\alpha\beta\alpha^{-1},
\alpha\gamma\alpha^{-1}}\,\varphi_{ \alpha }:  A_{\beta\gamma} \to
A_{\alpha\beta\alpha^{-1}} \otimes A_{ \alpha\gamma\alpha^{-1}};$$

(11.2.7) $\varphi $ is an action of $\pi$, i.e., $\varphi_{ 
\alpha\alpha'}=\varphi_{
\alpha}\varphi_{ \alpha'}$ for all $\alpha,\alpha' \in \pi$.

It is clear that
 $\varphi_{
\alpha}(1_\beta)=1_{\alpha\beta\alpha^{-1}}$ for all $\alpha, \beta\in \pi$
and $\varphi_{\alpha}(A_1)=A_1$. 
Restricting the action $\varphi $   of $\pi$ to   $A_1$ 
we obtain an action of $\pi$ on $A_1$ by Hopf algebra endomorphisms.
  
We end this subsection with  two examples of crossed Hopf $\pi$-coalgebras. Both
examples are 
derived from  an action of $\pi$
on a  Hopf
algebra $(A, \Delta, \varepsilon, s)$    over $K$  
  by Hopf algebra endomorphisms.  Set  $A^{\pi}=\{A_\alpha\}_{\alpha\in
\pi}$ where for  each $\alpha\in \pi$, the algebra $A_\alpha$ is a copy of $A$.
Fix an  identification isomorphism of algebras 
  $i_\alpha: A\to A_\alpha$.   For $\alpha, \beta\in \pi$, we define a
comultiplication $\Delta_{\alpha,\beta}:  A_{\alpha \beta} \to A_{\alpha}\otimes
A_{\beta}$   by  $$\Delta_{\alpha,\beta} (i_{\alpha\beta} (a))= \sum_{(a)}
i_\alpha (a' ) \otimes i_\beta (a'' )  $$ where $a \in A$ and    $\Delta  (a )=
\sum_{(a)}  a' \otimes  a''$ is the given    comultiplication in $A$ written in
Sweedler's sigma notation. The
counit $  \varepsilon_1 :A_1\to K$ is defined by $\varepsilon_1 (i_1(a))=
\varepsilon (a)\in K$ for $a\in A$. For ${\alpha  \in \pi}$,
the antipode  $s_{ \alpha }:  A_{\alpha } \to 
A_{
\alpha^{-1}}$ is given by $$s_{ \alpha } (i_\alpha(a ))= i_{\alpha^{-1}}( s(a)
) $$ where $a\in A$.  For ${\alpha, \beta \in \pi}$, the homomorphism
$\varphi_{ \alpha }:  A_{\beta } 
\to
 A_{\alpha \beta \alpha^{-1}}$ is defined by
$\varphi_{ \alpha } (i_\beta (a ))= i_{\alpha \beta \alpha^{-1}}  (\alpha(a))$.
All the axioms of a crossed Hopf $\pi$-coalgebra  for $A^{\pi}$ follow directly from  
definitions. 

The second example differs only by the definition of the comultiplication and
the antipode. Let $\overline A^{\pi}$ be the same family of algebras
$\{A_\alpha=A\}_{\alpha\in \pi}$ with the same counit,  the same action
$\varphi$ of $\pi$,   the
comultiplication $\overline \Delta_{\alpha,\beta}:  A_{\alpha \beta} \to
A_{\alpha}\otimes A_{\beta}$  and the antipode  $\overline s_{ \alpha }:  A_{\alpha
} \to  A_{
\alpha^{-1}}$ defined by  $$\overline \Delta_{\alpha,\beta} (i_{\alpha\beta} (a))=
\sum_{(a)} i_\alpha (\beta (a') ) \otimes i_\beta (a'' ),  $$ 
  $$\overline s_{ \alpha }
(i_\alpha(a ))= i_{\alpha^{-1}}(\alpha (s(a)) )=i_{\alpha^{-1}}(s(\alpha ( a))
)$$ where $a\in A$. The axioms of a crossed Hopf $\pi$-coalgebra  for 
$\overline A^{\pi}$ follow   from  
definitions.   Both   $A^{\pi}$ and $\overline A^{\pi}$ are extensions of $A$ since
$A^{\pi}_1=\overline A^{\pi}_1=A_1$ as Hopf algebras. 

In particular, if $G$ is a Lie group with Lie algebra
$\,\bold     g$ then      the universal enveloping
algebra $U(\,\bold     g)$ has a canonical  structure of a Hopf algebra 
 and  $G$ acts on  $U(\,\bold     g)$ by Hopf algebra endomorphisms induced by   the
group conjugation.  The  constructions above give     crossed Hopf $G$-algebras
 $({U}(\,\bold     g))^G= \{U(\,\bold     g)_\alpha\}_{\alpha\in G}$ 
and $\overline {({U}(\,\bold     g))}^G= \{U(\,\bold     g)_\alpha\}_{\alpha\in G}$ where    each
$U(\,\bold     g)_\alpha$ is a copy of $U(\,\bold     g)$ sitting at  $\alpha\in G$.

\skipaline \noindent {\bf 11.3.  Quasitriangular   Hopf $\pi$-coalgebras.}
Let $A=(\{A_\alpha\}, \Delta , \varepsilon_1,s,\varphi)$ be a crossed Hopf
$\pi$-coalgebra.  A {\it universal $R$-matrix} in $A$ is a family of invertible
elements $$R=\{R_{\alpha,\beta } \in A_{\alpha}\otimes 
A_{\beta}\}_{\alpha,\beta
\in \pi} \leqno (11.3.a) $$ satisfying the following conditions:

(11.3.1) for any $\alpha,\beta \in \pi, a\in A_{\alpha\beta}$, $$R_{\alpha,\beta 
}
\,\Delta_{\alpha,\beta } (a)= 
\Perm_{ \beta,\alpha } ( (
\varphi_{\alpha^{-1}}\otimes \id_{A_\alpha} ) \Delta_{\alpha \beta
\alpha^{-1},\alpha } (a))
 \,R_{\alpha,\beta } $$
where $\Perm_{ \beta,\alpha }$ is the flip $A_{\beta}\otimes A_{\alpha}\to
A_{\alpha}\otimes A_{\beta}$;

(11.3.2) for any $\alpha,\beta, \gamma \in \pi$, 
$$ ( \id_{A_\alpha} \otimes
\Delta_{\beta,\gamma }) (R_{\alpha, \beta \gamma })= (R_{ \alpha, \gamma })_{1
\beta 3} \,(R_{ \alpha, \beta })_{12\gamma} $$
 and 
$$ ( \Delta_{\alpha,\beta
}\otimes \id_{A_\gamma}) (R_{\alpha \beta,\gamma })
=
 ((\varphi_{\beta} \otimes
 \id_{A_\gamma}) (R_{\beta^{-1} \alpha \beta,  \gamma}))_{1\beta 3} \, (R_{
\beta, \gamma })_{\alpha 23} $$
 where for $K$-modules $P,Q$ and $r=\sum_j
p_j\otimes q_j \in P\otimes Q$ we set $$ r_{12\gamma} =r \otimes 1_\gamma \in
P\otimes Q \otimes A_\gamma,\,\,\,\,\, r_{\alpha 23} =1_\alpha \otimes r\in
A_\alpha \otimes P \otimes Q $$ and $$r_{1 \beta 3}=\sum_j p_j\otimes 1_\beta
\otimes q_j \in P\otimes A_\beta \otimes Q;$$

(11.3.3) the family (11.3.a) is invariant under the endomorphisms $\varphi_{ 
\alpha
}$, i.e., $$ (\varphi_{\alpha } \otimes \varphi_{\alpha }) (R_{ \beta, \gamma 
})
=R_{\alpha \beta\alpha^{-1},\alpha \gamma\alpha^{-1}}.$$

A crossed Hopf $\pi$-coalgebra endowed with a universal $R$-matrix is said to be
{\it quasitriangular}. It is easy to deduce from  (11.3.1) - (11.3.3)  the
Yang-Baxter equality for $R$:
$$(R_{ \alpha, \beta })_{12\gamma}
\, ((\varphi_{\beta} \otimes
 \id_{A_\gamma}) (R_{\beta^{-1} \alpha \beta,  \gamma}))_{1\beta 3}
\, (R_{
\beta, \gamma })_{\alpha 23} 
=  (R_{
\beta, \gamma })_{\alpha 23}  \,
 (R_{ \alpha, \gamma })_{1
\beta 3} \,(R_{ \alpha, \beta })_{12\gamma}. $$

\skipaline \noindent {\bf 11.4.  Ribbon  Hopf $\pi$-coalgebras.}  Let $A $ 
be
a quasitriangular (crossed) Hopf $\pi$-coalgebra with universal $R$-matrix (11.3.a).
A {\it twist} in $A$ is a collection of invertible elements $\{\theta_\alpha\in
A_\alpha\}_{\alpha\in \pi}$ such that

(11.4.1) $\varphi_{ \alpha }(a)=\theta_{\alpha}^{-1} a \theta_{\alpha} $ for all 
$
\alpha\in \pi, a\in A_\alpha$;

(11.4.2) $s_\alpha(\theta_\alpha)=\theta_{ \alpha^{-1}}$ for all $\alpha \in 
\pi$;

(11.4.3) for all $\alpha,\beta \in \pi$, 
$$\Delta_{\alpha,\beta}
(\theta_{\alpha\beta})= 
(\theta_\alpha \otimes 
\theta_\beta)\,\Perm_{\beta,\alpha
} ( (\id_{A_\beta} \otimes \varphi_{\alpha }  ) R_{  \beta,\alpha } )\, R_{ \alpha,
\beta } ;$$

(11.4.4) $\varphi_{ \alpha }(\theta_\beta)=\theta_{\alpha \beta \alpha^{-1}}$ for
 all $\alpha, \beta \in \pi$.

A quasitriangular  crossed   Hopf $\pi$-coalgebra endowed with a twist is said to be
{\it ribbon}. 

It follows from definitions that the neutral component $A_1$  of $A$
endowed with $ R_{ 1,1 }\in A_1\otimes A_1, \theta_1\in A_1$ is 
a ribbon   Hopf  algebra in the sense of   Section 11.1.
In particular, the equality
$\varphi_1=\id$ implies that $\theta_1$ lies in the center of $A_1$. 

For $\pi=1$,  the notions
introduced   in Sections 11.3 and 11.4  boil down to the standard notions of   
 quasitriangular and ribbon  Hopf algebras.

\skipaline \noindent {\bf 11.5.  Examples.}  We give here 
two examples of ribbon  Hopf $\pi$-coalgebras.
Let
  $(A, \Delta, \varepsilon, s,  R,v)$  be a  ribbon  Hopf
algebra  over $K$. An element $\alpha \in A$ is {\it group-like} if 
$\Delta(\alpha)= \alpha\otimes \alpha$ and $\varepsilon (\alpha)=1\in K$. Any
group-like element $\alpha$ is  invertible and
$  s(\alpha)=\alpha^{-1}$. The group-like elements
of $A$ form a group, $\pi=\pi(A)$,  under multiplication in $A$.      For   $\alpha\in
\pi$,   the formula $a\mapsto  \alpha a   \alpha^{-1}  $ with $a\in A$
  defines a Hopf algebra endomorphism of $A$. This gives  an action of $\pi$ on $A$
by Hopf algebra endomorphisms.  Applying the  constructions   of
Section 11.2 to   this action we obtain   crossed Hopf $\pi$-coalgebras $A^{\pi}$ and
$\overline A^{\pi}$.  We  define a universal $R$-matrix and twist   in
 $A^{\pi}$
by   $$R_{\alpha,\beta } =(i_\alpha \otimes i_\beta) (  (1_A \otimes \alpha^{-1})\,
R ) \in A_{\alpha}\otimes  A_{\beta}\,\,\,\,\,
{\text {and}}  \,\,\,\, \, \theta_\alpha=i_\alpha (v \alpha^{-1})
\in A_{\alpha} $$ where   $\alpha,\beta\in
\pi$. We  define a universal $R$-matrix and twist in
 $\overline  A^{\pi}$
by   $$\overline  R_{\alpha,\beta } =(i_\alpha \otimes i_\beta) (    R
\,(\beta^{-1} \otimes 1_A) ) \in A_{\alpha}\otimes  A_{\beta}  \,\,\,\, \,
{\text {and}}   \,\,\,\,\, \overline \theta_\alpha=i_\alpha (v \alpha^{-1}) \in
A_{\alpha} $$ where   $\alpha,\beta\in \pi$. 
A direct computation shows that  
$(A^{\pi}, \{R_{\alpha,\beta }  \}_{\alpha,\beta \in \pi},
 \{\theta_\alpha \}_{\alpha\in \pi})$ and $(\overline A^{\pi}, \{\overline
R_{\alpha,\beta }  \}_{\alpha,\beta \in \pi},
 \{\overline \theta_\alpha \}_{\alpha\in \pi})$   are ribbon  
Hopf $\pi$-coalgebras.

Group-like elements of a quantum universal enveloping algebra $A=U_q (\,\bold     g)$
 are well-known.  For instance,  if $\,\bold     g=sl(N+1)$ and $q$ is generic then 
  $\pi(A)=\bold Z^{N}$ is a free abelian group of rank $N$ generated by
the canonical group-like elements
$K_1,..., K_N\in A$. If $q$ is a primitive root of unity of order $\ell$, then 
one usually considers a version $A^{{\text {res}}}$ of $A=U_q (sl(N+1))$ with
$K_i^{\ell}=1$ for all $i=1,...,N$ (see [KRT]). Then $\pi(A^{{\text {res}}})=
(\bold Z/\ell \bold Z)^N$. 
 
\skipaline \noindent {\bf 11.6.  Operations on crossed Hopf group-coalgebras.}
In analogy with the  pull-back of group-categories, we can pull back   a   Hopf
$\pi$-coalgebra $A$   along a   group homomorphism
$q:\pi'\to \pi$. This gives  a   Hopf
$\pi'$-algebra $A'=q^*(A)$ defined by $A'_{\alpha}= A_{q(\alpha)}$  for any
$\alpha\in \pi'$. If $A$ is crossed (resp. quasitriangular,  ribbon), then $A'$
has the structure of a crossed (resp. quasitriangular,  ribbon)  Hopf
$\pi'$-algebra obtained by lifting the   data from $A$ to $A'$ in the
obvious way.  Taking $\pi=1$ and choosing as $A$ any quantum group, we obtain 
an example of a ribbon Hopf $\pi'$-algebra for any $\pi'$.

The constructions  of  direct sum, tensor product, and
transfer discussed in Sections 1.4,   2.5, 13  have their analogues for 
crossed, quasitriangular and ribbon Hopf group-coalgebras.  We leave the details to
the reader. (A   transfer for $\pi$-algebras  is discussed in a related 
setting  in [Tu3]). 

Given a  crossed Hopf $\pi$-coalgebra 
$A=(\{A_{\alpha}\}_{\alpha\in \pi},
\Delta , \varepsilon_1, s, \varphi)$ we define its mirror 
$\overline A=(\{\overline A_{\alpha}\}_{\alpha\in \pi},
\overline\Delta , \overline \varepsilon_1, \overline s, \overline \varphi)$.  For
$\alpha \in \pi$, set $\overline A_\alpha=A_{\alpha^{-1}}$.
For $a\in  \overline A_{\alpha \beta}$, set
$$\overline \Delta_{\alpha, \beta} (a) = 
(\varphi_\beta\otimes
\id_{A_{\beta^{-1}}}) \Delta_{\beta^{-1} \alpha^{-1}\beta, \beta^{-1}} (a)
$$
$$= (\id_{A_{\alpha^{-1}}} \otimes \varphi_{\beta^{-1}}) 
 \Delta_{\alpha^{-1}, \beta^{-1}} \varphi_\beta(a)\in A_{\alpha^{-1}}\otimes
A_{\beta^{-1}} =\overline A_\alpha \otimes \overline A_\beta.$$
This defines a comultiplication $\overline \Delta_{\alpha, \beta}:
\overline A_{\alpha \beta}\to \overline A_{\alpha  }\otimes \overline A_{ 
\beta}$. Set $\overline \varepsilon_1=\varepsilon_1:\overline A_1=A_1\to K$.
For $\alpha\in \pi$, set
$$ \overline s_\alpha=\varphi_\alpha s_{\alpha^{-1}}: \overline A_\alpha=A_{\alpha^{-1}}
\to A_{\alpha }=\overline A_{\alpha^{-1}}.$$
Finally set $\overline \varphi_\alpha=\varphi_\alpha$ for all $\alpha \in \pi$.
A direct computation shows that 
$\overline A=(\{\overline A_{\alpha}\}_{\alpha\in \pi},
\overline\Delta , \overline \varepsilon_1, \overline s, \overline \varphi)$ 
is a  crossed Hopf $\pi$-coalgebra. Moreover, if $R, \theta$ are a universal
$R$-matrix and a twist in $A$ then the formulas 
$${\overline
R}_{\alpha, \beta}=(\Perm (R_{\beta^{-1}, \alpha^{-1}}))^{-1}\in  
\overline A_{\alpha  }\otimes \overline A_{ 
\beta},\,\,\,\, {\overline
\theta}_{\alpha}= (\theta_{\alpha^{-1}})^{-1} \in  \overline A_{\alpha  } $$  define a
universal
$R$-matrix and  a twist in ${\overline A}$.  It is easy to see that   $\overline
{\overline A}=A$.  

For example, the crossed Hopf $\pi$-coalgebras   $A^{\pi}$ and 
$\overline A^{\pi}$ defined in  Section 11.2 are mirrors of each other.
The ribbon Hopf $\pi$-coalgebras $A^{\pi}$ and $\overline A^{\pi}$  constructed in
Section 11.5  are related as follows: $\overline A^{\pi}$ is the mirror 
of $B^{\pi}$ where $B=(A,(\Perm
(R))^{-1},v^{-1})$ and $\pi=\pi(A)=\pi(B)$.

\skipaline \noindent {\bf 11.7.  Categories of representations.}  We shall
associate with every   Hopf $\pi$-coalgebra $ A=(\{A_\alpha\}, \Delta ,
\varepsilon_1,s,\varphi)$ a category of representations $\Rep (A) $ which has a
natural structure of a  $\pi$-category.  Moreover, if $A$ is crossed
(resp. quasitriangular,  ribbon) then $\Rep (A) $ is crossed (resp. braided, 
ribbon).

By an $A_\alpha$-module we mean a left $A_\alpha$-module whose underlying
$K$-module is projective of finite type.  (The unit of $A_\alpha$ is supposed 
to
act as the identity.)  The category $\Rep (A) $ is the disjoint union of the
categories $\{\Rep (A_\alpha) \}_{\alpha\in \pi}$ where $\Rep (A_\alpha) $ is the
category of $A_\alpha$-modules and $A_\alpha$-linear homomorphisms.  The tensor
product and the unit object in $\Rep (A) $ are defined in the usual way using the
comultiplication $\Delta$ and the counit $\varepsilon_1$.  The associativity
morphisms are the standard identification isomorphisms for modules $(U\otimes
V)\otimes W= U\otimes (V \otimes W)$; they will be supressed from the notation.
The same applies to the morphisms $l,r$ which are the standard identifications
$U\otimes K= U=K\otimes U$.

For $U\in \Rep (A_\alpha) $, we have $U^*=\Hom_K(U,K)\in \Rep (A_{\alpha^{-1}})$ 
where
$a\in A_{\alpha^{-1}}$ acts as the transpose of $x\mapsto s_{\alpha^{-1}}(a)
(x):U\to U$.  The duality
morphism $d_U:  U^*\otimes U \to \1=K$ is the  evaluation
pairing; it gives rise to $b_U$ in the usual way, cf.  [Tu2, Chapter XI].  The second
line of
  equalities in (11.2.2) implies that $d_U, b_U$ are $A_1$-linear.

The automorphism $\varphi_\alpha$ of $A$ defines an automorphism,
$\Phi_\alpha$,  of $\Rep (A) $.   If
$U\in \Rep (A_\beta)$ then $\Phi_\alpha(U)$ has the same underlying $K$-module
as $U$ and each $a\in A_{\alpha\beta\alpha^{-1}}$ acts as multiplication by
$\varphi_\alpha^{-1}(a)\in A_\beta$.  Every $A_\beta$-homomorphism $U\to U'$ 
is
mapped to itself considered as a $A_{\alpha\beta\alpha^{-1}}$-homomorphism.  It 
is
easy to check that $\Rep (A) $ is a crossed $\pi$-category.

A    universal $R$-matrix  (11.3.a) in $A$ induces a braiding in $\Rep (A) $ 
as
follows.  For $U\in \Rep (A_\alpha) , V\in \Rep (A_\beta)$, the braiding
$c_{V,W}:V\otimes W \to {}^V\!  W \otimes V$ is the composition of 
multiplication
by $R_{\alpha,\beta }$, permutation $V\otimes W \to W \otimes V$ and the
$K$-isomorphism $W\otimes V = {}^V\!  W \otimes V$ which comes from the fact 
that
$W = {}^V\!  W$ as $K$-modules.  The conditions defining a universal $R$-matrix
ensure that $\{c_{V,W}\}_{V,W}$ is a braiding.

A twist  $ \theta $ in $A$ induces a twist in $\Rep (A) $: for any
$A_\alpha$-module $V$, the morphism $\theta_V: V\to {}^V\!  V$ is   the
composition of multiplication by $\theta_\alpha\in A_\alpha$ and the
$K$-isomorphism $ V \to {}^V\!  V$ which comes from the fact that ${}^V\!  V=V$
as $K$-modules.  Conditions (11.4.1), (11.4.4)  imply that the  resulting
homomorphism $ V \to {}^V\!  V$ is $A_\alpha$-linear.  Condition (2.3.1) 
follows
from definitions, conditions (2.3.2) - (2.3.4) follow from (11.4.2) - (11.4.4),
respectively.  Thus, $\Rep (A) $ is a ribbon  crossed  $\pi$-category.

As an exercise the reader may check that $\overline { \Rep (A) }=\Rep (\overline
{A})$. 

\skipaline \noindent {\bf 11.8.  Remarks.}   
1. The idea of a  Hopf group-coalgebra  comes from the following 
observation. Consider a topological group $G$. For   $\alpha\in G$, denote  by
$C_\alpha=C_\alpha (G)$ the algebra of   germs of continuous functions in  
$\alpha\in G$. The group multiplication  $G\times G\to G$ induces an algebra
homomorphism $ C_1  \to  C_1\hat {\otimes} C_1 $ which   makes $C_1 $   a
topological Hopf algebra. Similarly, the group
multiplication  in $ G$ induces an algebra homomorphism
 $\Delta_{\alpha,\beta}: C_{\alpha\beta}   \to  C_\alpha \hat {\otimes} 
C_\beta$ for any $\alpha, \beta\in G$. This makes the  system of algebras
$\{C_\alpha   \}_{\alpha\in G}$   a (topological) Hopf $G$-coalgebra. 

In this example  we can   compute 
$\Delta_{\alpha,\beta} $ via $\Delta_{1,1}$
 and the   adjoint action of $G$ on $C_1 $
as follows.  Observe first that   left multiplication by $  \alpha $ induces an
algebra isomorphism, $i_\alpha: C_{1}  \to C_\alpha $. 
Let $x,y$ be two elements of $G$ close  to $1\in G$. Then $\alpha x, \beta
y$ are close   to $\alpha, \beta$, respectively. For   $f\in C_{1}  $,
$$
(i_\alpha^{-1} \otimes i_\beta^{-1}) \Delta_{\alpha,\beta} (i_{\alpha\beta} (f))
( x,  y)= \Delta_{\alpha,\beta} (i_{\alpha\beta} (f))
(\alpha x, \beta y)$$
$$= i_{\alpha\beta} (f) (\alpha x  \beta y)=f( (\alpha\beta)^{-1} 
\alpha x  \beta y)= f(  \beta^{-1}    x  \beta y)=
\Delta_{1,1} (f)
( \beta^{-1}    x  \beta,y).$$
This  
computation leads to the second example of  Hopf group-coalgebras in Section 11.2. 
 
2.  Many
aspects of the   theory  of Hopf algebras generalize to  crossed  Hopf
$\pi$-coalgebras. The notion of a modular Hopf algebra [RT] can be extended to this
setting which gives rise to modular crossed $\pi$-categories.  
The Drinfeld's double construction  and the Drinfeld's
theory of quasi-Hopf algebras   also generalize  to this setting. In particular, the 
double of a  
 crossed  Hopf
$\pi$-coalgebra   is a quasitriangular crossed  Hopf
$\pi$-coalgebra, see [Zu].   

3.    By now there is quite a number of constructions of topological invariants  of 3-manifolds from Hopf algebras.  Note
in particular the constructions introduced by  Turaev-Viro [TV],
Hennings [He],  and  Kuperberg [Ku]. All of them     can
be
generalized to our  setting. This allows to derive invariants of 3-dimensional
$\pi$-manifolds from certain  crossed  Hopf
$\pi$-coalgebras,  see Appendix 2 and [Vi]. 

4.   There is a special case of the definitions above where all
 the crossed isomorphisms
$\{\varphi_\alpha\}_{\alpha\in \pi}$ are the identity maps.  Namely, assume that  $A$ is a Hopf $\pi$-coalgebra  over
an abelian group $\pi$. Then the trivial homomorphism
$\varphi=1:\pi
\to
\Aut (A)$ makes
$A$  a crossed Hopf $\pi$-coalgebra.  In this case the definitions of the universal $R$-matrix and the twist in $A$ considerably
simplify. Such  quasitriangular crossed  Hopf
$\pi$-coalgebras were first considered by T. Ohtsuki [Oh1], [Oh2].  He calls them 
colored Hopf algebras and derives   examples of such algebras
from $U_q(sl_2)$.

\skipaline \centerline {\bf 12.  Canonical extensions}

\skipaline
We show in  this section 
that a monoidal (resp. braided, ribbon) category canonically extends
to a certain crossed  (resp. braided, ribbon) group-category.

\skipaline \noindent {\bf 12.1.  Extensions of a monoidal category.}
Let ${\Cal C}$ 
be
 a $K$-additive monoidal category with left duality. Observe
that ${\Cal C}$ is a $\{1\}$-category where $\{1\}$ is a trivial group.
For any group $\pi$, the
category ${\Cal C}$ gives rise to a  $\pi$-category ${\Cal C}^{\pi}$ obtained by pulling
back ${\Cal C}$ along the trivial homomorphism $\pi \to 1$. By definition, 
${\Cal C}^{\pi}=\amalg_{\alpha\in \pi} {\Cal C}^{\pi}_\alpha$ where the objects of
${\Cal C}^{\pi}_\alpha$ are   pairs $(U\in  {\Cal C}, \alpha )$. A morphism $(U ,
\alpha  )\to  (V , \beta  )$ with $\alpha,\beta\in \pi$  is  $0 \in \Hom_{\Cal C}(U,V)$
 if $\alpha\neq \beta$ and any element of
$  \Hom_{\Cal C}(U,V)$  if $\alpha=\beta$. The operations on objects and the unit object
are defined by $$ (U , \alpha  )\otimes  (V , \beta  ) =(U \otimes V, \alpha    \beta),  
\,\,(U , \alpha  )^*=(U^*, \alpha^{-1}), \,\, \1_{{\Cal C}^{\pi}}= (\1_{{\Cal C}},1).$$
The structural morphisms (1.1.a,b,e) and composition and tensor product of
morphisms are induced by the corresponding operations in ${\Cal C}$ in the obvious way. 
 
Assume now that $\pi $ acts on ${\Cal C}$ by   automorphisms.   Such an action is
determined by a group homomorphism $\pi\to \Aut({\Cal C})$ where $\Aut  ({\Cal
C})$ is  the group  of automorphisms of ${\Cal C}$  defined in 
Section 2.1.  Then   ${\Cal C}^{\pi}$ acquires
the structure of a crossed $\pi$-category as follows. For  $\alpha\in \pi$ and  $(V ,
\beta  )\in {\Cal C}^{\pi}_\beta $, set $$\varphi_{\alpha} (V , \beta  ) =
(\alpha(V), \alpha \beta \alpha^{-1}) \in {\Cal C}^{\pi}_{\alpha \beta
\alpha^{-1}}.$$  For a morphism 
$f:(V , \beta  )\to  (W ,\gamma  )$ in  ${\Cal C}^{\pi}$ with $f\in \Hom_{\Cal
C}(V,W)$ set  $$\varphi_{\alpha} (f  )=\alpha(f)\in  \Hom_{\Cal
C}(\alpha(V),\alpha(W)) =\Hom_{{\Cal C}^{\pi}} 
(\varphi_{\alpha} (V , \beta  ), \varphi_{\alpha} (W ,\gamma  )).$$
  All the axioms of a crossed $\pi$-category follow from definitions. 
Applying this construction to $\pi=\Aut  ({\Cal C})$, we obtain a canonical extension of
${\Cal C}$ to a crossed $\Aut  ({\Cal C})$-category.

\skipaline \noindent {\bf 12.2.   The group} $\Aut_0({\Cal C})$. 
Let ${\Cal C}$ 
be
 a $K$-additive monoidal category with left duality. Denote by $\id_{\Cal C}$
the identity functor ${\Cal C}\to {\Cal C}$ sending each object and each morphism of ${\Cal C}$ into
itself.  We  introduce a
 group   $\Aut_0({\Cal C})$, formed by  monoidal 
  equivalences of    $\id_{\Cal C}$ to   automorphisms
  of ${\Cal C}$. 
More precisely, 
an element  of $ \Aut_0({\Cal C})$   is a pair  ($\alpha\in \Aut ({\Cal C})$, an invertible
monoidal morphism of functors $F :\id_{\Cal C} \to \alpha$).  The latter means that for
each object  $U\in {\Cal C}$ we have  an invertible morphism $F_U:U \to \alpha (U)$ in
${\Cal C}$ such that

(i) for any  morphism $f:U\to V$ in
${\Cal C}$ the following diagram is commutative: $$ \CD
 U  @> f >>  V \\ 
  @V {F_U} VV    @VV {F_V} V        \\
\alpha (U) @>\alpha (f)>> \alpha (V) ;
\endCD $$

(ii) $F_{\1}=\id_{\1}$ and for any $U,V\in {\Cal C}$, we have $F_{U\otimes V}=
F_{U}\otimes F_{V}$.

The product of
two   pairs  $(\alpha, F ), (\alpha', F') \in \Aut_0 ({\Cal C})$ is the pair $(\alpha \alpha', F 
  F')$ where for each
   $U\in {\Cal C}$ we have   $(F 
  F')_U= F_{\alpha'(U)} F'_U:U\to (\alpha \alpha')(U)$.  It is clear
that  $\Aut_0({\Cal C})$ is a group with respect to this multiplication.
The pair $(\alpha=\id_{\Cal C}, F=\{\id_U\}_{U\in {\Cal C}})$   is the unit object of $ \Aut_0({\Cal C})$.

Forgetting $F$, we obtain a group homomorphism
$i_0:\Aut_0({\Cal C})\to \Aut ({\Cal C})$ whose 
image   consists of the automorphisms of ${\Cal C}$ monoidally
equivalent to   $\id_{\Cal C} $.  The key property of the elements
of the image  is given by the following lemma.

  \skipaline \noindent {\bf 12.2.1.  Lemma.  } {\sl If ${\Cal C}$ is braided then all
elements of $ i_0(\Aut_0({\Cal C})) \subset \Aut ({\Cal C})$ preserve the braiding. 
If ${\Cal C}$ is ribbon then all
elements of $ i_0(\Aut_0({\Cal C}))$ preserve both the braiding and the twist.}

\skipaline {\sl Proof.}  Let  $\{c_{U,V}:U\otimes V \to   V \otimes U\}_{U,V\in {\Cal C}}
 $ be a braiding in ${\Cal C}$. Let $(\alpha, F)\in \Aut_0({\Cal C})$.
Then for any $U,V\in {\Cal C}$ we have a commutative diagram 
$$ \CD
 U\otimes V@> c_{U,V} >>  V\otimes U  \\ 
  @VF_{U\otimes V} VV    @VV F_{V\otimes U} V        \\
\alpha (U) \otimes \alpha(V) @>\alpha (c_{U,V})>> \alpha (V) \otimes \alpha(U).
\endCD \leqno (12.2.a)$$
On the other hand, by the naturality of the braiding we have
$$F_{V\otimes U} \,c_{U, V}= (F_{V}\otimes F_{U}) \,c_{U, V}=
c_{\alpha (U), \alpha(V)}(F_{U}\otimes F_{V})=
c_{\alpha (U), \alpha(V)} F_{U\otimes V}.$$
Thus, if we replace in   (12.2.a) the bottom arrow by 
$c_{\alpha (U), \alpha(V)}$ we obtain a commutative diagram. 
Since  all the arrows in (12.2.a) are invertible morphisms,  $c_{\alpha (U),
\alpha(V)}=\alpha (c_{U,V})$.  Similarly, if $ \{\theta_{U}:U \to    U\}_{U\in {\Cal C}}
$ is a twist in ${\Cal C}$ then for any $U \in {\Cal C}$ we have a commutative diagram 
$$ \CD
 U @> \theta_{U} >>   U  \\ 
  @VF_{U } VV    @VV F_{  U} V        \\
\alpha (U)   @>\alpha (\theta_{U})>>   \alpha(U).
\endCD  $$
By the naturality of the twist we have
 $ F_{U}  \theta_{U}=\theta_{\alpha (U) } F_{U } $
and therefore $\alpha (\theta_{U})=\theta_{\alpha (U) }$.

\skipaline \noindent {\bf 12.3.  The canonical extension of a ribbon category.}
Let ${\Cal C}$ 
be
 a braided   $K$-additive monoidal category with left duality. 
Set $\pi=\Aut_0({\Cal C})$. We  
define here a canonical extension of ${\Cal C}$ to a braided   
crossed  $\pi$-category.
Applying the constructions of Section 12.1 to the 
group homomorphism
$i_0:\pi=\Aut_0({\Cal C})\to \Aut ({\Cal C})$ we obtain a crossed $\pi$-category ${\Cal C}^{\pi}$.
The   braiding $\{c_{U,V}:U\otimes V \to   V \otimes U\}_{U,V\in {\Cal C}}$ in ${\Cal
C}$ induces a braiding in ${\Cal C}^{\pi}$ as follows. Let $u=(U,(\alpha,F)),  
v=(V,(\beta, G))$ be objects of ${\Cal C}^{\pi}$ where $U,V\in {\Cal C}$ and $(\alpha,F),(\beta,
G)\in \pi$. Observe that 
$$u\otimes v=(U\otimes V, (\alpha,F) (\beta, G)),\,\, \varphi_{(\alpha,F)}
(v)=(\alpha(V), (\alpha,F) (\beta, G) (\alpha,F)^{-1}),$$
$$ {}^u\! v \otimes  u=\varphi_{(\alpha,F)}
(v) \otimes u=( \alpha(V)\otimes U, (\alpha,F) (\beta, G)).$$
The invertible morphism 
$(F_V\otimes \id_U) \,c_{U,V}: U\otimes V \to \alpha(V)\otimes
U$ in ${\Cal C}$ defines an invertible morphism  in   ${\Cal C}^{\pi}$ 
  $$c_{u, v}:
u\otimes v\to  {}^u\! v \otimes u.\leqno
(12.3.a)$$

 If ${\Cal C}$ is a ribbon category, then the 
 twist   $ \{\theta_{U}:U \to    U\}_{U\in {\Cal C}}
$   in ${\Cal C}$ induces a twist in ${\Cal C}^{\pi}$ as follows. Let 
$u=(U,(\alpha,F))\in
{\Cal C}^{\pi}$. Then the invertible morphism $F_U \theta_U:U\to \alpha(U)$ in ${\Cal C}$
defines an invertible morphism in   ${\Cal C}^{\pi}$
$$\theta_{u}:u \to  {}^u\! u=\varphi_{(\alpha,F)}
(u)=(\alpha(U), (\alpha,F)  ). \leqno
(12.3.b)$$

  \skipaline \noindent {\bf 12.3.1.  Theorem.  } {\sl The    
morphisms (12.3.a) with $u,v\in {\Cal C}^{\pi}$ form  a braiding in ${\Cal C}^{\pi}$. 
If ${\Cal C}$ is ribbon, then the  
  morphisms (12.3.b) with $u\in {\Cal C}^{\pi}$ form  a twist in ${\Cal C}^{\pi}$.
If ${\Cal C}$ is modular, then so is ${\Cal C}^{\pi}$.}
 \skipaline
The proof goes by  a routine verification  of  the  axioms, we leave it to
the reader. Note that the neutral component of ${\Cal C}^{\pi}$
is   ${\Cal C} $ with its original braiding and twist.

 \skipaline \noindent {\bf 12.4.  Remarks and examples.  }
1. Condition (12.2.i) on a pair $(\alpha,F)\in \Aut_0({\Cal C})$ shows that the action of
$\alpha$ on morphisms is completely determined by $F$. Setting $F(U)=\alpha(U)$
for $U\in {\Cal C}$, we can   reformulate the definition  of $ \Aut_0({\Cal C})$  entirely in
terms of $F$.  An element  of $
\Aut_0({\Cal C})$ is thus   described as  a pair (a bijection
  $F$ from the set of objects of ${\Cal C}$ into itself, a system of invertible
morphisms $\{F_U:U\to F(U)\}_{U\in {\Cal C}}$) such that 

(i) $F(\1)=\1$ and
$F_{\1}=\id_{\1}$; 

(ii) for any $U,V\in {\Cal C}$, we have $F(U\otimes V)=
F(U)\otimes F(V)$ and $F_{U\otimes V}=
F_{U}\otimes F_{V}$, 

(iii) for any $U\in {\Cal C}$, we have $F(U^*)=(F(U))^*$
and
$F_{U^*}= ((F_U)^*)^{-1}$ where $(F_U)^*:(F(U))^*\to U^*$ is the transpose of $F_U$.
 
2. Assume that in Example 2.6, the
group $\pi$ is abelian.  Then the corresponding category
${\Cal C}$  is a ribbon monoidal category
in the usual sense of the word.  Any automorphism of ${\Cal C}$ monoidally
equivalent to
$\id_{\Cal C}$ is
equal to $\id_{\Cal C}$   since  all non-zero morphisms in 
${\Cal C}$ are proportional to the identity  morphisms of objects.  An
element $(\alpha,F)\in \Aut_0({\Cal C})$ is  therefore completely determined by the 
map $U\mapsto F_U\in \Hom_{\Cal C} (U, U)=K$ where $U$ runs over the
elements of $\pi$.
  The inclusion  $(\alpha,F)\in \Aut_0({\Cal C})$ is equivalent to the
condition that this map is a group homomorphism 
  $ \pi\to K^*\subset K$. Hence, 
$\Aut_0({\Cal C})=\Hom (\pi,  K^*)= \pi^*$. By Theorem 12.3.1, ${\Cal C}$  gives rise to
a ribbon crossed $\pi^*$-category ${\Cal C}^{\pi^*}$.

3. Let $G$ be a topological group and ${\Cal C}=\Rep (G)$ be the category of finite
dimensional linear representations of $G$ and $G$-linear homomorphisms. It is
clear that ${\Cal C}$ is a $K$-additive monoidal category with left duality. Moreover, this
category is braided  (in fact symmetric) with braiding   given by the  
flips (permutations) $U\otimes V\to V\otimes U$. The
category ${\Cal C}$ is ribbon with twist $\theta_U=\id_U:U\to U$ for all $U\in {\Cal C}$. 
There is a  
group homomorphism $g\mapsto (\alpha^g, F^g):G\to  \Aut_0({\Cal C})$ defined as
follows. Given $g\in G$, the     functor $\alpha^g : {\Cal C}\to {\Cal C}$  sends a 
$G$-module $(U,\rho_U:G\to \End (U))$  into  the same
$K$-module $U$ where  each $ h\in G$ acts as $\rho_U ( g^{-1}hg)$. The     functor
$\alpha^g$ sends any  $G$-linear homomorphism $f:U\to V$  into itself. It is clear
that $\alpha^g \in \Aut({\Cal C})$.
  The   morphism $F^g_U:U\to \alpha^g (U)$ in ${\Cal C}$ is defined by $u\mapsto
\rho_U(g^{-1}) (u):U\to U$ where $u$ runs over $U$.
The family $\{F^g_U\}_{U\in {\Cal C}}$ satisfies the conditions of Section 12.2.
Thus,
  $(\alpha^g, F^g)\in  \Aut_0({\Cal C}) $.  Pulling back the canonical extension
${\Cal C}^{ \Aut_0({\Cal C})}$ along the  group homomorphism $g\mapsto (\alpha^g, F^g):G\to 
 \Aut_0({\Cal C}) $ we obtain a ribbon crossed $G$-category $(\Rep (G))^G$. 

If $G$ is a semisimple complex connected Lie group with Lie algebra
$\,\bold     g$, then the category $\Rep (G)$ coincides with the category $U(\,\bold    
g)$-mod  of finite dimensional $U(\,\bold     g)$-modules. It  would be most  
interesting to find a quantum deformation of the ribbon crossed $G$-category $(\Rep
(G))^G= (U(\,\bold     g)-\mod)^G$ generalizing the deformation $U_q(\,\bold     g)$-mod 
 of $U(\,\bold     g)$-mod  arising in the theory of quantum groups.

4.  One can generalize the constructions of this section to the case where
the initial category ${\Cal C}$ is a (braided or ribbon) crossed group-category; one
should then consider only those automorphisms of ${\Cal C}$ which are compatible
 with the given group action. We shall not pursue this line here.

\skipaline \centerline {\bf 13.  Transfer of categories}

\skipaline
 Throughout this section 
$G\subset
\pi$ is a subgroup of a group $\pi$. We shall show that each $G$-category ${\Cal C}$ 
gives rise to  a $\pi$-category $\tilde {\Cal C}$ via a natural transfer. If ${\Cal C}$ is crossed
(resp. braided, ribbon) then $\tilde {\Cal C}$ is crossed
(resp. braided, ribbon).

\skipaline \noindent {\bf 13.1.  Transfer for group-categories.}  Fix a $G$-category
${\Cal C}$. We shall construct a   $\pi$-category $ \tilde {\Cal C}$  called the {\it
transfer} of ${\Cal C}$.

Fix a representative $\omega_i\in \pi$ for each right coset class $i\in
G\backslash \pi$ so that $i=G\omega_i \subset \pi$.  For $\alpha\in \pi$, set
$$N(\alpha)= \{i\in G\backslash \pi\, \vert \, \omega_i \alpha \omega_i^{-1} 
\in
G\} \subset G\backslash \pi.  $$ An object $U$ of $\tilde{\Cal C}$ is a triple
$(\alpha\in \pi$, a subset $A$ of $N(\alpha)$, a family $\{U_i\in
{\Cal C}_{\omega_i\alpha\omega_i^{-1}}\}_{ i\in A}$).  The set $A$ (which may be void)
will be denoted by $\vert U\vert$.  The morphisms in $\tilde {\Cal C}$ are defined by
$$\Hom_{\tilde {\Cal C}} (U, U')= \prod_{i\in \vert U\vert\cap \vert U'\vert}
\Hom_{\Cal C}(U_i,U'_i).   $$
Thus a morphism $f:U\to U'$ in $\tilde {\Cal C}$ is a family $\{f_i:U_i\to U'_i\}_{i\in
\vert U\vert\cap \vert U'\vert}$ of morphisms in ${\Cal C}$.  We view each $f_i$ as the 
$i$-th
coordinate of $f$.  The $K$-linear structure in $\Hom_{\tilde {\Cal C}} (U, U')$ is
coordinate-wise.  The composition of two morphisms $f:U\to U'$ and $f':U'\to 
U''$
is defined in coordinates by $$(f'f)_i= \cases f'_i f_i:U_i\to U''_i,~ {
{if}}\,\,\, i\in  \vert U\vert\cap \vert U'\vert \cap \vert U''\vert, \\ 
0:U_i\to
U''_i,~ { {if}}\,\,\, i\in (\vert U\vert\cap \vert U''\vert) \backslash \vert
U'\vert.\endcases $$ This is an associative composition in $\tilde{\Cal C} $ because
the composition of morphisms in ${\Cal C}$ is associative and the composition of a 
zero
morphism in ${\Cal C}$ with any other morphism is again a zero morphism.  This defines
$\tilde {\Cal C}$ as a $K$-additive category.  We have a splitting $\tilde
{\Cal C}=\amalg_{\alpha\in \pi} \tilde {\Cal C}_\alpha$ where $\tilde {\Cal
C}_\alpha$ is the full subcategory of $\tilde {\Cal C}$ with objects $(\alpha , A ,
\{U_i \}_{i\in A})$.

The unit object of $\tilde {\Cal C}$ is  the triple $(1\in \pi,
A=G\backslash \pi, \{U_i=\1_{{\Cal C}} \}_{ i\in A})$.  The duality and tensor product
for objects of $\tilde {\Cal C}$ are defined by $$ (\alpha , A , \{U_i \}_{i\in A})^*=
(\alpha^{-1} , A , \{U^*_i \}_{i\in A}), $$ $$(\alpha , A , \{U_i \}_{i\in A})
\otimes (\beta , B , \{V_j \}_{j\in B})=(\alpha\beta , A\cap B , \{U_i\otimes 
V_i
\}_{ i\in A\cap B}).  $$ Note that $N(\alpha)=N(\alpha^{-1})$ and that the
inclusions $A\subset N(\alpha), B\subset N(\beta)$ imply that  $A\cap B
\subset N(\alpha\beta)$.  Observe the identities $\vert U^* \vert = \vert 
U\vert
, \vert U\otimes V \vert =\vert U \vert \cap \vert V \vert $.

The tensor product of morphisms $f:U\to U'$ and $g:V\to V'$ is defined   by
$$(f\otimes g)_i= f_i \otimes g_i:U_i \otimes V_i\to U'_i  \otimes
V'_i$$ for all $ i\in \vert U\vert\cap \vert V\vert\cap \vert U'\vert\cap \vert
V'\vert$.  It is a simple exercise  to check the identity
$(f'\otimes g') (f\otimes g)=f'f\otimes g'g$.

The structural morphisms $a,l,r,b,d$ in $\tilde {\Cal C}$ are all defined coordinate-wisely
and their coordinates are the corresponding morphisms in ${\Cal C}$.  In particular, for
every $U \in \tilde {\Cal C}$ we define $b_U:\1\to U\otimes U^*$ by $(b_U)_i=
b_{U_i}:\1\to U_i\otimes (U_i)^*$ where $i\in \vert U\vert$.  Similarly,
$(d_U)_i= d_{U_i}, (l_U)_i= l_{U_i}, (r_U)_i= r_{U_i}$ where $i\in \vert
U\vert$.  The associativity morphisms are defined by
$(a_{U,V,W})_i=a_{U_i,V_i,W_i}$ for all $i\in \vert U\vert\cap \vert V\vert\cap
\vert W\vert$.  The naturality of $a,l,r$ and the identities (1.1.c,d,f,g) 
follow
from the corresponding properties of ${\Cal C}$.
In general,  the isomorphism class of $\tilde {\Cal C}$   depends on
the choice of $\omega_i\in i$.

\skipaline \noindent {\bf 13.2.  Transfer for crossed group-categories.}   For each 
crossed $G$-category $({\Cal C}, \varphi:G\to \Aut ({\Cal C}))$, its transfer    $
\tilde {\Cal C}$  is a crossed $\pi$-category.   We need only to define the action   of
$\pi$ on $\tilde {\Cal C}$.  Fix 
the representatives $\omega_i\in  i$ for   $i\in
G\backslash \pi$ needed in the construction of $ \tilde {\Cal C}$.
Consider the    left action of
$\pi$ on $G\backslash \pi$ defined  by $ {\alpha} (i) = i\alpha^{-1}$ for $\alpha\in
\pi, i \in G\backslash \pi$.   
We have   $G \omega_{\alpha(i)}=G\omega_i \alpha^{-1}$ so that 
$\alpha_i=\omega_{\alpha(i)} \alpha
\omega_i^{-1} \in G$ for all $\alpha\in \pi, i \in G\backslash \pi$.

  For  $\beta\in \pi$, the map $j\mapsto {{\alpha(j)}}$ sends 
bijectively $N(\beta)$ onto $ N(\alpha \beta \alpha^{-1})$.  For every $j\in
N(\beta)$, we have the functor $$\varphi_{\alpha_j}:  {\Cal C}_{\omega_j \beta
\omega_j^{-1}} \to  {\Cal C}_{\omega_{\alpha(j)}
\alpha \beta \alpha^{-1} (\omega_{\alpha(j)})^{-1}}.  $$ Given an object 
$V=(\beta , B , \{V_j \}_{j\in B})$ of $\tilde {\Cal C}_\beta$, we apply
$\{\varphi_{\alpha_j}\}_{j\in B}$ coordinate-wisely to obtain an object
$\tilde \varphi_\alpha (V)=(\alpha \beta \alpha^{-1}, 
\alpha(B)
,\{\varphi_{\alpha_j}(V_j)\}_{j\in B})$ of $\tilde {\Cal C}_{ \alpha \beta \alpha^{-1}}$.
More precisely,  $$\tilde \varphi_\alpha (V)=(\alpha \beta \alpha^{-1}, 
\alpha(B)
, \{\varphi_{\alpha_{\alpha^{-1}(i)}} (V_{\alpha^{-1}(i)}) \}_{i\in \alpha(B)}).  
$$
Clearly, $\vert \tilde \varphi_\alpha (V) \vert=\alpha(\vert V\vert)$.  The 
action
of $\tilde \varphi_\alpha$ on a morphism $f:V\to V'$ is given by $$\tilde
\varphi_\alpha (f) = \{\varphi_{\alpha_{\alpha^{-1}(i)}} (f_{\alpha^{-1}(i)}
:V_{\alpha^{-1}(i)}\to V'_{\alpha^{-1}(i)})\}_{i \in \alpha (\vert V\vert)\cap
\alpha(\vert V'\vert) }.$$ The functor $\tilde \varphi_\alpha$ preserves the
structural morphisms in $\tilde {\Cal C}$ because by assumptions the functors $\{ 
\varphi_{\alpha_j}\}_j$ preserve the structural morphisms in $ {\Cal C}$.  The category
$\tilde {\Cal C}$ with action $\tilde \varphi$ of $\pi$ satisfies all the axioms of a
crossed $\pi$-category.  The isomorphism class of $\tilde {\Cal C}$ does not depend on
the choice of $\omega_i\in i$.

\skipaline \noindent {\bf 13.3.  Transfer for braided and ribbon group-categories.} 
If ${\Cal C}$ is a braided 
 crossed  $G$-category then its transfer    $ \tilde {\Cal C}
$  acquires the structure of    a braided $\pi$-category as follows. 
Let
$U=(\alpha , A , \{U_i \}_{i\in A})$ and $V= (\beta , B , \{V_j \}_{j\in B})$
be objects of $\tilde {\Cal C}$. 
By definition,
$$U\otimes V= (\alpha \beta, A \cap B , \{U_i \otimes V_i\}_{i\in A\cap B}),$$
$${}^{U }\!  V=  \tilde \varphi_\alpha (V)=(\alpha \beta \alpha^{-1}, 
\alpha(B)
, \{\varphi_{\alpha_{\alpha^{-1}(i)}} (V_{\alpha^{-1}(i)}) \}_{i\in \alpha(B)}),$$
$${}^{U }\!  V\otimes U=
(\alpha \beta  , 
\alpha(B) \cap A
, \{\varphi_{\alpha_{\alpha^{-1}(i)}} (V_{\alpha^{-1}(i)}) \otimes U_i\}_{i\in
\alpha(B)\cap A}).$$
The latter expression simplifies if we observe that for any $i\in N(\alpha)$,
we have $G\omega_i=G\omega_i \alpha^{-1}$ and therefore
$\alpha(i)=i$. This implies that the map $i\mapsto \alpha(i): G\backslash \pi\to
G\backslash \pi$ is the identity on $A$. Hence,
$\alpha(B) \cap A=B \cap A=A\cap B$ and
$${}^{U }\!  V\otimes U=
(\alpha \beta  , 
A\cap B
, \{\varphi_{\alpha_{i}} (V_{i}) \otimes U_i\}_{i\in A\cap B}).$$
Now, we define a morphism $c_{U,V}:U\otimes V \to 
{}^{U }\!  V\otimes U$ by its coordinates
$(c_{U,V})_i = c_{U_i,V_i}:U_i\otimes V_i\to 
\varphi_{\alpha_{i}} (V_{i}) \otimes U_i$
where $i$ runs over $A\cap B$ and $c_{U_i,V_i}$ is the given braiding in ${\Cal C}$. 
Note that $\varphi_{\alpha_{i}} (V_{i})=\varphi_{\omega_{i} \alpha\omega_i^{-1}}
(V_{i})= {}^{U_i }\!  (V_i)$. The morphisms $\{c_{U,V}\}_{U,V}$ satisfy all the
axioms of a braiding in the crossed $\pi$-category $\tilde {\Cal C}$; this is easily
verified coordinate-wisely.

Similarly, if  ${\Cal C}$ is a ribbon
  crossed   $G$-category then its transfer    $ \tilde {\Cal C}
$  has the structure of    a ribbon $\pi$-category.
For an object $U=(\alpha , A , \{U_i \}_{i\in A})$   of $\tilde {\Cal C}$ we have
$ {}^{U }\!  U= (\alpha , A , \{\varphi_{\alpha_{i}} (U_{i}) \}_{i\in A})$.
Here $\alpha_{i}=\omega_{i} \alpha\omega_i^{-1}$ so that
$\varphi_{\alpha_{i}} (U_{i})={}^{U_i }\!  (U_i)$.
The twist $\theta_U:U\to {}^{U }\!  U$ is defined 
by its coordinates $(\theta_U)_i=\theta_{U_i}:U_i\to {}^{U_i }\!  (U_i)$
where $i$ runs over $A$. 
The morphisms $\{\theta_{U}\}_{U}$ satisfy all the axioms
of a twist in the braided $\pi$-category $\tilde {\Cal C}$.

If $G\neq \pi$ then the transfer $\tilde {\Cal C}$ of a modular $G$-category 
${\Cal C}$ is not
modular:    the algebra of endomorphisms of the unit object of $\tilde {\Cal C}$
equals $(\End_{\Cal C} (\1_{\Cal C}))^{[\pi:G]}=K^{[\pi:G]}\neq K$.

\skipaline \centerline {\bf   Appendix 1. Quasi-abelian  cohomology of groups}

\skipaline

Let  
  $a=\{a_{\alpha,\beta,\gamma} \in K^*\}_{\alpha,\beta,\gamma\in \pi}$ be a
3-cocycle of a group $\pi$ with values in   $K^*$ and  
$c=\{c_{\alpha,\beta}\in
K^*\}_{\alpha,\beta  \in
\pi}$  be a family of elements of $K^*$.
 Equations   (2.6.a,c,e,f)  on the pair $(a,c)$ are equivalent to the equations
introduced by C.  Ospel  [Os]
from   a different viewpoint (in his notation $ c_{\alpha, \beta }=\Omega_{
\beta, \alpha}$ and $a_{\alpha,\beta,\gamma}= f(\alpha,\beta,\gamma)$). 
Following Ospel, we call  such pairs $(a,c)$  
satisfying  (2.6.a,c,e,f) {\it 
quasi-abelian
3-cocycles} on $\pi$.   Examples of such cocycles are provided by pairs $(a,c)$
where  $a=1$ and $c $ is   a bilinear form
$(\alpha,\beta)\mapsto c_{\alpha, \beta }:  H_1(\pi) \times H_1(\pi) \to K^*$.

Quasi-abelian
3-cocycles $(a,c)$   on $\pi$ form a
commutative group under pointwise multiplication with unit $a=1,c=1$.   
This group   contains the coboundaries of the
conjugation invariant 2-cochains.  A conjugation invariant 2-cochain   is 
a
map $(\alpha,\beta)\mapsto \eta_{\alpha,\beta}:\pi\times \pi \to K^*$ such that
$\eta_{\delta\alpha\delta^{-1}, \delta \beta \delta^{-1}}= \eta_{\alpha, \beta 
}$
for any $\alpha, \beta , \delta\in \pi$.  Its coboundary is defined by
$$a_{\alpha, \beta, \gamma}=\eta_{\alpha, \beta } \eta_{\alpha \beta, \gamma}
(\eta_{\alpha, \beta \gamma})^{-1} (\eta_{ \beta, \gamma})^{-1},\,\,\, 
c_{\alpha,
\beta }=\eta_{\alpha, \beta } (\eta_{\beta, \alpha })^{-1}.$$ A direct
computation shows that this is a quasi-abelian
3-cocycle.  The 
quotient
of the group of quasi-abelian 3-cocycles by the subgroup of coboundaries is the
group of quasi-abelian cohomology of $\pi$ denoted $H^3_{qa}(\pi; K^*)$.
Involution  (2.6.i)  transforms 
the subgroup of coboundaries into itself and defines   an involution on 
$H^3_{qa}(\pi; K^*)$.  Forgetting   $c$ we obtain a homomorphism to the usual
cohomology $H^3_{qa}(\pi; K^*)\to H^3 (\pi; K^*)$. 

The constructions of Section 2.6 associate with any  quasi-abelian 3-cocycle of
$\pi$ and a conjugation invariant family $\{b_\alpha\in K^*\}_{\alpha\in
\pi}$ a certain braided crossed $\pi$-category. A twist 
in this category is  determined by a family $\{\theta_\alpha\in K^*\}_{\alpha\in
\pi}$ satisfying equations (2.6.g,h). The next lemma   explicitly describes all the
solutions to (2.6.g,h).   

 \skipaline \noindent {\bf    Lemma.  } {\sl Let $(a,c)$ be
 a quasi-abelian 3-cocycle of
$\pi$.  Then $\{\theta_\alpha=c_{\alpha,\alpha}\}_{\alpha\in
\pi}$ satisfies (2.6.g,h).  A general solution to (2.6.g,h) is the product of 
this
solution with a group homomorphism $\pi \to \{k\in K^*\,\vert\, k^2=1\}$.  }

\skipaline {\sl Proof.}  It is obvious that any two solutions to (2.6.g) (with
given $c$) are obtained from each other by multiplication by a group 
homomorphism
$\pi \to K^*$.  The condition $\theta_{{\alpha}^{-1}}=\theta_{\alpha}$ implies
that this homomorphism takes values in the group $\{k\in K\,\vert\,
k^2=1\}\subset K^*$.

Let us prove the first claim.  A direct computation using first (2.6.e)
and then (2.6.f), (2.6.a) yields $$c_{\alpha \beta, \alpha \beta} =
c_{\alpha,\alpha} \,c_{\beta,\beta} \,c_{\alpha,\beta} \,c_{\beta,\alpha}\times
$$ $$\times [ ( a_{\alpha \beta, \alpha, \beta}\, a_{\alpha, \beta, \alpha 
\beta}
)^{-1} \,a_{\alpha \beta \alpha \beta^{-1}\alpha^{-1}, \alpha \beta, \beta} \,
a_{ \beta \alpha \beta^{-1}, \beta, \alpha }\, a_{\alpha \beta \alpha^{-1},
\alpha , \beta}].$$ Formula (2.6.g) follows from the fact that the expression 
in
the square brackets is equal to 1.  This can be obtained from (1.3.a) by 
applying
substitution $\alpha\mapsto \beta \alpha \beta^{-1}, \gamma\mapsto \alpha,
\delta\mapsto \alpha^{-1}\beta \alpha$ and using (2.6.a).

   It remains to check the identity $c_{\alpha^{-1},\alpha^{-1}}
=c_{\alpha,\alpha}$.  Applying (2.6.e) to $\beta= \gamma=1$ we obtain
$c_{\alpha,1}=a_{\alpha,1,1}\, a_{1,1,\alpha}$.  Applying (2.6.e) to
$\beta=\alpha, \gamma=\alpha^{-1}$ we obtain $$c_{\alpha,\alpha^{-1}}=
(c_{\alpha,\alpha})^{-1} \,c_{\alpha,1}\, a_{\alpha,\alpha^{-1},\alpha}=
(c_{\alpha,\alpha})^{-1} \,a_{\alpha,1,1}\,
a_{1,1,\alpha}\,a_{\alpha,\alpha^{-1},\alpha} .$$ Similarly, substituting
$\beta=\delta=1, \gamma=\alpha$ and $\beta=\alpha, \gamma=\alpha^{-1},
\delta=\alpha^{-1}$ in (2.6.f) we obtain $c_{1,\alpha}=(a_{\alpha,1,1} \,
a_{1,1,\alpha})^{-1}$ and $$c_{\alpha^{-1},\alpha^{-1}}= c_{ 1,\alpha^{-1}}\,
(c_{\alpha,\alpha^{-1}}\,
a_{\alpha^{-1},\alpha,\alpha^{-1}})^{-1}=(c_{\alpha,\alpha^{-1}} \,
a_{\alpha^{-1},1,1}\,a_{1,1,\alpha^{-1}}\,a_{\alpha^{-1},\alpha,\alpha^{-1}})^{
-1}.$$ Substituting here the expressions for $c_{ \alpha,\alpha^{-1}} $ 
obtained
above and using the fact that the righ-hand sides of (1.3.b) and (1.3.c) are
equal we obtain $c_{\alpha^{-1},\alpha^{-1}} =c_{\alpha,\alpha}$.

\skipaline \centerline {\bf   Appendix 2. State sum invariants of 3-dimensional
$\pi$-manifolds}

\skipaline \noindent {\bf 1. Spherical categories.} We begin with  a definition
of a pivotal monoidal category, cf.  [FY], [BW2].  Let $\Cal C$ be a $K$-additive monoidal category
with (left)  duality $b,d$ and structural morphisms $a,l,r$, as in Section 1. 
Using these morphisms one can define  for any
objects $V,W\in \Cal C$   a canonical isomorphism  $\gamma_{V,W}:
 V^*\otimes W^*\to (W\otimes V)^*$   (see, for instance, [Tu2, p.31]).
In particular, if $\Cal C$ is strict then
$$\gamma_{V,W}=
(d_V\otimes \text {id}_{{(W\otimes V)^*}})(\text {id}_{  V^*}\otimes
d_W\otimes \text {id}_V\otimes \text {id}_{{(W\otimes V)^*}} )(\text {id}_{ 
V^*}\otimes \text {id}_{W^*}\otimes b_{{W\otimes V}}),$$
$$\gamma_{V,W}^{-1}=
(d_{{W\otimes V}}\otimes  \text {id}_{V^*}\otimes  \text {id}_ {W^*})(\text
{id}_{{(W\otimes V)^*}}\otimes  \text {id}_W \otimes  b_V\otimes  \text
{id}_{W^*}) (\text {id}_{{(W\otimes V)^*}}\otimes  b_W) .$$
Recall also  that the duality defines  a contravariant  functor  $*:\Cal C \to \Cal
C$ (see Section 2.4).

The data
needed to make   $\Cal C$ {\it pivotal} is a system of invertible morphisms
$\{\tau_V:V\to V^{**}\}_{V\in \Cal C}$  
satisfying the following four conditions:

(1) $\tau$ is a natural transformation  $\id_{\Cal C} 
\to **$, i.e., for any morphism $f:V\to W$ in $\Cal C$, the following diagram is
commutative:
$$\CD V @>f>>   W\\  @V{\tau_V}VV    @VV{\tau_W}V\\
  V^{**} @> f^{**} >> W^{**} ; \endCD  $$

(2) for  any objects $V,W\in \Cal C$, 
the following diagram is
commutative:
$$\CD V\otimes W @>\tau_{V\otimes W} >>   ({V\otimes W})^{**}\\  @V{\tau_V
\otimes {\tau_W}}VV    @VV{\gamma_{W,V}^*}V\\
  V^{**} \otimes W^{**}@> \gamma_{V^*,W^*} >> (W^* \otimes V^*)^* ; \endCD  $$

(3) for any object  $V\in  \Cal C$, we have $\tau_{V^*}= (\tau^{-1}_V)^*:V^*\to
V^{***}$;

(4) the morphism $\nu= 
(r_{\1^*})^{-1}\,b_{\1}:\1\to \1^{*}$ is invertible and 
 $\tau_{\1}=(\nu^{-1})^* \nu$.

Let $\Cal C$ be a pivotal category. For a morphism $f:V\to V$  in $\Cal C$,  the
trace  $\tr (f)\in \End_{\Cal C} (\1)$ is defined by $$\tr(f)=d_{V^*}   (\tau_V
f\otimes \id_{V^*}) \,b_V:\1\to \1.$$ The trace verifies 
$\tr(fg)=\tr(gf)$ for any morphisms $f:V\to W, g:W\to V$. 

Following [BW2], we call a pivotal category $\Cal C$   {\it spherical}   if
$\tr(f)=\tr(f^*)$ for any endomorphism $f$ of an object  in $\Cal C$. 
In a spherical category we always have $\tr(f\otimes g)=\tr(f) \,\tr(g)$ where on the
right-hand side we use multiplication in  $\End_{\Cal C}(\1)$ defined by composition.  This
multiplication makes $\End_{\Cal C}(\1)$ a commutative $K$-algebra with unit $\id_{\1}$. 
In a spherical category, $\tr(\id_{\1})=\id_{\1}$. 

Barrett  and 
 Westbury [BW1] (see also [GK]) extended the state sum approach to 3-manifold invariants introduced in [TV]  
to spherical categories satisfying a few extra conditions. Namely, they showed
that the state sum approach works  for  any  finite semisimple  
spherical category $\Cal C$ over a field $K$ such that $\End (\1_{\Cal C})=K$. 
Thus, such a category  gives rise to a $K$-valued
topological invariant  of  closed  oriented 3-manifolds.

\skipaline \noindent {\bf 2. Remarks. } 1. A pivotal category $\Cal C$ has  a right duality
$ \tilde b, \tilde d$
defined   by
$$\tilde b_V= (\id_{V^*} \otimes \tau_V^{-1})\, b_{V^*}:\1\to V^*\otimes V,\,\,\,
\tilde d_V= d_{V^*} (\tau_V\otimes \id_{V^*} ):V\otimes V^*\to \1.$$
One can  check that     the 
contravariant functor    defined as in Section 2.4  but using $ \tilde b, \tilde d$ coincides with the duality functor  ${}^*$ 
defined in  Section 2.4 using $b,d$. Similarly, the  natural isomorphisms $ 
 V^*\otimes W^*\to (W\otimes V)^*$ defined using the right (resp. left) duality
coincide. Note also the   equality $\nu=r_{\1^*}^{-1} b_{\1}=l_{\1^*}^{-1}
\tilde b_{\1}:\1\to \1^*$.

2. Our definition of a pivotal category    differs from but is  
essentially equivalent to the one in [BW2].  Barrett  and 
 Westbury  define a pivotal category as a
monoidal category  $\Cal C$  equipped with  a contravariant functor $*:\Cal C\to
\Cal C$,  an isomorphism $\nu:\1\to \1^*$,    morphisms $\{b_V:\1\to V\otimes
V^*,  \tau_V:V\to V^{**}\}_{V\in \Cal C}$, and invertible morphisms 
$\{\gamma_{V,W}:
 V^*\otimes W^*\to (W\otimes V)^*\}_{V,W\in \Cal C}$  satisfying certain conditions.
Among their conditions are
  axioms (1) - (3) above and the equality $\tau_{\1}= (\nu^{-1})^* \nu$ replacing
axiom (4) above. Barrett  and 
 Westbury do not suppose the existence of   morphisms $d $
forming  together with $b $ a duality in $\Cal C$. 
However,  this  can be   deduced from their  axioms.  We prefer to assume  from the very beginning that $\Cal C$ has a
duality and extract from it the duality functor  $\ast$ and  the isomorphisms $\gamma$
and   $\nu $.
Conversely, our axioms   imply all the Barrett-Westbury axioms.    Only Condition 
(4)  in [BW2, p.\ 362] is somewhat involved; the proof   uses 
  the   right duality $ \tilde b, \tilde d$  defined above. 

\skipaline \noindent {\bf 3. Spherical crossed $\pi$-categories. }  We say that a crossed $\pi$-category 
$\Cal C$ is {\it spherical} if it is
spherical in the sense of the definition above and the given system of morphisms
$\{\tau_V:V\to V^{**}\}_{V\in \Cal C}$ is invariant under the action of $\pi$ on  $\Cal C$.  A crossed $\pi$-category  $\Cal
C$ is {\it finite semisimple} if it satisfies axioms  (6.1.1) - (6.1.4).   The methods of [TV], [BW1] allow to derive from a
finite semisimple spherical  crossed $\pi$-category $\Cal C$ a topological  invariant $\vert M \vert_{\Cal C}\in K$ of a closed
oriented 3-dimensional
$\pi$-manifold $M$. The construction goes as follows.  A {\it local order}  on a triangulation $T$  of $M$  is a compatible  choice
for each simplex
$\Delta \in T$ of a total order on the set of    vertices of $\Delta$. By  \lq\lq compatible" we mean  that the order on the
vertices of any  subsimplex of any $\Delta\in T$
  is induced by the order on the vertices of $\Delta$. For instance, a  total order on the set of all vertices of
$T$ induces a local order on $T$.  Note that a local order on $T$ determines an orientation of all edges of $T$.
Fix now  a triangulation
$T$ of
$M$ and endow it with a local  order.   Choose in the given homotopy class of maps
$M\to K(\pi, 1)$ a representative,
$g$,  sending all the vertices to the base point of $K(\pi,1)$.   We  assign to each   edge $e$ of    $T$ the 
element  
$g_e\in \pi$ represented by the loop $g\vert_e$. In
terminology of [Tu3, Sect. 7.2],  the function
$e\mapsto g_e, e^{-1} \mapsto g_e^{-1}$
 is a $\pi$-system   on $T$  representing the homotopy class of
$g$.
For $\alpha\in \pi$, denote by
$I_\alpha$ the set   of  the  isomorphism  classes of
simple objects of ${\Cal C}_\alpha$. 
 We define $\vert T \vert_{\Cal C}$     as in [BW1] with the difference that we involve only the labelling 
assigning to each   edge $e$ of $T$ an element of   $I_{g(e)}$. The rest of the construction is similar to the one in 
[BW1].  

 We claim that  $\vert M \vert_{\Cal C}=\vert T \vert_{\Cal C}$ is  a well-defined topological invariant of $M$, i.e., it  does not
depend on the choice of  
$T$, the choice of the local order on $T$,  and the choice  of $g$ in the given homotopy class.
We outline the   proof. 
Note first that a Pachner move on $T$ extends in a unique way to the $\pi$-systems. The extension is uniquely
determined by the condition that the values of the
$\pi$-system on all edges preserved by the move should be  preserved.
The  same arguments  as in [BW1] show that   $\vert T \vert_{\Cal C}$ does not depend of the choice of the local order
on $T$ and is invariant under   the  Pachner moves.  It remains to prove  that 
$\vert T
\vert_{\Cal C}$ does not depend on the choice of  $g$ in its homotopy class.
  Any two  $\pi$-systems  on $T$ representing homotopic maps $M\to K(\pi,1)$
can be related by homotopy moves at the vertices of  $T$, cf. [Tu3, Sect. 7.2]. 
It would be nice to give a direct proof of the invariance of $\vert T \vert_{\Cal C}$ under these moves. 
The proof outlined below is based on the theory of skeletons of $M$. 
A {\it skeleton}  of $M$ is a finite simple
2-polyhedron in
$M$ whose complement is a disjoint union of  open 3-balls. For instance, the closed 2-cells in $M$ dual to the edges of 
$T$ form a skeleton $T^*\subset M$.  Using this observation,  we can   dualize the
notion of a $\pi$-system  from triangulations to   skeletons of $M$.    (A $\pi$-system on a skeleton assigns 
elements of
$\pi$ to   oriented 2-faces). If $T$ is locally ordered then  
we can provide
 the 2-cell dual to  an edge $e$  with an orientation   which 
determines together with the orientation of $e$   the given orientation of $M$. 
This makes    $T^*$ an oriented branched 2-polyhedron in the sense of [BP]. Now we observe that the state sum
invariant can be defined for a $\pi$-system on an oriented branched skeleton of $M$.  This allows us to switch from the
language of  state sums on triangulations  to the language of state sums on oriented branched skeletons.
In the latter language, a homotopy move on  $\pi$-systems can be decomposed into a composition of more elementary
Matveev-Piergallini moves and bubblings, cf. [TV] and 
 [Tu3, Sect. 8].  The invariance of the state sum under these moves is verified as in [BW1]. 

\skipaline \noindent {\bf 4. Spherical algebras. }
The notion of a spherical category has its counterpart in the theory of Hopf algebras. Following [BW2] we call a {\it spherical
Hopf algebra} over
$K$ any pair $(A,w)$ where $A$ is a Hopf algebra over $K$ and $w\in A$ is a group-like element   such that
 the square of the antipode in $A$ equals the conjugation by $w$ and for any $A$-linear endomorphism  $f$ of an
$A$-module (in the sense of Section 11.7) we have $\Tr(fw)=\Tr (fw^{-1})$
where $\Tr$ denotes the
  trace of a $K$-endomorphism of a projective $K$-module.  
 We can extend this definition  to crossed Hopf
$\pi$-coalgebras.  A {\it spherical crossed Hopf
$\pi$-coalgebra} is a pair  consisting of a crossed Hopf
$\pi$-coalgebra $(\{A_{\alpha}\}_{\alpha\in \pi},
\Delta , \varepsilon_1, s,\varphi)$ and invertible  elements $\{w_\alpha\in A_\alpha\}_{\alpha\in \pi}$
satisfying the following conditions:
$$ s_{\alpha^{-1}} s_\alpha (a) =w_\alpha a w_{\alpha}^{-1},\,\,\, \Delta_{\alpha, \beta} (w_{\alpha\beta})=w_\alpha
\otimes w_\beta, \,\,\,s_\alpha (w_\alpha)=w_{\alpha^{-1}}^{-1},
$$
$$\varepsilon_1 (w_1)=1,\,\,\,
\varphi_\alpha (w_\beta)=w_{\alpha \beta \alpha^{-1}}, \,\,\, \Tr(fw_\alpha)=\Tr (fw_{\alpha}^{-1})$$
for any $\alpha, \beta \in \pi, a\in A_\alpha$
  and any $A_\alpha$-linear endomorphism  $f$ of an
$A_\alpha$-module (in the sense of Section 11.7).

  It is clear that an involutory
crossed Hopf
$\pi$-coalgebra (so that $ s_{\alpha^{-1}} s_\alpha=\id$ for all $\alpha\in \pi$) is   spherical with $w_\alpha=1_\alpha\in
A_\alpha$ for all $\alpha\in \pi$.  A ribbon crossed Hopf
$\pi$-coalgebra $(A,\theta)$ is sperical with $w_\alpha=\theta_\alpha u_\alpha$ where
$u_\alpha\in A_\alpha$ is the (generalized) Drinfeld element of $R$ (see [Vi]). 

It is proven in [BW2] that 
for  a spherical Hopf algebra $(A,w)$, the monoidal category $\Rep (A) $ is spherical.
For $V\in \Rep (A) $, the morphism $\tau_V:V\to V^{**}$ is defined as the standard    identification $V= V^{**}$
followed by multiplication by $w^{-1}$. 
Similarly, for a spherical crossed Hopf
$\pi$-coalgebra $(A, w)$, the category $\Rep (A) $ is spherical: 
for $V\in \Rep (A_\alpha) $, the morphism $\tau_V:V\to V^{**}$ is defined as the standard identification $V= V^{**}$
followed by multiplication by $w_{\alpha}^{-1}$.

\skipaline \centerline {\bf   Appendix 3. Open problems}

\skipaline \noindent {\it 1. Classification of crossed $\pi$-algebras.} Classify
crossed $\pi$-algebras for sufficiently simple groups $\pi$, say, cyclic, abelian, finite, etc. 
Find interesting  examples of non-semisimple crossed $\pi$-algebras.

\skipaline \noindent {\it 2. Relations between various approaches.}  In the  theory of quantum invariants
of 3-dimensional  manifolds,  there
is a narrow relationship between the surgery approach and the state sum approach, see [Tu2], [Ro].
It would be interesting to generalize it to  HQFT's. A related question: generalize the state sums considered
in Appendix 2 to shadows, cf.  [Tu2].

\skipaline \noindent {\it 3. Modular  $\pi$-categories.}   Can one systematically produce  interesting modular $\pi$-categories ? 
By \lq\lq interesting" I mean those yielding
 non-trivial topological invariants allowing to distinguish 3-dimensional $\pi$-manifolds.
In particular,  if $\pi$ is a subgroup of a semisimple Lie group $G$ with Lie algebra
$\,\bold     g$ and $q$ is a  complex root of unity, can one use the quantum group $U_q({\,\bold     g})$ to construct modular 
crossed Hopf
$\pi$-algebras ?
 
\skipaline \noindent {\it 4. Invariants of spin-structures.} Instead of maps from  a  manifold $M$ 
   to a fixed target space one can consider maps   whose source is  the oriented   frame bundle of $M$ (the principal bundle of
positively oriented bases in the tangent spaces of points). When the  target space is $K(\pi,1)$, this should lead to  new 
algebraic notions  generalizing  $\pi$-algebras and crossed $\pi$-categories.    This theory should include the   quantum
invariants of spin-structures on closed oriented 3-manifolds, see  [Bl], [KM], [Tu1].

\skipaline \noindent {\it  5. Factors and subfactors.}
It was established by A. Ocneanu that the subfactors of type $II_1$ give rise to 3-manifold invariants via the state sum approach
(see for instance [EK]).
What is the counterpart of the theory introduced above in the setting of 
  subfactors ? 

\skipaline \noindent {\it  6. Perturbative aspects.}
What are the perturbative aspects  of HQFT's ? Are there perturbative  
invariants of 3-dimensional $\pi$-manifolds generalizing the Le-Murakami-Ohtsuki invariant [LMO] ?

\skipaline \noindent {\it  7. Higher-dimensional generalizations.}  
The  quantum invariants of 3-manifolds have  4-dimensional counterparts, see  [CKY],
[CKS], [CKJLS], [Mac]  and references therein.  The role of categories in these constructions is typically
played by 2-categories. The homotopy quantum field theory should have similar  4-dimensional and
high-dimensional versions and give rise in particular to the notion of a crossed 2-category.

In this paper we considered mainly   target spaces of type $K(\pi,1)$.
In dimension 4 it can be easier and more appropriate  to consider target spaces of type $K(H,2)$ where $H$ is an abelian group.
Note that 2-dimensional HQFT's  with  such target space  were studied in  [BT].

It is most interesting to  analyze the algebraic data yielding invariants of $\Spin^c$-structures on 4-manifolds.
Here as the  sources of maps one takes the  oriented frame bundles of 4-manifolds and  as the target space 
one takes   $K(H,2)$.
Can one include the Seiberg-Witten invariants in this framework ?

\skipaline \noindent {\it  8. Miscellaneous questions.}  Study the representations of   subgroups of the
mapping class groups resulting from HQFT's. Find an interpretation of HQFT's in terms of algebraic
geometry, i.e., in terms of  sections of  bundles over moduli spaces. Study   the
2-dimensional homotopy modular functors and their relations to quantum computations (cf.
[FKW]).  Study relations with number theory, cf. [LZ].

\skipaline

\skipaline \centerline {\bf References}

\skipaline

[BW1] Barrett, J., Westbury, B., Invariants of piecewise-linear 3-manifolds.
Trans. Amer. Math. Soc.  348 (1996), 3997--4022.

[BW2] Barrett, J., Westbury, B.,  Spherical categories. Adv. Math.  143 (1999),
357--375. 

[BP] Benedetti, R.,  Petronio, C., Branched Standard Spines of $3$-manifolds. Lect. Notes in Math., 1653.
Springer-Verlag, Berlin, 1997. 

[Bl]  Blanchet, C., Invariants of
three-manifolds with spin structure. Comm. Math. Helv.  67 (1992),
406--427.

[Br]  Brieskorn, E., Automorphic sets and singularities. Contemp. Math. 78 (1988), 
45--115.

[BT]  Brightwell, M. Turner, P., 
Representations of the homotopy surface category of a simply connected space. Preprint 
math. AT/9910026.

[Bru] Brugui\`eres, A., Cat\'egories pr\'emodulaires, modularisations et invariants
des vari\'et\'es de dimension 3.  Math. Ann. 316
(2000), 215--236. 

[CKS]  Carter, J.,  Kauffman, L., Saito, M., Structures and Diagrammatics of Four 
Dimensional Topological Lattice Field Theories. 
 Adv. Math. 146 (1999), no. 1, 39--100.

[CKJLS] Carter, J.,  Kamada, S., Jelsovsky, D., Langford, L., Saito, M., 
State-sum invariants of knotted curves and surfaces from quandle cohomology. Electron. Res. Announc. Amer. Math.
Soc. 5 (1999), 146--156 (electronic). 

[CKY] Crane, L.,  Kauffman, L., Yetter, D. State-sum invariants of 4-manifolds. J. Knot Theory Ramifications 6 (1997), 177--234.

[EK] ]   Evans, D.,  Kawahigashi, Y.,  Quantum symmetries on operator algebras. Oxford Mathematical Monographs. Oxford Science Publications.
The Clarendon Press, Oxford University Press, New York, 1998.

[FR1]  Fenn, R.,  Rourke, C., On Kirby's calculus of links.
Topology 18 (1979), 1--15.

[FR2] Fenn, R., Rourke, C., Racks and Links in Codimension Two. J. of Knot Theory and
Ramifications 1 (1992), 343--406.

[FKW]  Freedman, M.,   Kitaev, A.,   Wang, Z., 
Simulation of topological field theories by quantum computers. Preprint quant-ph/0001071.

[FY] Freyd, P., Yetter, D., Braided compact closed categories
 with applications
to low-dimensional topology. Adv. Math.  77 (1989), 156--182.

[GK]  Gelfand, S., Kazhdan, D. Invariants of 3-dimensional manifolds. Geom.
Funct. Anal. 6 (1996),  268--300.

[He]  Hennings, M., Invariants of links and $3$-manifolds obtained from Hopf algebras. 
J. London Math. Soc. (2) 54 (1996), 594--624. 

[KRT]    Kassel, C.,   Rosso, M. , Turaev, V.,  Quantum Groups  and Knot Invariants,
Panoramas et Synth\`eses 5, Soci\'et\'e Math\'ematique de France, 115 p., 1997.

[Ki]  Kirby,  R., A calculus of framed links in $S^3$. Invent. Math.
45 (1978), 35--56.

[KM] Kirby, R.,  Melvin, P., On the 3-manifold invariants of Witten and
Resheti\-khin-Turaev for $sl_2({\bold C})$. Invent. Math. 105 (1991),
473--545.

[Ku] Kuperberg, G.,  Involutory Hopf algebras and $3$-manifold invariants. Internat. J. Math. 2 (1991),  
41--66. 

[LZ] Lawrence, R.,  Zagier, D.,  Modular forms and quantum invariants of $3$-manifolds.   Asian J. Math. 3
 (1999),  93--107.

[LMO] Le, T., Murakami, J.,  Ohtsuki, T.,  On a universal perturbative invariant of $3$-manifolds. Topology 37 (1998),  539--574.

[LT] Le, T., Turaev, V.,  Quantum groups and abelian HQFT's, in preparation.

[Mac] Mackaay, M.,  Spherical $2$-categories and $4$-manifold invariants. 
Adv. Math. 143 (1999),   288--348. 

[Ma] MacLane, S., Categories for the working mathematician.  Graduate Texts in
Math.  5.  Springer-Verlag, New York-Heidelberg-Berlin 1971.

[Oh1]  Ohtsuki, T.,  Colored ribbon Hopf algebras and universal invariants of framed
links. J. of Knot Theory and Ramifications 2 (1993),   211--232. 

[Oh2] Ohtsuki, T., 
 Invariants of $3$-manifolds derived from universal invariants of framed
links. Math. Proc. Cambridge Philos. Soc. 117 (1995),   259--273.

[Os] Ospel, C., Tressages et th\'eories cohomologiques pour les alg\`ebres de Hopf.
Application aux invariants des 3-vari\'et\'ets. Ph. D. thesis (Strasbourg, 1999).

[RT]   Reshetikhin, N., Turaev, V., Invariants of 3-manifolds via
link polynomials and quantum groups. Invent. Math. 103 (1991), 547--598.

[Ro] Roberts, J., Skein theory and Turaev-Viro invariants. Topology 34 (1995),  771--787.

[Th] Thys, H,. Description topologique des repr\'esentations de $U_q(sl_2)$, Preprint
math. QA/9810088, to
appear in Ann. Fac. Sci. Toulouse.

[Tu1] Turaev, V., State sum models in low-dimensional topology.   Proc.
ICM, Kyoto, 1990  vol. 1 (1991), 689--698.

[Tu2] Turaev, V., Quantum invariants of knots and 3-manifolds.  Studies in
Mathematics 18, Walter de Gruyter, 588 p., 1994.

[Tu3] Turaev, V.,  Homotopy  field  theory in dimension 2 and group-algebras.
Preprint QA/9910010.

[TV]  Turaev, V.,  Viro, O., State sum invariants of 3-manifolds
and quantum $6j$-symbols.   Topology 31 (1992), 865--902.

[Vi] Virelizier, A., Crossed Hopf $\pi$-coalgebras and invariants of links and
3-manifolds, in preparation.

[Zu] Zunino, M.  Center of $\pi$-categories and double of 
crossed Hopf $\pi$-coalgebras, in preparation.

\skipaline

 {Institut de Recherche Math\'ematique Avanc\'ee, Universit\'e
Louis Pasteur - CNRS,
7 rue Ren\'e Descartes, 67084 Strasbourg Cedex, France} 

\skipaline

turaev$\@$math.u-strasbg.fr

\end

[At] Atiyah, M., Topological quantum field theories.  Publ.  Math.  IHES 68
(1989), 175-186.

\skipaline \noindent {\bf 8.1.      Homotopy quantum field theory.}
The notion of a  homotopy quantum field theory (HQFT)   introduced in [Tu2] is a 
version of  the notion of a topological    quantum field theory (TQFT). A
$(d+1)$-dimensional  HQFT   may be   described as a TQFT   for  closed 
$d$-dimensional manifolds and $(d+1)$-dimensional cobordisms endowed with
homotopy classes of maps into a given space. Throughout this paper  $d=2$ and the
target space is the Eilenberg-MacLane space $X=K(\pi,1)$ associated with the group
$\pi$. 

We fix an  Eilenberg-MacLane space $X=K(\pi,1)$ associated with   $\pi$
and a base point $x\in X$. We shall always assume that $X$ is a CW-space and $x$
is a 0-cell of $X$.  We fix a modular crossed $\pi$-category $\Cal C$. In this
section we shall    derive  from $\Cal C$  a 3-dimensional HQFT with target   
$X=K(\pi,1)$.  

We shall use the language of   pointed homotopy theory. A
topological space is {\it pointed} if all its connected components are provided
with base points.  
A map  between pointed spaces   is a continuous map  sending base points
into base points and considered up to homotopy constant on
the base points.  Note that maps from a pointed connected CW-complex $Y$ with
base point $y$ into  $X$ bijectively correspond to group homomorphisms
$\pi_1(Y,y)\to \pi=\pi_1(X,x)$.

Let $(M, \partial_- M, \partial_+ M)$  be an extended    $\pi$-cobordism and let
$\Upsilon\subset M\backslash \partial M$ be a closed surface splitting $M$ into two
compact submanifolds $M_1,M_2$ such that
$\partial (M_1)=\partial_- M \cup \U$ and 
$\partial (M_2)=\partial_+  M \cup \U$

Let $(M_1, \partial_-
(M_1), \partial_+ (M_1))$ and $(M_2, \partial_-
(M_2), \partial_+ (M_2))$ be two extended    $\pi$-cobordisms 
such that the surfaces $\partial_+ (M_1)$ and  $\partial_-
(M_2)$ are non-void and 
  $e$-homeomorphic. We can glue   these cobordisms along an
$e$-homeomorphism $p:\partial_+ (M_1)\to   \partial_-
(M_2)$ to obtain  an  extended    $\pi$-cobordism $M$.
 The cobordism $M$ is obtained from  $M_1\amalg
M_2$ by identification of $\partial_+ (M_1)$ with $   \partial_- (M_2)$ along $p$.
The bases    
$\partial_-
(M_1) $ and $  \partial_+ (M_2)$ of $M$ keep their structure of an extended
$\pi$-surface determined by $M_1,M_2$.

 For $r=1,2$, let $(\Omega_r\subset M_r, g_r:M_r\backslash \Omega_r\to X,
u_r,v_r)$ be  the given colored $\pi$-graph in $M_r$. We  form a ribbon graph
$\Omega=\Omega_1\cup \Omega_2 \subset M$ and provide it
with the structure of a colored $\pi$-graph.  
Clearly, it suffices to consider the case where $M$ is connected and we assume it in
the sequel.

 Since $p$ is an $e$-homeomorphism, we
can choose representatives of the homotopy classes of maps $
g_r:M_r\backslash \Omega_r\to X$  (where $r=1,2$) so that on the
complement of marked points in $\partial_+ (M_1)$ we have $g_1 =g_2 p$. Then we
  glue these representatives into a map  $g:M\backslash \Omega \to X$.
(This map   may depend  on the choice of the
representatives, but using the facts that $X=K(\pi,1)$ 
and that we operate in the pointed category one can show that the  homotopy class
of $g$ is well-defined. We however do not need this.)
 Now, we need to consider two cases. 

Case $  \partial_-
(M_1)=  \partial_+ (M_2)=\emptyset$. Then $M$ is a closed connected 3-manifold.
As a base point   $z\in M\backslash \Omega$ we choose the base point of any
component of $\partial_+ (M_1)=\partial_-
(M_2)\subset M$. By assumptions $g(z)=x\in X$. 
The triple $(\Omega\subset M, z\in M\backslash \Omega, g_\#:\pi_1(M\backslash
\Omega)\to \pi)$ is then a $\pi$-graph in $M$ in the sense of Section  7.5. One
can check that there is a unique coloring $(u,v)$ of $\Omega$ such that for
both $r=1,2$ and any path $\gamma$ in $M_r\backslash \Omega_r$ leading from $z$
to a point on $\tilde {\Omega}_r\subset \tilde \Omega$
the object $u_\gamma$ (or the mporphism $v_\gamma$) 
the given coloring of $\Omega_r\subset M_r$ determines a object  $u_\gamma$
associated to . This function extends uniquely
to a coloring, $u$,  of $\Omega$. 

(10.2.5) Let $(M, \partial_-
M, \partial_+ M)$ be an  extended  (3-dimensional) $\pi$-cobordism
obtained from  extended    $\pi$-cobordisms 
$(M_1, \partial_-
(M_1), \partial_+ (M_1))$ and $(M_2, \partial_-
(M_2), \partial_+ (M_2))$ by gluing along a weak $e$-homeomorphism
$p:\partial_+ (M_1)\to  \partial_-
(M_2)$.